\documentclass[11pt,reqno]{amsart}

\usepackage[
    a4paper,
    margin=1in
]{geometry}
\usepackage{setspace}
\usepackage{needspace}
\usepackage{float}

\usepackage{appendix}

\usepackage{amssymb}
\usepackage{mathtools}

\numberwithin{equation}{section}

\usepackage[
    setpagesize=false,
    colorlinks=true,
    linkcolor=BrickRed,
        citecolor=OliveGreen,
        urlcolor=black,
    pdfencoding=auto,
    psdextra,
]{hyperref}
\usepackage{cleveref}

\usepackage{etoolbox} 
\usepackage[shortlabels]{enumitem}
\usepackage{xparse} 
\usepackage{blkarray,multirow}

\usepackage{graphicx}
\usepackage[dvipsnames]{xcolor}
\usepackage{tikz}
\usetikzlibrary{patterns}
\usetikzlibrary{decorations.markings, arrows.meta, calc, positioning}

\newcommand{\matchingfigurestyles}{%
\definecolor{MBlue}{RGB}{47,79,107}
\definecolor{MRed}{RGB}{157,57,44}
\definecolor{MOlive}{RGB}{102,114,54}
\definecolor{MGray}{RGB}{117,122,127}
\tikzset{
  vertex/.style={circle,fill=black,inner sep=1.55pt,outer sep=0pt},
  vacant/.style={circle,draw=black!65,fill=white,line width=.6pt,inner sep=2pt},
  source/.style={draw=MBlue,fill=MBlue!5,line width=.75pt,rounded corners=5pt},
  alternative/.style={draw=MRed,fill=MRed!10,line width=1pt,dash pattern=on 3pt off 2pt,rounded corners=4pt},
  arr/.style={-{Stealth[length=5pt,width=5pt]},line width=.8pt},
  panel/.style={font=\small\bfseries,align=center},
  note/.style={font=\small,align=center},
  every node/.style={font=\small},
}
}

\usepackage{todonotes}
\IfFileExists{ulem.sty}{\usepackage[
    deletedmarkup=sout,
    commentmarkup=uwave,
    authormarkuptext=name
]{changes}}{}

\usepackage{comment}

\makeatletter
\renewcommand{\section}{\@startsection{section}{1}%
    \z@{1.3\linespacing\@plus.5\linespacing\@minus.2\linespacing}%
    {.6\linespacing}%
    {\normalfont\scshape\centering}}
\renewcommand{\paragraph}{\@startsection{paragraph}{4}%
    \z@{.75\linespacing\@plus.2\linespacing\@minus.1\linespacing}%
    {-.6em}%
    {\normalfont\bfseries}}
\newcommand{\proofparagraph}[1]{%
    \paragraph*{$\bullet$ #1}\ignorespaces}
\makeatother

\setlist[enumerate,1]{label=(\arabic*), ref=(\arabic*)}
\setlist[enumerate,3]{label=(\roman*), ref=(\roman*)}

\theoremstyle{plain}
\newtheorem{theorem}{Theorem}[section]
\newtheorem{lemma}[theorem]{Lemma}
\newtheorem{corollary}[theorem]{Corollary}
\newtheorem{proposition}[theorem]{Proposition}

\newtheorem{fact}[theorem]{Fact}
\newtheorem{conjecture}[theorem]{Conjecture}

\newtheorem{question}[theorem]{Question}

\newtheorem{claim}{Claim}[section]
\newtheorem*{claim*}{Claim}

\makeatletter
\newenvironment{claimproof}[1][Proof]{\par
    \pushQED{\qed}%
    
    \normalfont \topsep6\p@\@plus6\p@\relax
    \trivlist
    \item[\hskip\labelsep
    \textit{#1}\@addpunct{.}~]\ignorespaces
}{%
    \popQED\endtrivlist\@endpefalse
}
\makeatother

\theoremstyle{definition}
\newtheorem{definition}[theorem]{Definition}
\newtheorem{remark}[theorem]{Remark}

\newcommand{\calH}{\mathcal{H}}
\newcommand{\calM}{\mathcal{M}}
\newcommand{\calP}{\mathcal{P}}
\newcommand{\calQ}{\mathcal{Q}}

\newcommand{\calS}{\mathcal{S}}
\newcommand{\calT}{\mathcal{T}}

\title[Matching statistics in regular hypergraphs and applications]{Matching statistics in regular hypergraphs and applications}

\author{Hyunwoo Lee}
\address{Department of Mathematical Sciences, KAIST, South Korea and Extremal Combinatorics and Probability Group (ECOPRO), Institute for Basic Science (IBS)}
\email{hyunwoo.lee@kaist.ac.kr}
\thanks{Supported by the National Research Foundation of Korea (NRF) grant funded by the Korea government (MSIT) No. RS-2023-00210430, and the Institute for Basic Science (IBS-R029-C4).}
\date{\today}

\usepackage[square,sort,comma,numbers]{natbib}
\usepackage{booktabs}
\allowdisplaybreaks[2]
\newcommand{\E}{\mathbb E}
\newcommand{\Prob}{\mathbb P}
\newcommand{\Var}{\operatorname{Var}}
\newcommand{\HH}{\mathsf H}
\newcommand{\eps}{\varepsilon}
\newcommand{\calU}{\mathcal U}
\newcommand{\dd}{\,\mathrm d}
\newcommand{\qbinom}[2]{\genfrac{[}{]}{0pt}{}{#1}{#2}_q}

\begin{document}

\begin{abstract}
    For fixed $k\ge3$, we determine the first sublinear correction to the logarithm of the matching partition function for a class of $n$-vertex $d$-regular $k$-uniform hypergraphs, uniformly for activities in each fixed compact interval of $(0,\infty)$.
    Our hypotheses include every regular linear hypergraph of polynomial degree and a class of nonlinear hypergraphs described by higher codegrees and intersection patterns, which we introduce here.
    The weighted counting estimate implies an averaged version of Kahn's conjecture: under the uniform measure, the expected number of uncovered vertices is $(1+o(1))nd^{-1/k}$, with concentration and large deviations on this scale.
    The proof combines a novel variant of random greedy counting, based on short stages and regular completions, with a new entropy method that uses intersection patterns to control reconstruction after deleting matching edges.

    We highlight the following applications.
    \begin{enumerate}[(i)]
        \item \emph{Partial transversals in Latin squares.}
            Every Latin square of order $n$ has $(n/e^2)^n\exp((3+o(1))n^{2/3})$ partial transversals. Also, almost all partial transversals have $(1+o(1))n^{2/3}$ missing cells.
        \item \emph{Linear $3$-uniform hypergraphs.}
            There are $(n/e^2)^{n^2/6}\exp((1/2+o(1))n^{5/3})$ linear $3$-uniform hypergraphs on $n$ vertices, and almost all have $n^2/6-(1/6+o(1))n^{5/3}$ edges.
        \item \emph{Matroids of fixed rank.}
            For each fixed $r\ge4$, the number $m(n,r)$ of rank-$r$ matroids on $[n]$ satisfies
            \[
                \log m(n,r)=\frac{n^{r-1}}{r!}(\log n-r+1)
                +\left(\frac1{(r-1)!}+o(1)\right)n^{r-1-1/r},
            \]
            refining the estimates of van der Hofstad, Pendavingh and van der Pol.
    \end{enumerate}

\end{abstract}

\maketitle

\clearpage
\begingroup
\normalsize
\singlespacing
\makeatletter
\renewcommand{\l@section}{\@tocline{1}{4pt}{1pc}{}{}}
\renewcommand{\l@subsection}{\@tocline{2}{0pt}{3pc}{3pc}{}}
\patchcmd{\@tocline}
    {\hfil\hbox}
    {\nobreak\hskip.5em
     \leaders\hbox{\normalfont\kern.2em.\kern.2em}\hfill
     \hskip.5em\hbox}
    {}{\PackageError{matching-statistics}{Cannot set contents leaders}{Check the AMS contents definition.}}
\makeatother
\tableofcontents
\endgroup
\clearpage


\section{Introduction}\label{sec:intro}

A matching in a hypergraph is a collection of pairwise vertex-disjoint edges.
We use weighted enumeration to determine the mean and typical size of a random matching in a regular hypergraph, bound the variance of its size, and obtain large-deviation estimates for the number of uncovered vertices.
For a finite simple hypergraph $H$, let $\calM(H)$ denote the set of all matchings, including the empty matching.
For $\lambda>0$, define the matching partition function by
\[
    Z_H(\lambda)=\sum_{M\in\calM(H)}\lambda^{|M|}.
\]
The \emph{Gibbs matching measure} assigns probability
\[
    \Prob_\lambda(M)=\frac{\lambda^{|M|}}{Z_H(\lambda)}
\]
to each $M\in\calM(H)$.
The parameter $\lambda$ is called the \emph{activity}, or fugacity.
At activity one, the Gibbs measure is uniform and $Z_H(1)=|\calM(H)|$.

Our main result determines the first correction to the previously known leading estimate.
Since the hypotheses of our general counting theorem (Theorem~\ref{thm:weightedcount}) involve both higher codegrees and \emph{intersection patterns}, which we introduce in this paper, we first state the result for linear hypergraphs, in which any two distinct edges have at most one common vertex.

\begin{theorem}[Weighted counting in linear hypergraphs]\label{thm:linearcount}
    For fixed $k\ge3$ and $\eps>0$, let $H$ be an $n$-vertex $d$-regular linear $k$-uniform hypergraph with $d\ge n^\eps$.
    Uniformly over these hypergraphs and $\lambda$ in each fixed compact interval of $(0,\infty)$,
    \[
        \log Z_H(\lambda)=\frac nk\bigl(\log(\lambda d)-k+1\bigr)
        +(1+o(1))n(\lambda d)^{-1/k}.
    \]
\end{theorem}

In particular, the total number of matchings is
\[
    |\calM(H)|=\left(\frac d{e^{k-1}}\right)^{n/k}
    \exp\bigl((1+o(1))nd^{-1/k}\bigr).
\]
The positive correction $nd^{-1/k}$ gives the scale of the typical deficit from a perfect matching, whether or not a perfect matching exists.
To determine the mean deficit and its large deviations, we need the dependence of this correction on $\lambda$: changing the activity gives exponential moments of the number of uncovered vertices.
An error larger than the correction would obscure these moments.

\paragraph{Earlier counting results.}
Many problems in design theory can be expressed as matching problems in regular hypergraphs.
For example, a Steiner system is a perfect matching in the hypergraph on the $t$-sets in which each possible block gives an edge consisting of its $t$-subsets.
Similarly, representing each cell of a Latin square by its row, column and symbol identifies transversals with perfect matchings in a regular $3$-uniform hypergraph.
In both examples arbitrary matchings correspond to the respective partial objects.
Erd\H{o}s and Hanani~\cite{EH} asked whether, for fixed $t<k$, one can pack $k$-sets so that all but a vanishing proportion of the $t$-sets are covered, each at most once.
R\"odl~\cite{Rodl} proved this using the nibble. Frankl and R\"odl~\cite{FR} and Pippenger extended the method to nearly perfect matchings in hypergraphs; see also Pippenger and Spencer~\cite{PS}.
Quantitative matching bounds were obtained by Alon, Kim and Spencer~\cite{AKS}, Kang, K\"uhn, Methuku and Osthus~\cite{KKMO}, and Gould and Kelly~\cite{GK}.
For full designs, Keevash~\cite{Keevash} and Glock, K\"uhn, Lo and Osthus~\cite{GKLO} established existence under the necessary divisibility conditions for sufficiently large orders; Keevash~\cite{KeevashII} also treated more general decompositions, including resolvable designs.
These results concern the construction of large matchings; enumeration additionally requires control of the number of choices throughout the construction.

For graphs, Schrijver~\cite{Schrijver} proved a lower bound on the number of perfect matchings in regular bipartite graphs.
For lower bounds on the number of matchings of prescribed size in bipartite graphs, see~\cite{Csikvari,Lelarge}.
Carroll, Galvin and Tetali~\cite{CGT} and Ilinca and Kahn~\cite{IK} obtained entropy bounds for matchings of prescribed size, while Davies, Jenssen, Perkins and Roberts~\cite{DJPR} proved sharp bounds for the matching partition function.
See also~\cite{DJPRStability} for sharp coefficient bounds obtained from stability of the occupancy estimates.
For hypergraphs, Asratian and Kuzjurin~\cite{AsratianKuzjurin} counted nearly perfect matchings in almost regular hypergraphs with small codegrees.
Related entropy estimates appear in~\cite{Luria,DDK,KSW}.
The theory of matchings avoiding prescribed forbidden submatchings was developed by Delcourt and Postle~\cite{DP} and by Glock, Joos, Kim, K\"uhn and Lichev~\cite{GJKKL}.
The latter authors also proved a counting lower bound~\cite[Theorem~3.5]{GJKKL}; with no forbidden configurations, its logarithmic error is $O(nd^{-\eta^4})$ for a sufficiently small fixed $\eta>0$, which does not resolve the correction $nd^{-1/k}$.

Recently, Bennett and Frieze~\cite[Theorem~4(b)]{BF} proved
\[
    \log|\calM(H)|=\frac nk(\log d-k+1)+o(n)
\]
under a small pair-codegree assumption, also allowing the uniformity to grow.
For fixed $k$, their hypothesis is $\Delta_2(H)/d=o((\log n)^{-10})$, where $\Delta_2(H)$ is the maximum number of edges containing a given pair of vertices.
Earlier counting and entropy arguments already give this leading estimate in several of our applications.
The error $o(n)$ need not be small compared with $nd^{-1/k}$.
Our estimate determines the coefficient of $n(\lambda d)^{-1/k}$ uniformly on compact intervals of activities. This precision gives the expected number of uncovered vertices, concentration on that scale and the large-deviation rate.

For full Steiner triple systems, Linial and Luria~\cite{LinialLuria} obtained an entropy upper bound, and Keevash~\cite{KeevashCount} proved a matching lower bound on the logarithmic scale as part of his work on counting designs.
For partial systems, Asratian and Kuzjurin~\cite{AKPartial} counted packings without prescribing their size, and van der Hofstad, Pendavingh and van der Pol~\cite{HPP} obtained sharper estimates in several cases.
We compare these counts with ours in Section~\ref{subsec:designapplications}.

\paragraph{Vacancy and Kahn's conjecture.}
To study how far a random matching is from covering every vertex, we keep track of the vertices it misses.
If $M$ is sampled from the Gibbs matching measure on $H$ at activity $\lambda$, write
\[
    \calU_\lambda(H)=V(H)\setminus V(M)
\]
for its random set of uncovered vertices.
At activity one we omit the subscript and write $\calU(H)$.
For each $v\in V(H)$, define its \emph{vacancy probability} by
\[
    p_H(v)=\Prob(v\in\calU(H)).
\]
Then $\E[|\calU(H)|]=\sum_v p_H(v)$ and $\E[|M|]=(n-\E[|\calU(H)|])/k$.
In the 1990s, Kahn conjectured that regularity and linearity determine this probability to first order at every vertex.

\begin{conjecture}[Kahn~\cite{Kahn95,Kahn97}]\label{conj:Kahn}
    For each fixed $k\ge2$, if $H$ is a $d$-regular linear $k$-uniform hypergraph, then
    \[
        p_H(v)=(1+o_d(1))d^{-1/k}
    \]
    for every $v\in V(H)$, where the error tends to zero as $d\to\infty$, uniformly over $H$ and $v$.
\end{conjecture}

The scale $d^{-1/k}$ can be explained by a local independence heuristic; see~\cite[Section~1]{Lee}.
Inserting an edge gives $\Prob(e\in M)=\Prob(e\subseteq\calU(H))$, and hence
\begin{equation}\label{eq:vacancyheuristic}
    p_H(v)^{-1}=1+\sum_{\substack{e\in E(H)\\v\in e}}
    \Prob(e\setminus\{v\}\subseteq\calU(H)\mid v\in\calU(H)).
\end{equation}
If the conditional vacancies were approximately independent with common probability $x$, this would give $x=(1+dx^{k-1})^{-1}$, whose positive solution is asymptotic to $d^{-1/k}$.
Kahn and Kim~\cite{KK} proved the pointwise vacancy estimate for graphs and also showed that $\Var[|M|]\sim n/(4\sqrt d)$.
They also proposed the asymptotic formula $\Var[|M|]\sim n/(k^2d^{1/k})$ for linear hypergraphs.

For every $k\ge3$, the author~\cite{Lee} disproved Conjecture~\ref{conj:Kahn} by constructing regular linear hypergraphs with vacancy probabilities ranging from $(1+o(1))/(d+1)$ to $1-(1+o(1))d^{-(k-2)}$.
The construction iterates the rooted vacancy recursion
\[
    g_d(x)=\frac1{1+(d-1)x^{k-1}}.
\]
At its fixed point $\alpha_d\sim d^{-1/k}$, one has $g_d'(\alpha_d)=-(k-1)(1-\alpha_d)$.
Thus the fixed point is attracting for $k=2$, but repelling for $k\ge3$ and large $d$; the second iterate has two attracting fixed points forming a two-cycle for $g_d$.
The resulting finite counterexamples have order very large compared with degree.
This led the author to ask whether polynomial degree and small codegrees suffice~\cite[Conjecture~5.2]{Lee}.
For $S\subseteq V(H)$, let $d_H(S)$ count the edges containing $S$, and let $\Delta_j(H)=\max_{|S|=j}d_H(S)$.

\begin{conjecture}[Lee~{\cite[Conjecture~5.2]{Lee}}]\label{conj:Lee}
    For fixed $k\ge3$ and $0<\delta,\eps<1$, suppose that an $n$-vertex $k$-uniform hypergraph $H$ satisfies
    \[
        (1-n^{-\delta})d\le d_H(v)\le(1+n^{-\delta})d
    \]
    for every $v\in V(H)$, where $d>n^\eps$, and suppose that $\Delta_2(H)\le n^{-\delta}d$.
    Then $p_H(v)=(1+o_n(1))d^{-1/k}$ for every $v\in V(H)$, uniformly over these hypergraphs and vertices.
\end{conjecture}

Conjecture~\ref{conj:Lee} is false even after averaging; Appendix~\ref{sec:counterexample} gives an exactly regular counterexample.
Theorem~\ref{thm:statistics} establishes an averaged version of Kahn's conjecture for the hypergraphs covered by Theorem~\ref{thm:weightedcount}, including all regular linear hypergraphs of polynomial degree.

Without a relation between degree and order, Kahn and Kim~\cite[p.~204]{KK}, crediting Johansson, recorded $\E[|\calU(H)|]=o_d(n)$ for regular linear hypergraphs; this also follows from~\cite{AsratianKuzjurin}.
Gupta~\cite{Gupta} obtained $\E[|\calU(H)|]/n\le(\max\{3,k-1\}+o_d(1))\log\log d/\log d$.
Under polynomial degree, we identify the scale $d^{-1/k}$ and its coefficient.
For symmetric examples, the averaged result also gives pointwise vacancy estimates: these include partial Steiner and subspace designs, group multiplication tables, and toroidal queens configurations.
In the partite examples, transitivity within each class and equal occupancy of the classes suffice; see Section~\ref{sec:applications}.

\begin{remark}[Vacancy and occupancy]
The same statistic can also be expressed in terms of occupancy.
The edge occupancy fraction used by Davies, Jenssen, Perkins and Roberts~\cite{DJPR} is
\[
    \alpha_H^M(\lambda)=\frac{\E_\lambda[|M|]}{|E(H)|}.
\]
Since $|E(H)|=nd/k$, its relation to the average vacancy probability is
\begin{equation}\label{eq:occupancyvacancy}
    \frac{\E_\lambda[|\calU_\lambda(H)|]}{n}=1-d\alpha_H^M(\lambda).
\end{equation}
The expected proportion of covered vertices is therefore $d\alpha_H^M(\lambda)$.
The occupancy method estimates the logarithmic derivative of the partition function and then integrates; see also Cohen, Perkins, Sarantis and Tetali~\cite{CPST} for related hypergraph models.
We first estimate the partition function and then use changes of activity to obtain vacancy moments and tail bounds.
We use vacancy probabilities because the expected proportion of uncovered vertices is of order $(\lambda d)^{-1/k}$, the scale of the correction to $n^{-1}\log Z_H(\lambda)$.
Moreover, $|\calU_\lambda(H)|/k=n/k-|M|$ directly measures the deficit from a perfect matching.
\end{remark}


\subsection{The general counting criterion}\label{subsec:mainresults}

Here we state the general counting theorem, of which Theorem~\ref{thm:linearcount} is the linear case, and the matching statistics that follow from it.
To determine the first correction, we need two errors in the counting argument to be smaller than $nd^{-1/k}$.

The first arises in the lower bound.
In Section~\ref{sec:greedy} we introduce an iteration of weighted counting inequalities based on short random greedy stages and regular completions.
Building on the degree estimates of Bennett and Bohman~\cite{BB}, this gives Theorem~\ref{thm:greedycount}, a lower bound for the partition function over a polynomial range of activities.
Higher codegrees control how rapidly the available degrees can change during the process; see also the matching bounds of Gould and Kelly~\cite{GK}.
The parameter governing this error is
\begin{equation}\label{eq:Gamma}
    \Gamma(H)=\sqrt{\frac{\Delta_2(H)}d}
    +\max_{3\le j\le k}\left(\frac{\Delta_j(H)}d\right)^{1/(j-1)}.
\end{equation}
The condition $\Gamma(H)=O(d^{-1/k}n^{-\eps})$ is equivalent to
\[
    \Delta_2=O\bigl(d^{1-2/k}n^{-2\eps}\bigr)
\]
and, for every $3\le j\le k$,
\[
    \Delta_j=O\bigl(d^{1-(j-1)/k}n^{-(j-1)\eps}\bigr).
\]
These inequalities bound the number of edges containing a small set of vertices.
The stronger condition on the pair codegree comes from the square-root error in the greedy estimate.

The second error arises when we delete edges from a matching sampled at large activity.
To use the entropy gained by deletion, we must bound the number of ways to recover the original matching from the edges that remain.
An alternative edge on the newly uncovered vertices can be used only if every matching edge that it meets has been deleted.
The entropy estimates in Section~\ref{sec:patterns} therefore distinguish alternative edges according to these intersections.

\begin{definition}[Intersection coefficients]\label{def:coefficients}
    Let $\Pi_k$ consist of the partitions $\pi=(a_1,\ldots,a_m)$ of $k$ with $a_1\ge\cdots\ge a_m\ge1$ and $2\le m\le k-1$, and write $m(\pi)=m$.
    For a matching $P$, let $N_\pi(P)$ count edges $e\notin P$ contained in $V(P)$ whose nonempty intersections with the edges of $P$ have sizes $\pi$.
    Define
    \[
        b_\pi(H)=\sup_{P\ne\varnothing}\frac{N_\pi(P)}{|P|}.
    \]
    If $H$ has no edges, set $b_\pi(H)=0$.
\end{definition}

The pattern $(k)$ corresponds to an edge of the matching itself, whereas $(1,\ldots,1)$ corresponds to an edge meeting $k$ different matching edges.
The coefficients $b_\pi(H)$ count the remaining patterns, which do not occur in a linear hypergraph.
If each matching edge is deleted independently with probability $d^{-1/k}$, an edge meeting $m$ matching edges is contained in the deleted set with probability $d^{-m/k}$.
This explains the scale $d^{m/k}$ in the bound on $b_\pi(H)$.
We require a fixed polynomial margin in both conditions to absorb the logarithmic factors in the proof.

\begin{theorem}[Weighted counting]\label{thm:weightedcount}\label{thm:counting}
    For fixed $k\ge3$, $\eps>0$ and $C\ge1$, let $H$ be a simple $n$-vertex $d$-regular $k$-uniform hypergraph with $d\ge n^\eps$ satisfying
    \begin{equation}\label{eq:admissibility}
        \begin{gathered}
            \Gamma(H)\le C d^{-1/k}n^{-\eps},\\
            b_\pi(H)\le C d^{m(\pi)/k}n^{-\eps}\quad(\pi\in\Pi_k).
        \end{gathered}
    \end{equation}
    For every fixed compact interval $I\subset(0,\infty)$, uniformly over $H$ and $\lambda\in I$,
    \begin{equation}\label{eq:weightedcount}
        \log Z_H(\lambda)=\frac nk\bigl(\log(\lambda d)-k+1\bigr)
        +(1+o(1))n(\lambda d)^{-1/k}.
    \end{equation}
    The $o(1)$ term tends to zero as $n\to\infty$ for fixed $k,\eps,C,I$.
\end{theorem}

The same exponent $\eps$ is used for convenience; it can be decreased to accommodate any fixed positive polynomial margins.
At $\lambda=1$, the theorem gives
\begin{equation}\label{eq:totalcount}
    |\calM(H)|=\left(\frac d{e^{k-1}}\right)^{n/k}
    \exp\bigl((1+o(1))nd^{-1/k}\bigr).
\end{equation}
The error is in the logarithm of the count; a multiplicative $1+o(1)$ estimate is not asserted.

\begin{theorem}[Matching statistics]\label{thm:statistics}
    Let $H$ satisfy the hypotheses of Theorem~\ref{thm:weightedcount}.
    The following conclusions hold uniformly over these hypergraphs as $n\to\infty$.
    \begin{enumerate}[(a)]
        \item Uniformly for $\lambda$ in each fixed compact interval of positive activities,
            \[
                \frac{|\calU_\lambda(H)|}{n(\lambda d)^{-1/k}}\longrightarrow1\quad\text{in }L^2.
            \]
        \item For the uniform measure and every fixed $x>1$,
            \[
                \log\Prob(|\calU(H)|\ge xnd^{-1/k})
                =-\bigl(x\log x-x+1+o(1)\bigr)nd^{-1/k}.
            \]
            For every fixed $0<x<1$,
            \[
                \log\Prob(|\calU(H)|\le xnd^{-1/k})
                =-\bigl(x\log x-x+1+o(1)\bigr)nd^{-1/k}.
            \]
    \end{enumerate}
\end{theorem}

Part~(a) gives, uniformly on each fixed compact activity interval,
\[
    \E_\lambda[|\calU_\lambda(H)|]=(1+o(1))n(\lambda d)^{-1/k}
\]
and
\[
    \Var_\lambda\bigl[|\calU_\lambda(H)|\bigr]=o\bigl(n^2(\lambda d)^{-2/k}\bigr).
\]
In particular, under the uniform measure, almost all matchings have size $n/k-(1+o(1))nd^{-1/k}/k$, where \emph{almost all} means a proportion tending to one.
The sharper variance asymptotic and finer distributional questions are discussed in Section~\ref{sec:concluding}.

Letting $x\downarrow0$ after $n\to\infty$ in part~(b) gives, for a uniformly chosen matching $M$,
\begin{equation}\label{eq:perfect-probability}
    \Prob(M\text{ is perfect})\le\exp\bigl(-nd^{-1/k}+o(nd^{-1/k})\bigr).
\end{equation}
A corresponding lower bound requires additional assumptions; see Section~\ref{sec:concluding}.

Both $\Gamma(H)$ and the intersection coefficients admit bounds in terms of the maximum pair codegree alone.
This gives the following sufficient condition, which includes the linear case and the bounded-codegree hypergraphs used for toroidal queens.

\begin{corollary}[A pair-codegree condition]\label{cor:pairrange}
    For fixed $k\ge3$, $\eps>0$ and $0\le\gamma<1/k$, every simple $n$-vertex $d$-regular $k$-uniform hypergraph with $d\ge n^\eps$ and $\Delta_2(H)\le d^\gamma$ satisfies the hypotheses of Theorem~\ref{thm:counting} after decreasing $\eps$ if necessary.
    Hence the counting and statistical conclusions above hold uniformly over these hypergraphs.
\end{corollary}

\begin{proof}[Proof of Corollary~\ref{cor:pairrange}]
    Since $\Delta_j(H)\le\Delta_2(H)\le d^\gamma$ for $2\le j\le k$, definition~\eqref{eq:Gamma} and $d\ge n^\eps$ give
    \[
        \Gamma(H)\le2d^{-(1-\gamma)/(k-1)}
        \le2d^{-1/k}n^{-\eps(1-k\gamma)/(k(k-1))}.
    \]
    Every $\pi\in\Pi_k$ has a part of size at least two, so each edge counted by $N_\pi(P)$ contains a pair of vertices from one edge of $P$.
    There are $\binom{k}{2}|P|$ such pairs, each contained in at most $\Delta_2(H)$ edges.
    Dividing by $|P|$ and taking the supremum over nonempty matchings $P$, we obtain
    \[
        b_\pi(H)\le\binom{k}{2}d^\gamma
        \le\binom{k}{2}d^{m(\pi)/k}n^{-\eps(2/k-\gamma)},
    \]
    where the last inequality uses $m(\pi)\ge2$ and $d\ge n^\eps$.
    Both exponents of $n$ are negative because $\gamma<1/k$, so~\eqref{eq:admissibility} holds after decreasing $\eps$, with a constant depending only on $k$.
    The conclusions follow from Theorems~\ref{thm:counting} and~\ref{thm:statistics}.
    This completes the proof.
\end{proof}

For partial designs, the intersection patterns give information that the pair codegree alone does not capture.
Section~\ref{sec:applications} gives examples of both kinds, ranging from Latin squares and toroidal queens to partial designs and matroids.

\begin{remark}[Extending the codegree range]\label{rem:codegree-extension}
    Substantially extending the codegree range of Theorem~\ref{thm:weightedcount} would have consequences for quantitative maximum-matching problems.
    Indeed, even the one-sided estimate $\E[|\calU(H)|]=O(nd^{-1/k})$ under the weaker hypotheses considered in Appendix~\ref{sec:conditionallifting} would improve bounds on the maximum matching size in regular linear hypergraphs of polynomial degree.
    For recent advances in matching bounds using higher codegrees, see Gould and Kelly~\cite{GK}.
    We do not claim that the present criterion is optimal; the remaining pair-codegree gap is discussed in Section~\ref{sec:concluding}, under ``The pair-codegree range'', and the precise conditional reduction is proved in Appendix~\ref{sec:conditionallifting}.
\end{remark}


\subsection{Proof overview}\label{subsec:overview}

We first prove the weighted counting formula and then deduce the matching statistics.
The main difficulty is to determine the correction $n(\lambda d)^{-1/k}$ with an error smaller than this quantity.
Our proof combines a new greedy counting argument based on regular completions with an entropy method based on the intersection patterns introduced above.

We begin with Theorem~\ref{thm:greedycount}, proved in Section~\ref{sec:greedy}, which gives a lower bound for $Z_H(\lambda)$ uniformly over a polynomial range of activities.
Each greedy stage stops when half of its vertices remain, so the degree estimate of Bennett and Bohman~\cite{BB} is used only over a short interval.
We also control the higher codegrees by applying Freedman's inequality after dividing each codegree process by its deterministic decay factor.
We count the resulting matching prefixes, assign each prefix $P$ the weight $\lambda^{|P|}$, and sum over all matching extensions, accounting for the different prefixes contained in the same matching.
We then complete each residual hypergraph to an exactly regular one and repeat, using an encoding to control the cost of transferring weighted counts through the completion.
The key point is that, after subtracting the leading term and dividing by the number of vertices, the remainder of the next stage enters with coefficient $1/2+o(1)$.
Consequently, errors arising later in the iteration receive geometrically decreasing weights.
Together with the higher-codegree estimates, this bounds the total logarithmic loss by $O(n\Gamma(H)(\log n)^L)$, which is smaller than the desired correction under the hypotheses of Theorem~\ref{thm:weightedcount}.
This preliminary estimate may be of independent interest: it holds over a polynomial activity range without assumptions on intersection patterns.
In particular, we can apply it at activity $\lambda d$, giving the source count needed for the thinning argument below.

To obtain the positive correction, sample a Gibbs matching $M$ at activity $\lambda d$ and delete each edge independently with an appropriate probability $\theta$, leaving a matching $J$.
Writing $\HH$ for entropy, the entropy chain rule gives
\[
    \HH(J)=\HH(M)+\E[|M|]h(\theta)-\HH(M\mid J),
\]
where $h$ is the binary entropy function.
The middle term records the deletion choices, and the last term accounts for different original matchings producing the same $J$.
The Gibbs variational inequality bounds $\log Z_H(\lambda)$ from below by $\HH(J)+\E[|J|]\log\lambda$.
Thus, a reconstruction bound for $\HH(M\mid J)$ converts the deletion choices into a weighted counting gain.

Our intersection patterns control this reconstruction.
After specifying which vertices were already uncovered in $M$, the deleted edges form a perfect matching on the newly uncovered vertices.
An alternative edge meeting $m$ distinct edges of $M$ is available only when all $m$ are deleted, an event of probability $\theta^m$.
Its pattern records the sizes of these intersections: the number of parts determines this probability, while their sizes determine which codegrees bound the number of alternatives.
This records the dependence between vertices of the same matching edge.
We expose the matching in a random order, following Radhakrishnan's entropy proof of the Bregman--Minc bound on permanents~\cite{Radhakrishnan} and its use by Linial and Luria for Steiner triple systems~\cite{LinialLuria} and by Luria for hypergraph perfect matchings~\cite{Luria}.
Here the coefficients $b_\pi(H)$ bound the number of alternatives after averaging over vertices and source matchings; integrating the logarithm of this bound estimates the conditional entropy of the deleted matching.

The source matching may already leave vertices uncovered, and its average vacancy probability is unknown.
The higher activity contributes a term proportional to the number of these original vacancies.
A tangent-line bound for binary entropy uses this term to control the entropy of distinguishing original vacancies from vertices uncovered by deletion.
The resulting lower bound is uniform in this unknown average vacancy probability.
Later, we will show that choosing $\theta=(\lambda d)^{-1/k}$ gives the correction $n(\lambda d)^{-1/k}$.

For the upper bound, we describe a Gibbs matching by its uncovered set and a perfect matching on the remaining vertices.
The same random-order exposure and pattern estimates bound the available matching edges.
The coefficient one requires retaining the degree loss caused by uncovered vertices.
Several vacancies in one edge make these losses overlap.
We express the degree loss in terms of the number of edges meeting the uncovered set.
Shearer's inequality and the Gibbs variational inequality then bound the entropy of that set together with this loss, one edge at a time.
Choosing the auxiliary weights in terms of $(\lambda d)^{-1/k}$ makes the terms involving the mean vacancy density cancel.
Summing the resulting edge estimates gives the correction $n(\lambda d)^{-1/k}$, with a Taylor error of order $n(\lambda d)^{-2/k}$.

Finally, changes of activity express the exponential moments of the vacancy count as ratios of partition functions.
The counting estimate then gives tail bounds and $L^2$ convergence, yielding both the sharp mean and concentration.
Concentration under the corresponding tilted laws gives the large-deviation lower bounds.


\section{Applications}\label{sec:applications}

Our main results give refined counting formulae and typical-size estimates for several families of combinatorial objects.
We illustrate these applications with partial Latin squares and transversals, partial toroidal queens configurations, matchings in Steiner systems, partial designs over sets and vector spaces, and matroids of fixed rank.
All objects and their underlying points are labelled.
For partial objects, unless a host is explicitly fixed, the uniform measure is taken over all valid collections, including the empty collection, without prescribing their size.
The proofs are collected in Section~\ref{sec:applicationproofs}.

We use $N,K,D$ for the order, uniformity and degree of an auxiliary hypergraph $\calH$, and all structural parameters are fixed as the underlying size tends to infinity.
The partial objects and matchings in fixed designs are encoded bijectively as matchings in $\calH$, preserving their counts and uniform measures.
The auxiliary degrees are at least a fixed positive power of the order, so checking the codegrees and intersection coefficients allows us to apply Theorems~\ref{thm:weightedcount} and~\ref{thm:statistics}.
For matroids, the bijection applies to sparse paving matroids; additional counting comparisons give the result for paving and arbitrary matroids.
The typical-size conclusions below are stated as $L^2$ convergence.
They therefore give the expectation asymptotic and show that the variance is of smaller order than the square of the typical deficit.
In the leading term, we replace the exact auxiliary degree by an asymptotic equivalent only when the resulting error is $o(ND^{-1/K})$.
For the families encoded bijectively as matchings, Theorem~\ref{thm:statistics}(b) also gives large deviations for the number of uncovered auxiliary vertices.
If the deficit is $|\calU(\calH)|/K$, the logarithmic tail estimate has leading term $-ND^{-1/K}(x\log x-x+1)$.


\subsection{Latin squares and toroidal queens}\label{subsec:latinqueens}

In these examples, each entry or cell occupies one vertex for each constraint, so compatible entries or cells correspond to disjoint edges.

\paragraph{Partial Latin squares.}

A Latin square of order $n$ is an $n\times n$ array on $n$ symbols in which each symbol occurs exactly once in every row and column.
The number of Latin squares of order $n$ is $((1+o(1))n/e^2)^{n^2}$, as follows from permanent bounds; see~\cite[Theorem~17.3]{vanLintWilson}.
For exact enumeration, see McKay and Wanless~\cite{McKayWanless}.
For entropy bounds in the more general setting of high-dimensional permutations, see~\cite{LinialLuriaPermutations}.
A \emph{partial Latin square} allows empty cells and requires each symbol to occur at most once in every row and column.
Its size is the number of filled cells.
We count all such arrays on labelled rows, columns, and symbols, without requiring a completion to a Latin square.

Ghandehari, Hatami and Mahmoodian~\cite[Theorem~3]{GHM} gave an upper bound for the number of partial Latin squares with a prescribed number of filled cells.
For exact enumeration of partial Latin rectangles, including small orders and small numbers of filled cells, see Falc\'on and Stones~\cite{FalconStones}.
The leading logarithmic estimate $n^2 \log (n/e^2) +o(n^2)$ also follows from the recent result of Bennett and Frieze~\cite[Theorem~4(b)]{BF}.
We note that summing the bound of Ghandehari, Hatami and Mahmoodian over all sizes already gives the upper bound with correction $(3+o(1))n^{5/3}$ below.

\begin{corollary}[Partial Latin squares]\label{cor:partiallatin}
    Let $P$ be chosen uniformly from the family $\mathcal{P}_n$ of all partial Latin squares of order $n$.
    As $n\to\infty$,
    \[
        \frac{n^2-|P|}{n^{5/3}}\longrightarrow1\quad\text{in }L^2
    \]
    and
    \[
        |\mathcal{P}_n| =  e^{(3 + o(1))n^{5/3}} \left( \frac{n}{e^2} \right)^{n^2}.
    \]
\end{corollary}

Our result gives the matching lower bound and determines the typical number of empty cells.
In particular, the mean number of empty cells is $(1+o(1))n^{5/3}$ and its variance is $o(n^{10/3})$.
Every fixed cell is empty with probability $(1+o(1))n^{-1/3}$.

\paragraph{Partial transversals in Latin squares.}

A \emph{partial transversal} of a Latin square is a set of cells with no repeated row, column or symbol.
A transversal has $n$ cells.

The conjectures of Ryser and of Brualdi and Stein assert the existence of a transversal in every Latin square of odd order and a partial transversal with at least $n-1$ cells in every Latin square of order $n$; see~\cite{Ryser,Stein}.
After the bound of Keevash, Pokrovskiy, Sudakov and Yepremyan~\cite{KPSY}, Montgomery~\cite{Montgomery} proved the latter assertion for all sufficiently large orders.
For full transversals, Glebov and Luria~\cite{GL} proved that the maximum over order-$n$ Latin squares is $((1+o(1))n/e^2)^n$; Kwan~\cite[Theorem~1.3]{Kwan} showed that almost all Latin squares attain this leading count.
The leading estimate $\log|\calT(L)|=n\log(n/e^2)+o(n)$ for the family $\calT(L)$ of all partial transversals follows from earlier matching-counting and entropy arguments~\cite{GJKKL,Luria} and directly from Bennett and Frieze~\cite[Theorem~4(b)]{BF}.
Tang and Zhao~\cite{TZ} count partial transversals with $n-1$ cells with logarithmic error $O(n^{1-c})$ for some $c>0$; their stated bounds do not specify error $o(n^{2/3})$ throughout the typical-deficit window considered here.

\begin{corollary}[Partial transversals]\label{cor:latin}
    Let $L$ be a Latin square of order $n$, and let $T$ be chosen uniformly from the family $\calT(L)$ of all partial transversals.
    As $n\to\infty$, uniformly over $L$,
    \[
        \frac{n-|T|}{n^{2/3}}\longrightarrow1\quad\text{in }L^2
    \]
    and
    \[
        |\calT(L)| = e^{(3 + o(1))n^{2/3}} \left(\frac{n}{e^2} \right)^n.
    \]
\end{corollary}

The typical deficit $(1+o(1))n^{2/3}$ is substantially larger than the smallest possible deficit.
If $L$ is a group multiplication table, each fixed row, column or symbol is missed with probability $(1+o(1))n^{-1/3}$.

\paragraph{Partial toroidal queens configurations.}

The classical $n$-queens problem asks for placements of $n$ nonattacking queens on the ordinary $n\times n$ board.
Simkin~\cite{Simkin} proved that their number is $((1+o(1))ne^{-\alpha})^n$, where $\alpha\approx1.94$ is determined by a variational problem.
On the toroidal board, a cell $(i,j)\in\mathbb Z_n^2$ belongs to row $i$, column $j$, and diagonals indexed by $i+j$ and $i-j$.
A partial configuration is a set of cells containing no two queens on the same row, column or either diagonal.
P\'olya's existence criterion for a full configuration is $\gcd(n,6)=1$; see~\cite{BK}.
Luria~\cite{Luria} obtained the counting upper bound, and Bowtell and Keevash~\cite{BK} proved that the number of full configurations in those orders is $((1+o(1))n/e^3)^n$.

\begin{corollary}\label{cor:queens}
    Let $Q$ be the number of queens in a uniformly chosen partial toroidal queens configuration, and let $\calQ_n$ denote the family of all partial configurations.
    As $n\to\infty$ through arbitrary integers,
    \[
        \frac{n-Q}{n^{3/4}} \to 1 \quad \text{in } L^2
    \]
    and
    \[
        |\calQ_n| = e^{(4 + o(1))n^{3/4}} \left( \frac{n}{e^3} \right)^n.
    \]
\end{corollary}

Every fixed row, column or diagonal is empty with probability $(1+o(1))n^{-1/4}$.


\subsection{Steiner systems and partial designs}\label{subsec:designapplications}

We distinguish matchings of point-disjoint blocks in a fixed design from partial designs chosen from all valid block collections.
These are different uniform measures, and we specify each family before applying the matching estimates.

\paragraph{Matchings in Steiner systems.}\label{subsec:steiner}

A \emph{Steiner system} $S(n,k,t)$ is a collection of $k$-subsets, called blocks, of an $n$-element set such that each $t$-subset belongs to exactly one block.
A Steiner triple system is an $S(n,3,2)$.
We first fix such a system and regard its points as vertices and its blocks as hyperedges.
A matching is a collection of point-disjoint blocks; a perfect matching is a parallel class.
The common degree is
\begin{equation}\label{eq:steiner-degree}
    d= \binom{n-1}{t-1}/\binom{k-1}{t-1}.
\end{equation}

Kwan~\cite[Theorems~1.1 and~1.2]{Kwan} proved that almost all Steiner triple systems of order $n\equiv3\pmod6$ have $((1+o(1))n/(2e^2))^{n/3}$ perfect matchings, which is also the asymptotic maximum.
Montgomery~\cite{Montgomery} proved that every sufficiently large Steiner triple system has a matching leaving at most four points uncovered.

\begin{corollary}\label{cor:steiner}
    For fixed $2\le t<k$ with $t<k/2+1$ and $\gcd(k,t-1)=1$, let $H$ be any Steiner system $S(n,k,t)$, and choose $M$ uniformly from $\calM(H)$.
    With $d$ as in~\eqref{eq:steiner-degree}, uniformly over $H$,
    \[
        \frac{n-k|M|}{nd^{-1/k}}\longrightarrow1\quad\text{in }L^2.
    \]
    The total number of matchings satisfies
    \[
        \log|\calM(H)|=\frac nk(\log d-k+1)+nd^{-1/k}+o(nd^{-1/k}),
    \]
    uniformly over $H$.
    The same conclusions hold for regular $k$-uniform hypergraphs with $d=\Theta(n^{t-1})$ and $\Delta_t=O(1)$.
\end{corollary}

For example, in every Steiner triple system, almost all matchings have
\[
    \frac n3-\left(\frac{2^{1/3}}3+o(1)\right)n^{2/3}
\]
blocks.
Our conclusion describes a uniformly chosen matching in every fixed system, including those with no perfect matching.
The corollary also covers $S(n,k,3)$ for every odd $k\ge5$.
The coprimality assumption ensures a strict margin in the pattern counts.

\paragraph{Uniformly chosen partial Steiner systems.}\label{subsec:partialsteiner}
A partial Steiner system $S(n,k,t)$ is a collection of $k$-subsets of $[n]$ in which each $t$-subset belongs to at most one block.
Write $\mathfrak S(n,k,t)$ for the family of all such collections.
The auxiliary hypergraph $\calH$ has the $t$-sets as vertices and the $t$-shadow of each $k$-set as an edge.
This is linear when $k=t+1$, but it is nonlinear for larger block sizes.
For fixed $2\le t<k$, its order is $N=\binom nt$, its uniformity is $K=\binom kt$, and its degree is $D=\binom{n-t}{k-t}$.

Asratian and Kuzjurin~\cite{AKPartial}, giving a simpler proof of an earlier result of Grable, obtained $\log|\mathfrak S(n,k,t)|\sim(N/K)\log D$.
Van der Hofstad, Pendavingh and van der Pol~\cite[Theorem~1.3]{HPP} proved the sharper upper bound
\[
    \log|\mathfrak S(n,k,t)|\le\frac NK(\log D-K+1)+o(N),
\]
and their Theorem~1.4 gives equality up to $o(N)$ when $k=t+1$.
In particular, their explicit formula for partial Steiner triple systems is
\[
    \log|\mathfrak S(n,3,2)|=\frac{n(n-1)}6\bigl(\log(n-2)-2\bigr)+o(n^2).
\]
For general fixed $t<k$, the equality with error $o(N)$ holds for every sufficiently large $n$; see~\cite[Theorem~4(b)]{BF}.
For related enumeration of full designs and sparse partial systems with a prescribed number of blocks, see~\cite{KeevashCount,McKayTian}.

\begin{theorem}\label{thm:partialdesigns}\label{cor:partialsteiner}
    For fixed $t\ge2$ and $k\in\{t+1,t+2,t+3\}$, let $N=\binom nt$ be the number of $t$-sets and $K=\binom kt$ the number of $t$-sets in each block.
    Each $t$-set is contained in $D=\binom{n-t}{k-t}$ possible blocks.
    Choose $\calS$ uniformly from $\mathfrak S(n,k,t)$, so that $\calU(\calH)$ is the set of $t$-sets contained in no block of $\calS$. As $n\to\infty$,
    \[
        \frac{|\calU(\calH)|}{ND^{-1/K}}\to1\quad\text{in }L^2.
    \]
    Furthermore,
    \[
        \log|\mathfrak S(n,k,t)| =\frac NK(\log D-K+1)+(1+o(1))ND^{-1/K}.
    \]
\end{theorem}

No divisibility assumption on $n$ is required, and almost all such systems have $\frac NK-\frac{1+o(1)}KND^{-1/K}$ blocks.
For partial Steiner triple systems, this gives
\[
    \log|\mathfrak S(n,3,2)|=\frac{n(n-1)}6\bigl(\log(n-2)-2\bigr)
    +\left(\frac12+o(1)\right)n^{5/3}.
\]
Thus it determines a positive correction, with its coefficient, that is not resolved by the earlier $o(n^2)$ error.

For $k=t+1$, the statistical conclusion gives
\[
    \E[|\calS|] =\frac1{t+1}\binom nt-\frac{1+o(1)}{t+1}\binom nt(n-t)^{-1/(t+1)}.
\]
In particular, a uniformly chosen partial Steiner triple system has
\[
    \E[|\calS|]=\frac{n(n-1)}6-\frac{1+o(1)}6n^{5/3}
\]
and
\[
    \Var[|\calS|]=o(n^{10/3}),
\]
and almost all such systems leave $(1+o(1))n^{5/3}/2$ pairs uncovered.
This is distinct from choosing a matching of point-disjoint blocks in a fixed Steiner system.
For partial systems generated by random greedy triangle removal, whose distribution differs from the uniform measure considered here, see~\cite{BFL}.

\begin{remark}[Linear $3$-uniform hypergraphs]\label{rem:linear-three-graphs}
    A partial Steiner triple system on $[n]$ is exactly a linear $3$-uniform hypergraph on $[n]$, since no pair of vertices belongs to two edges.
    Thus Theorem~\ref{thm:partialdesigns} counts
    \[
        (n/e^2)^{n^2/6}\exp\bigl((1/2+o(1))n^{5/3}\bigr)
    \]
    such hypergraphs, and almost all have
    \[
        \frac{n^2}{6}-\left(\frac16+o(1)\right)n^{5/3}
    \]
    edges.
\end{remark}

Since the auxiliary hypergraph is vertex-transitive, Theorem~\ref{thm:partialdesigns} also determines the vacancy probability of every fixed $t$-set:
\[
    \Prob(T\text{ belongs to no block})=(1+o(1))D^{-1/K}.
\]
Partial $S(n,r,r-1)$ systems are precisely the families of circuit-hyperplanes of rank-$r$ sparse paving matroids; we use this correspondence in Subsection~\ref{subsec:matroids}.

For $k=t+2$, at most one nonsingleton intersection occurs; for $k=t+3$, two such intersections of size $t+1$ already determine the candidate block.
When $k\ge t+4$, they need not do so, and more configurations must be counted.
We do not claim that the conclusion fails outside the stated range.

\paragraph{Partial subspace designs.}\label{par:subspaces}
There is a vector-space analogue of the packing problem above.
Write
\[
    \qbinom{h}{t}=\prod_{i=0}^{t-1}\frac{q^{h-i}-1}{q^{t-i}-1}
\]
for the number of $t$-dimensional subspaces of $\mathbb F_q^h$.
A \emph{partial subspace design} is a family of $(t+s)$-dimensional subspaces in which every $t$-dimensional subspace is contained in at most one block.
These objects are also called subspace packings; see Etzion, Kurz, Otal and \"Ozbudak~\cite{EKOO}.
Keevash, Sah and Sawhney~\cite{KSS} proved the existence and approximate enumeration of full subspace designs for all sufficiently large dimensions satisfying the necessary divisibility conditions.
For earlier explicit constructions of full subspace designs, see~\cite{BEOVW}.
Here we choose uniformly from partial designs of all sizes.
The auxiliary hypergraph $\calH$ has the $t$-dimensional subspaces as vertices, and one edge for the $t$-subspaces in each possible block.
Thus $\calU(\calH)$ consists of the $t$-subspaces covered by no chosen block.
The auxiliary hypergraph satisfies $\Delta_2(\calH)/D=O(q^{-h})$, so the leading counting estimate also follows from~\cite[Theorem~4(b)]{BF}.

\begin{theorem}[Partial subspace designs]\label{thm:subspaces}
    Fix a prime power $q$ and positive integers $t,s$ such that either $t\ge2$ and $s\in\{1,2,3\}$, or $t=1$ and $s\in\{1,2\}$.
    Let $\mathfrak D_q(h,t+s,t)$ be the family of all partial subspace designs with these parameters in $\mathbb F_q^h$.
    Write $N=\qbinom ht$, $K=\qbinom{t+s}t$ and $D=\qbinom{h-t}s$ for the order, uniformity and degree of the auxiliary hypergraph $\calH$.
    For a uniformly chosen $\calS\in\mathfrak D_q(h,t+s,t)$, as $h\to\infty$,
    \[
        \frac{|\calU(\calH)|}{ND^{-1/K}}\longrightarrow1\quad\text{in }L^2.
    \]
    Moreover,
    \[
        \log|\mathfrak D_q(h,t+s,t)| =\frac NK(\log D-K+1)+(1+o(1))ND^{-1/K}.
    \]
\end{theorem}

No divisibility assumption on $h$ is required, and almost all partial subspace designs have $\frac NK-\frac{1+o(1)}KND^{-1/K}$ blocks.

The action of $\operatorname{GL}(h,q)$ is transitive on the $t$-dimensional subspaces, so each is uncovered with probability $(1+o(1))D^{-1/K}$.
When $t=1$, the theorem covers partial spreads of two-dimensional and three-dimensional subspaces.
The restriction on $s$ has the same role as in the set version: for larger $s$, two nonsingleton intersections need not span the candidate block, and the pattern estimates require further work.
For $t=1,s=3$, the elementary bound on the number of nonsingleton intersections also does not give the required strict polynomial margin.


\subsection{Matroids of fixed rank}\label{subsec:matroids}

We refine the count of matroids of each fixed rank at least four on the labelled ground set $[n]$.

\Needspace{5\baselineskip}
A matroid on a finite set $E$ has an integer-valued rank function $\rho$ on its subsets, satisfying $0\le\rho(B)\le|B|$, monotonicity, and
\[
    \rho(B\cup C)+\rho(B\cap C)\le\rho(B)+\rho(C).
\]
Its rank is $\rho(E)$; a set $B$ is independent if $\rho(B)=|B|$, a basis is a maximal independent set, and a circuit is a minimal dependent set.
A hyperplane is a maximal set of rank $\rho(E)-1$.
We refer to Oxley~\cite{Oxley} for these standard definitions.

A rank-$r$ matroid is \emph{paving} if every circuit has at least $r$ elements, and \emph{sparse paving} if, in addition, every hyperplane has at most $r$ elements.
A paving matroid is determined by its hyperplanes: these are proper subsets of $[n]$ of size at least $r-1$, and every $(r-1)$-set is contained in exactly one of them.
For a sparse paving matroid, its hyperplanes of size $r$ are also circuits, called \emph{circuit-hyperplanes}, and form a partial $S(n,r,r-1)$ system.
Conversely, every such partial system specifies a unique sparse paving matroid; see~\cite{HPP,PvdP}.

For enumeration without a restriction on rank, see~\cite{Knuth,BPvdP,MNWW,PvdPComparison}.

Write $s(n,r)$, $p(n,r)$ and $m(n,r)$ for the numbers of sparse paving, paving and arbitrary matroids of rank $r$ on $[n]$, respectively.
Van der Hofstad, Pendavingh and van der Pol~\cite[Theorems~1.4, 1.5 and~1.7]{HPP} determined the logarithms of all three counts up to an error $o(n^{r-1})$ for fixed $r\ge4$.

Recently, Bennett and Frieze~\cite[Theorems~2 and~3]{BF} extended these leading estimates to certain growing ranks.
Kwan, Sah and Sawhney~\cite{KSSMatroids} showed that the mixture of three-point and four-point lines changes the leading constant for rank-$3$ matroids.

\begin{theorem}[Counting matroids of fixed rank]\label{thm:matroids}
    For each fixed $r\ge4$ and each $f\in\{s,p,m\}$, as $n\to\infty$,
    \[
        \log f(n,r)=\frac{n^{r-1}}{r!}(\log n-r+1)
        +\left(\frac1{(r-1)!}+o(1)\right)n^{r-1-1/r}.
    \]
\end{theorem}

The correction from Theorem~\ref{thm:partialdesigns} survives the passage from sparse paving matroids to all matroids.
For the paving upper bound, we fix the hyperplanes of size greater than $r$ and delete the $(r-1)$-sets they cover from the auxiliary hypergraph.
The hyperplanes of size $r$ form a matching in the remaining linear hypergraph, which has maximum degree at most $n-r+1$ but need not be regular.
Proposition~\ref{prop:maxdegree}, proved in Appendix~\ref{sec:maxdegreeproof}, bounds the number of such matchings.

The sparse paving formula also holds for $r=3$, but it does not extend to all matroids: the contribution of four-point hyperplanes is already of order $n^2$ in the upper bound.
We discuss the unresolved next correction in Section~\ref{sec:concluding}.

\begin{remark}
    For a fixed prime power $q$ and fixed $r\ge2$, a rank-$r$ sparse paving $q$-matroid on $\mathbb F_q^h$, with $h>r$, is uniquely determined by its $r$-dimensional circuits.
    These families of circuits are precisely the partial subspace designs in $\mathfrak D_q(h,r,r-1)$: each $(r-1)$-dimensional subspace lies in at most one member of the family; see~\cite[Definition~5 and Proposition~12]{JJP}.
    Theorem~\ref{thm:subspaces} with $t=r-1$ and $s=1$ therefore gives a refined logarithmic count of these matroids as $h\to\infty$.
\end{remark}


\section{Preliminaries}\label{sec:prelim}


Unless stated otherwise, $H$ is a simple $n$-vertex $d$-regular $k$-uniform hypergraph, where $k\ge3$ is fixed.
For $S\subseteq V(H)$, let $d_H(S)$ be the number of edges containing $S$, and let $\Delta_j(H)=\max_{|S|=j}d_H(S)$.
We use the codegree and intersection-pattern parameters defined in Section~\ref{subsec:mainresults}.
All logarithms are natural, and asymptotic statements are as $n\to\infty$ unless another limit is specified.
The parameters $k,C,\eps$, all positive constants appearing as exponents, and compact activity intervals are fixed unless otherwise specified.
The constant $L=L(k)$ is chosen sufficiently large for the greedy counting estimate; its value is not optimized.

We use the notation from Section~\ref{sec:intro}: for $M$ drawn from $\Prob_\lambda$ on $H$, the uncovered set is $\calU_\lambda(H)=V(H)\setminus V(M)$ and $|\calU_\lambda(H)|=n-k|M|$.

Throughout, $I\subset(0,\infty)$ denotes a fixed compact interval.


\subsection{Partition functions, the Gibbs measure and fluctuations}\label{subsec:physics}

A matching can be viewed as a collection of particles, each occupying the vertices of one edge.
The matching condition is a hard-core exclusion rule: two particles cannot occupy the same vertex.
On a graph, the occupied edges are called dimers and the uncovered vertices monomers.
This is the monomer-dimer model studied by Heilmann and Lieb~\cite{HL}.
The same interpretation applies to a $k$-uniform hypergraph, with each particle occupying $k$ vertices.

The activity $\lambda>0$, also called the fugacity, is the weight assigned to each particle.
Thus a configuration $M$ has weight $\lambda^{|M|}$.
The Gibbs measure is obtained by dividing each weight by the sum of all weights, the partition function $Z_H(\lambda)$.
This is a grand-canonical description, meaning that the number of particles is allowed to vary.

To see the counting content explicitly, let $c_m(H)$ be the number of matchings with exactly $m$ edges.
Then
\begin{equation*}
    Z_H(\lambda)=\sum_{m=0}^{\lfloor n/k\rfloor}c_m(H)\lambda^m
\end{equation*}
and
\begin{equation*}
    \Prob_\lambda(|M|=m)=\frac{c_m(H)\lambda^m}{Z_H(\lambda)}.
\end{equation*}
Conditional on $|M|=m$, the Gibbs measure is uniform over the matchings of that size.
At $\lambda=1$ every matching has the same probability, but sizes with more matchings occur more frequently.
Increasing $\lambda$ gives greater relative weight to larger matchings.

The logarithm of the partition function connects the count to the mean and variance.
Let $F_H(t)=\log Z_H(e^t)$, where $t=\log\lambda$.
The following standard identities explain why we estimate the partition function over an interval of activities.
The derivative identities are the usual cumulant identities for exponential families; see Wainwright and Jordan~\cite[Proposition~3.1]{WainwrightJordan}.

\begin{fact}[Activity identities]\label{fact:activity}
    For the Gibbs matching measure on a finite $k$-uniform hypergraph,
    \[
        F_H'(t)=\E_{e^t}[|M|]
    \]
    and
    \[
        F_H''(t)=\Var_{e^t}[|M|].
    \]
    For every $\lambda>0$ and $s\in\mathbb R$,
    \begin{equation}\label{eq:vacancy-mgf}
        \E_\lambda\bigl[e^{s|\calU_\lambda(H)|}\bigr]
        =e^{ns}\frac{Z_H(\lambda e^{-ks})}{Z_H(\lambda)}.
    \end{equation}
\end{fact}

For completeness, we include a proof.

\begin{proof}[Proof of Fact~\ref{fact:activity}]
    Differentiating the finite sum defining $Z_H(e^t)$ gives
    \[
        F_H'(t)=\frac{\sum_m m c_m(H)e^{tm}}{Z_H(e^t)}=\E_{e^t}[|M|].
    \]
    A second derivative gives the variance.
    Finally, multiplying the weight of a matching by $e^{s|\calU_\lambda(H)|}=e^{ns}e^{-ks|M|}$ replaces its activity by $\lambda e^{-ks}$, proving~\eqref{eq:vacancy-mgf}.
    This completes the proof.
\end{proof}

Thus the slope of $F_H$ records the mean number of particles, and its curvature records the variance. In particular, $F_H$ is convex.
Since $|\calU_\lambda(H)|=n-k|M|$, these identities also give
\[
    \frac{\dd}{\dd\log\lambda}\E_\lambda[|\calU_\lambda(H)|]
    =-\frac1k\Var_\lambda\bigl[|\calU_\lambda(H)|\bigr].
\]
This is an elementary fluctuation-response identity: mean vacancy decreases as activity increases, at a rate proportional to the variance.

For example, if $H$ consists of $q$ pairwise disjoint edges, then $Z_H(\lambda)=(1+\lambda)^q$.
Each edge is independently present with probability $\lambda/(1+\lambda)$, and
\[
    \E_\lambda[|M|]=\frac{q\lambda}{1+\lambda}
\]
and
\[
    \Var_\lambda[|M|]=\frac{q\lambda}{(1+\lambda)^2}.
\]
Fact~\ref{fact:activity} holds without independence between edges.
An asymptotic estimate for the log partition function cannot be differentiated term by term without controlling the derivatives of its error.
We instead use~\eqref{eq:vacancy-mgf} in Section~\ref{subsec:variance} to obtain $L^2$ convergence of the normalized vacancy count, which gives both the expectation asymptotic and the variance bound.


\subsection{Entropy tools}\label{subsec:entropytools}

We collect the entropy inequalities used in the proof; see~\cite[Chapter~2]{CoverThomas} and~\cite[Chapter~4]{BLM}.
All random variables below have finite ranges unless stated otherwise, and all logarithms are natural.

\paragraph{Definitions and basic properties.}
For a random variable $X$ taking values in a finite set $\Omega$, its \emph{entropy} is
\[
    \HH(X)=-\sum_{x\in\Omega}\Prob(X=x)\log\Prob(X=x),
\]
where $0\log0=0$.
The \emph{joint entropy} $\HH(X_1,\ldots,X_m)$ is the entropy of the random vector $(X_1,\ldots,X_m)$.
The \emph{conditional entropy} is
\[
    \HH(X\mid Y)=\sum_{y:\,\Prob(Y=y)>0}\Prob(Y=y)\HH(X\mid Y=y).
\]
When conditioning on auxiliary randomness such as continuous priorities, the sum is replaced by an expectation.
We also write $\HH(P)$ for the entropy of a probability law $P$.

\begin{proposition}[Basic entropy inequalities]\label{prop:entropybasics}
    The following properties hold for finite random variables.
    \begin{enumerate}[(i)]
        \item \textbf{Maximality of the uniform law.} If $X$ takes values in $\Omega$, then
            \[
                \HH(X)\le\log|\Omega|,
            \]
            with equality if and only if $X$ is uniform on $\Omega$.
        \item \textbf{Chain rule.} We have
            \[
                \HH(X_1,\ldots,X_m)=\sum_{i=1}^m\HH(X_i\mid X_1,\ldots,X_{i-1}).
            \]
        \item \textbf{Subadditivity.} We have
            \[
                \HH(X_1,\ldots,X_m)\le\sum_{i=1}^m\HH(X_i).
            \]
        \item \textbf{Conditioning reduces entropy.} We have
            \[
                \HH(X\mid Y,Z)\le\HH(X\mid Y).
            \]
    \end{enumerate}
\end{proposition}

The chain rule and subadditivity also hold conditionally, by averaging over the conditioning variable.
If the conditional range of $X$ given $Y$ is contained in $\Omega_Y$, the support bound and Jensen's inequality give
\begin{equation}\label{eq:conditionalsupport}
    \HH(X\mid Y)\le\E[\log|\Omega_Y|]\le\log\E[|\Omega_Y|].
\end{equation}

The binary entropy function is
\[
    h(p)=-p\log p-(1-p)\log(1-p).
\]
We use the expansion
\[
    h(p)=p\log(1/p)+p+O(p^2)\quad(p\to0),
\]
with the continuous interpretation at zero.

\paragraph{Shearer's inequality.}
We use the following form of Shearer's inequality; see also~\cite{MadimanTetali}.

\begin{lemma}[Shearer's inequality~\cite{CFGS}]\label{lem:shearer}
    Let $X=(X_1,\ldots,X_m)$ be a finite random vector, and let $\mathcal A$ be a family of subsets of $[m]$, allowing repetitions.
    Suppose every coordinate belongs to at least $r\ge1$ members of $\mathcal A$.
    Then
    \[
        \HH(X)\le\frac1r\sum_{A\in\mathcal A}\HH(X_A),
    \]
    where $X_A=(X_i)_{i\in A}$.
\end{lemma}

For a random subset $W$ of the vertex set of a $d$-regular hypergraph $H$, Lemma~\ref{lem:shearer} gives
\[
    \HH(W)\le\frac1d\sum_{e\in E(H)}\HH(W\cap e).
\]

\paragraph{Gibbs variational inequality.}
The following standard inequality relates entropy to weighted counts and follows directly from Jensen's inequality.

\begin{lemma}[Gibbs variational inequality]\label{lem:entropy-counting}
    Let $X$ take values in a finite set $\Omega$, and let $w:\Omega\to(0,\infty)$.
    Then
    \[
        \HH(X)+\E[\log w(X)]\le\log\sum_{x\in\Omega}w(x),
    \]
    with equality if and only if $\Prob(X=x)=w(x)/\sum_{y\in\Omega}w(y)$ for every $x\in\Omega$.
    In particular, for every random matching $J$ of $H$ and every $\lambda>0$,
    \[
        \HH(J)+\E[|J|]\log\lambda\le\log Z_H(\lambda),
    \]
    with equality if and only if $J$ has the Gibbs measure at activity $\lambda$.
\end{lemma}
\begin{proof}[Proof of Lemma~\ref{lem:entropy-counting}]
    Write $q(x)=\Prob(X=x)$.
    By concavity of the logarithm,
    \[
        \begin{aligned}
            \HH(X)+\E[\log w(X)]
            &=\sum_{x:\,q(x)>0}q(x)\log\frac{w(x)}{q(x)}\\
            &\le\log\left(\sum_{x:\,q(x)>0}w(x)\right)
            \le\log\sum_{x\in\Omega}w(x).
        \end{aligned}
    \]
    Equality requires full support and $w(x)/q(x)$ constant, which gives the stated law.
    Taking $\Omega=\calM(H)$ and $w(M)=\lambda^{|M|}$ gives the matching statement.
    This completes the proof.
\end{proof}

Thus a Gibbs matching $M$ satisfies $\HH(M)=\log Z_H(\lambda)-\E[|M|]\log\lambda$.
We apply the inequality to a thinned matching, whose law need not be Gibbs.

\paragraph{Concavity estimates.}
\begin{definition}\label{def:entropy-f}
    For nonnegative $x,y$ with $x+y>0$, define
    \[
        f(x,y)=(x+y)h\left(\frac{x}{x+y}\right),
    \]
    and let $f(0,0)=0$.
\end{definition}
When $x+y\le1$, disjoint events $A,B$ with probabilities $x,y$ satisfy $f(x,y)=\HH(\mathbf1_A\mid\mathbf1_{A\cup B})$.
The next lemma records the concavity of the function in Definition~\ref{def:entropy-f} and a consequence of the log-sum inequality. We use these estimates to bound conditional entropy after averaging over the source matching and the edge deletions.

\begin{lemma}[Concavity estimates]\label{lem:entropyconcavity}
    Let $A,B$ be nonnegative random variables with finite expectations.
    \begin{enumerate}[(i)]
        \item The function $f$ is concave, and
            \[
                \E[f(A,B)]\le f(\E[A],\E[B]).
            \]
        \item If $\E[A]>0$ and $B>0$ whenever $A>0$, then
            \[
                \E\left[A\log\frac BA\right]\le\E[A]\log\frac{\E[B]}{\E[A]},
            \]
            with the integrand defined as zero when $A=0$.
    \end{enumerate}
\end{lemma}


\section{Random greedy counting}\label{sec:greedy}

We develop a variant of random greedy counting that alternates short greedy stages with regular completions.
This gives a weighted counting lower bound over a polynomial range of activities, with an error depending only on the higher codegrees through $\Gamma(H)$, defined in~\eqref{eq:Gamma}.

\begin{theorem}[Counting over a polynomial activity range]\label{thm:greedycount}
    For each fixed $k\ge3$, there is a constant $L=L(k)>0$ such that every simple $n$-vertex $d$-regular $k$-uniform hypergraph satisfying $\Gamma(H)(\log n)^L=o(1)$ obeys
    \[
        \log Z_H(\lambda)\ge\frac nk(\log(\lambda d)-k+1)-O_{k,C_0}\bigl(n\Gamma(H)(\log n)^L\bigr),
    \]
    uniformly for $n^{-C_0}\le \lambda\le n^{C_0}$, for each fixed $C_0>0$.
\end{theorem}

The error bound in Theorem~\ref{thm:greedycount} has no direct dependence on the activity: for each fixed $C_0$, its implicit constant is uniform over $n^{-C_0}\le\lambda\le n^{C_0}$.
This uniformity is essential for the refined weighted count.
When the target activity $\lambda$ lies in a fixed compact interval, simplicity gives $d\le n^{k-1}$, so the theorem also bounds the source remainder at activity $\lambda d$ with the same codegree error.
Thinning these sources in Section~\ref{subsec:countlower} supplies the positive correction $n(\lambda d)^{-1/k}$.
Our goal in this section is to prove Theorem~\ref{thm:greedycount}.

A single large matching is not a sufficient starting point.
For a fixed matching $M_0$, the submatchings obtained by deleting its edges have total weight
\[
    \sum_{J\subseteq M_0}\lambda^{|J|} =(1+\lambda)^{|M_0|}\le(1+\lambda)^{n/k}.
\]
Its logarithm is only $O(n)$ for fixed $\lambda$, whereas the leading term of our count is of order $n\log d$.
We therefore first count a sufficiently rich collection of sources.
Thinning will add the smaller correction after this leading count has been secured.

The random greedy matching process starts with the empty matching.
At each step, it adds a uniformly chosen edge of the residual hypergraph and deletes the chosen vertices and all incident edges, stopping when no edges remain.
A \emph{prefix} is the matching formed by an initial segment of this process, with the order of its edges forgotten.
For classical analyses of this process, see~\cite{SpencerPacking,RodlThoma}.

The difficulty is to keep the loss in the leading count below the correction $nd^{-1/k}$ needed later.
In a single long trajectory, the available degree error in the estimates of Bennett and Bohman~\cite{BB} becomes comparable with the decreasing residual degree too early to give this precision.
For linear $3$-uniform hypergraphs, the resulting cutoff is $d^{-1/4}$ up to logarithmic factors, larger than the relevant vacancy scale $d^{-1/3}$.
Regularization with codegree control has precedents in Kang, K\"uhn, Methuku and Osthus~\cite[Lemma~8.1]{KKMO} and Gould and Kelly~\cite[Lemma~4.1]{GK}.
Gould and Kelly also iterate nibble stages while tracking higher codegrees~\cite[Sections~7--8]{GK}.
The new ingredient here is an iteration of weighted counting inequalities through exactly regular completions.
Each stage stops when half the vertices remain; we control the higher codegrees and restore exact regularity before starting the next stage.
This allows the existing degree estimate to be applied afresh at the current degree scale.
An encoding bounds the cost of transferring weighted counts from the completion to the residual hypergraph.
After subtracting the leading term and dividing by the number of vertices, the next-stage remainder enters with coefficient $1/2+o(1)$, reflecting the proportion of vertices retained in that stage.
The resulting geometric weights, together with the higher-codegree estimates, control the accumulated error on the finer scale required here.
Theorem~\ref{thm:greedycount} may also be of independent interest: it provides an explicit codegree error uniformly up to polynomially large activities, without assumptions on intersection patterns.

The accumulated error is governed by the full codegree parameter $\Gamma(H)$ defined in~\eqref{eq:Gamma}. The $j=3$ term is bounded by the first term. Since $H$ is nonempty and simple, $\Delta_k=1$ and $\Gamma(H)\ge d^{-1/(k-1)}$.

The following construction fills each degree deficit using new vertices, while preserving the higher codegree bounds.
These auxiliary hypergraphs can have a much larger order than $H$, so the completion and short-stage lemmas retain explicit codegree parameters.
The condition $\Gamma(H)(\log n)^L=o(1)$ is imposed only on the original hypergraph.
The partition-function inequality records the cost of the added vertices.

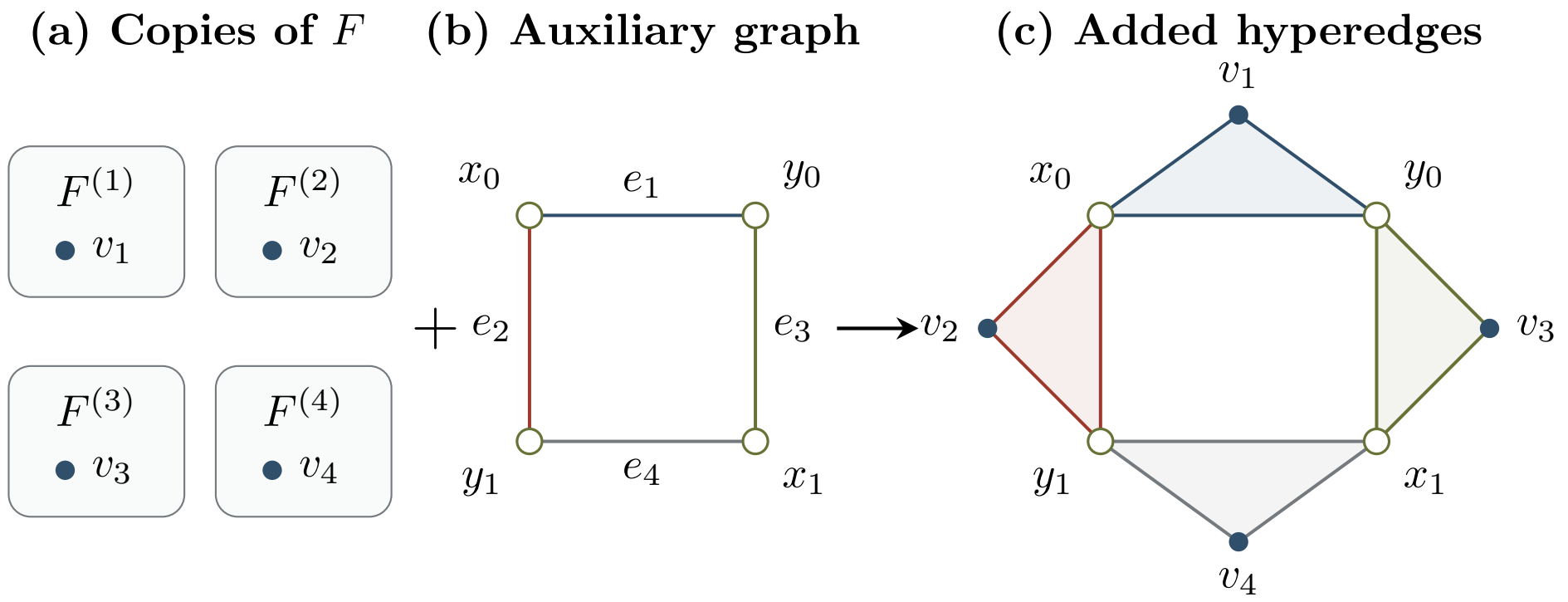
\begin{figure}
    \centering
    \singlespacing
    \matchingfigurestyles
\begin{tikzpicture}[x=1cm,y=1cm]
    \node[panel] at (1.8,4.35) {(a) Copies of $F$};
    \node[panel] at (5.35,4.35) {(b) Auxiliary graph};
    \node[panel] at (10.1,4.35) {(c) Added hyperedges};

    \draw[rounded corners=5pt,draw=MGray,fill=MGray!4] (.3,2.25) rectangle (1.7,3.45);
    \draw[rounded corners=5pt,draw=MGray,fill=MGray!4] (1.95,2.25) rectangle (3.35,3.45);
    \draw[rounded corners=5pt,draw=MGray,fill=MGray!4] (.3,.5) rectangle (1.7,1.7);
    \draw[rounded corners=5pt,draw=MGray,fill=MGray!4] (1.95,.5) rectangle (3.35,1.7);
    \node at (1,3.13) {$F^{(1)}$};
    \node at (2.65,3.13) {$F^{(2)}$};
    \node at (1,1.38) {$F^{(3)}$};
    \node at (2.65,1.38) {$F^{(4)}$};
    \node[vertex,fill=MBlue,label=right:$v_1$] at (.75,2.62) {};
    \node[vertex,fill=MBlue,label=right:$v_2$] at (2.4,2.62) {};
    \node[vertex,fill=MBlue,label=right:$v_3$] at (.75,.87) {};
    \node[vertex,fill=MBlue,label=right:$v_4$] at (2.4,.87) {};
    \node[font=\Large] at (3.7,2) {$+$};

    \draw[draw=MBlue,line width=.75pt] (4.45,2.9) -- node[above] {$e_1$} (6.25,2.9);
    \draw[draw=MRed,line width=.75pt] (4.45,2.9) -- node[left] {$e_2$} (4.45,1.1);
    \draw[draw=MOlive,line width=.75pt] (6.25,2.9) -- node[right] {$e_3$} (6.25,1.1);
    \draw[draw=MGray,line width=.75pt] (4.45,1.1) -- node[below] {$e_4$} (6.25,1.1);
    \node[vacant,draw=MOlive,label=above left:$x_0$] at (4.45,2.9) {};
    \node[vacant,draw=MOlive,label=above right:$y_0$] at (6.25,2.9) {};
    \node[vacant,draw=MOlive,label=below right:$x_1$] at (6.25,1.1) {};
    \node[vacant,draw=MOlive,label=below left:$y_1$] at (4.45,1.1) {};
    \draw[arr] (6.9,2) -- (7.55,2);

    \draw[draw=MBlue,fill=MBlue!8,line width=.75pt] (9,2.9) -- (11.2,2.9) -- (10.1,3.7) -- cycle;
    \draw[draw=MRed,fill=MRed!8,line width=.75pt] (9,2.9) -- (9,1.1) -- (8.1,2) -- cycle;
    \draw[draw=MOlive,fill=MOlive!8,line width=.75pt] (11.2,2.9) -- (11.2,1.1) -- (12.1,2) -- cycle;
    \draw[draw=MGray,fill=MGray!8,line width=.75pt] (9,1.1) -- (11.2,1.1) -- (10.1,.3) -- cycle;
    \node[vacant,draw=MOlive,label=above left:$x_0$] at (9,2.9) {};
    \node[vacant,draw=MOlive,label=above right:$y_0$] at (11.2,2.9) {};
    \node[vacant,draw=MOlive,label=below right:$x_1$] at (11.2,1.1) {};
    \node[vacant,draw=MOlive,label=below left:$y_1$] at (9,1.1) {};
    \node[vertex,fill=MBlue,label=above:$v_1$] at (10.1,3.7) {};
    \node[vertex,fill=MBlue,label=left:$v_2$] at (8.1,2) {};
    \node[vertex,fill=MBlue,label=right:$v_3$] at (12.1,2) {};
    \node[vertex,fill=MBlue,label=below:$v_4$] at (10.1,.3) {};
\end{tikzpicture}
    \caption{Regular completion for one degree deficit at $v$, illustrated with $k=3$ and $D=Q=2$. The four edges $e_i$ of the auxiliary $K_{2,2}$ are extended to $e_i\cup\{v_i\}$, where $v_i$ lies in the $i$th copy of $F$. Each $v_i$ gains one incident edge, and every new vertex retains degree two. A fresh auxiliary graph is used for each remaining deficit.}
    \label{fig:completion}
\end{figure}

\begin{lemma}[Regular completion]\label{lem:completion}
    Let $F$ be a simple $k$-uniform hypergraph on $s$ vertices with maximum degree at most the positive integer $D$, and let $S=\sum_{v\in V(F)}(D-d_F(v))$.
    There is a simple $D$-regular $k$-uniform hypergraph $K$ containing $a=O_k(D^2)$ disjoint copies of $F$ and $w$ additional vertices, where $w / a = (k-1)S/D$.
    One may require $2k\mid |V(K)|$, and $\Delta_j(K)\le\max\{\Delta_j(F),1\}$ for every $j\ge2$.
    Then for every $\lambda>0$,
    \begin{equation*}
        \log Z_F(\lambda)\ge a^{-1}\log Z_K(\lambda) -(w/a)\log(1+D\max\{1,\lambda\}).
    \end{equation*}
\end{lemma}
\begin{proof}[Proof of Lemma~\ref{lem:completion}]
    If $S=0$, take $2k$ disjoint copies of $F$.
    Otherwise choose a prime $Q=O_k(D)$ with $Q\ge\max\{D,k-1\}$.
    Choose $A\subseteq\mathbb F_Q$ of size $D$.
    On $\{0,\ldots,k-2\}\times\mathbb F_Q$, the sets
    \[
        \{(j,b+ja):0\le j<k-1\},
    \]
    indexed by $b\in\mathbb F_Q$ and $a\in A$, form a linear $D$-regular $(k-1)$-uniform hypergraph with $a_0=DQ$ edges.
    Label these edges by $[a_0]$.
    
    Represent the deficits of $F$ by $S$ stubs, and take $a_0$ disjoint copies of $F$.
    For each stub, take a fresh copy of the auxiliary hypergraph.
    Append to its edge $i$ the copy of the stub's vertex in the $i$th copy of $F$.
    Figure~\ref{fig:completion} illustrates this step for $k=3$ and $D=Q=2$.
    Every old and new vertex now has degree $D$.
    The number of new vertices is $(k-1)QS$, and no added edge contains two old vertices.
    Every pair involving a new vertex lies in at most one added edge.
    Taking $2k$ copies gives the divisibility assertion, without changing $w/a=(k-1)S/D$.

    Encode a matching by its restrictions to the old copies and, at each new vertex, its incident matching edge or its absence.
    The encoding is injective.
    Weighting each nonempty new-vertex record by $\max\{1,\lambda\}$ dominates the weight of every added matching edge.
    Hence
    \[
        Z_K(\lambda)\le Z_F(\lambda)^a(1+D\max\{1,\lambda\})^w,
    \]
    which proves the partition-function bound in Lemma~\ref{lem:completion}.
    This completes the proof.
\end{proof}

\subsection{Short greedy trajectories}\label{subsec:shortgreedy}
In the proof of Theorem~\ref{thm:greedycount}, we iteratively use the following theorem. 

\begin{theorem}[Bennett--Bohman~{\cite[Theorem~2.1]{BB}}, fixed-fraction form]\label{thm:BB}
    In an $N$-vertex $D$-regular $k$-uniform hypergraph with maximum pair codegree at most $L_2=o(D/(\log N)^5)$, the random greedy matching process through $N/(2k)$ steps, when this is an integer, has with high probability all surviving degrees
    \begin{equation*}
        Dp_i^{k-1}\pm O_k(\sqrt{DL_2\log N}),
    \end{equation*}
    where $p_i=1-ki/N$.
\end{theorem}
Only this fixed-fraction interval is used.
Its hypothesis requires a sufficiently large power of $\log N$, not a polynomial lower bound in $N$.
We use the following maximal inequality to control the higher codegrees.

\begin{lemma}[Freedman's inequality~\cite{Freedman}]\label{lem:freedman}
    Let $(X_i,\mathcal F_i)_{i\ge0}$ be a real-valued martingale with $X_0=0$.
    Suppose $|X_i-X_{i-1}|\le B$ almost surely for every $i\ge1$, where $B>0$.
    Write
    \[
        V_t=\sum_{i=1}^t\E\bigl[(X_i-X_{i-1})^2\mid\mathcal F_{i-1}\bigr].
    \]
    For all $x,v>0$,
    \[
        \Prob\bigl(|X_t|\ge x\text{ and }V_t\le v\text{ for some }t\ge0\bigr)
        \le2\exp\left(-\frac{x^2}{2(v+Bx)}\right).
    \]
\end{lemma}

\begin{lemma}[Higher codegrees in a short trajectory]\label{lem:short}
    Under the hypotheses of Theorem~\ref{thm:BB}, suppose $\Delta_j\le L_j$ with $L_j\ge1$ for $2\le j\le k$.
    With high probability the residual hypergraph $F$ after $N/(2k)$ steps satisfies the degree estimate in Theorem~\ref{thm:BB} and
    \begin{equation*}
        \Delta_j(F)\le2^{-(k-j)}L_j+ O_k\left(L_j\sqrt{L_2\log N/D}+\sqrt{L_jL_{j+1}\log N}+L_{j+1}\log N\right)
    \end{equation*}
    for every $2\le j<k$.
\end{lemma}
\begin{proof}[Proof of Lemma~\ref{lem:short}]
    Let $G_i$ be the residual hypergraph after $i$ steps, and let $\mathcal F_i$ record the first $i$ choices.
    Write $m=N/(2k)$ and let $\tau$ be the first time that the degree estimate in Theorem~\ref{thm:BB} fails, with $\tau=m$ if it holds throughout.
    Fix $2\le j<k$ and an original $j$-set $S$.
    Let $Y_S(i)$ count original edges $e\supseteq S$ for which all vertices of $e\setminus S$ survive at step $i$.
    We continue to count these edges even if a vertex of $S$ is removed, so removing $S$ does not itself cause a large jump.
    If $S$ survives, $Y_S(i)$ is its codegree in $G_i$.

    \proofparagraph{The one-step drift.}
    Before $\tau$, the degree estimate and $p_i\ge1/2$ give
    \[
        |E(G_i)|=\frac{ND}{k}p_i^k
        \left(1+O_k\left(\sqrt{L_2\log N/D}\right)\right).
    \]
    The next edge $e_{i+1}$ is uniform in $E(G_i)$.
    Hence, for each surviving vertex $u$,
    \[
        \Prob(u\in e_{i+1}\mid\mathcal F_i)
        =\frac{d_{G_i}(u)}{|E(G_i)|}
        =\frac{k}{Np_i}\left(1+O_k\left(\sqrt{L_2\log N/D}\right)\right).
    \]
    For two distinct surviving vertices,
    \[
        \Prob(u,v\in e_{i+1}\mid\mathcal F_i)
        \le\frac{L_2}{|E(G_i)|}
        =O_k\left(\frac{L_2}{ND}\right).
    \]
    An edge $e$ counted by $Y_S(i)$ ceases to be counted precisely when $e_{i+1}$ meets $e\setminus S$.
    The first two terms of inclusion--exclusion give the lower bound
    \[
        \begin{aligned}
            \Prob(e_{i+1}\cap(e\setminus S)\ne\varnothing\mid\mathcal F_i)
            &\ge \sum_{u\in e\setminus S}\Prob(u\in e_{i+1}\mid\mathcal F_i)\\
            &\hspace{8mm}-\sum_{\{u,v\}\in\binom{e\setminus S}{2}}
            \Prob(u,v\in e_{i+1}\mid\mathcal F_i).
        \end{aligned}
    \]
    Thus the removal probability of $e$ at step $i+1$ is at least
    \[
        \frac{k(k-j)}{Np_i}
        -C_k\left(\frac1N\sqrt{L_2\log N/D}+\frac{L_2}{ND}\right).
    \]
    The second error is absorbed by the first under the hypothesis on $L_2$.
    On the other hand, the union bound gives
    \allowdisplaybreaks
    \begin{align*}
        \Prob(e_{i+1}\cap(e\setminus S)\ne\varnothing\mid\mathcal F_i)
        &\le \sum_{u\in e\setminus S}\Prob(u\in e_{i+1}\mid\mathcal F_i)\\
        &\le \frac{k(k-j)}{Np_i}
        +\frac{C_k}{N}\sqrt{L_2\log N/D}.
    \end{align*}
    Thus the removal probability of $e$ is
    \[
        \frac{k(k-j)}{Np_i}
        +O_k\left(\frac1N\sqrt{L_2\log N/D}\right).
    \]
    Summing over the at most $L_j$ edges counted by $Y_S(i)$, and writing $\Delta Y_S(i)=Y_S(i+1)-Y_S(i)$, we obtain
    \begin{equation}\label{eq:codegree-drift}
        \E[\Delta Y_S(i)\mid\mathcal F_i]
        =-\frac{k(k-j)}{Np_i}Y_S(i)+r_i,
    \end{equation}
    where $r_i$ is $\mathcal F_i$-measurable and
    \[
        |r_i|\le C_k\frac{L_j}{N}\sqrt{L_2\log N/D}.
    \]

    \proofparagraph{Removing the linear drift.}
    The main term in~\eqref{eq:codegree-drift} gives a deterministic recurrence that multiplies $Y_S(i)$ by $1-k(k-j)/(Np_i)$ at each step.
    We therefore define
    \[
        g_t=\prod_{i=0}^{t-1}\left(1-\frac{k(k-j)}{Np_i}\right)
    \]
    for $0\le t\le m$, with $g_0=1$.
    Since $p_{i+1}=p_i-k/N$ and $p_i\ge1/2$,
    \[
        1-\frac{k(k-j)}{Np_i}
        =\left(\frac{p_{i+1}}{p_i}\right)^{k-j}
        \left(1+O_k(N^{-2})\right).
    \]
    Multiplying these identities gives
    \begin{equation}\label{eq:greedy-decay}
        g_t=p_t^{k-j}\left(1+O_k(N^{-1})\right).
    \end{equation}
    In particular, $g_t$ is bounded above and below by positive constants depending only on $k$.
    The reason for dividing by $g_t$ is the exact cancellation
    \[
        \frac{Y_S(i+1)}{g_{i+1}}-\frac{Y_S(i)}{g_i}
        =\frac{\Delta Y_S(i)+\frac{k(k-j)}{Np_i}Y_S(i)}{g_{i+1}}.
    \]
    By~\eqref{eq:codegree-drift}, this increment has conditional mean $r_i/g_{i+1}$ before $\tau$.
    Consequently,
    \[
        M_t=\sum_{i<t\wedge\tau}
        \frac{\Delta Y_S(i)-\E[\Delta Y_S(i)\mid\mathcal F_i]}{g_{i+1}}
    \]
    is a martingale, and
    \begin{equation}\label{eq:codegree-decomposition}
        \frac{Y_S(t\wedge\tau)}{g_{t\wedge\tau}}
        =Y_S(0)+\sum_{i<t\wedge\tau}\frac{r_i}{g_{i+1}}+M_t.
    \end{equation}
    The middle sum is the accumulated error in the conditional means.
    Since there are at most $N/(2k)$ summands and $g_i$ is bounded away from zero, uniformly in $t\le m$,
    \begin{equation}\label{eq:codegree-drift-error}
        \left|\sum_{i<t\wedge\tau}\frac{r_i}{g_{i+1}}\right|
        \le C_kL_j\sqrt{L_2\log N/D}.
    \end{equation}

    \proofparagraph{Martingale fluctuations.}
    To control the last term $M_t$ in~\eqref{eq:codegree-decomposition}, we prove the following estimate, which is strong enough for a union bound over all $j$-sets.
    \begin{claim}\label{clm:short-fluctuations}
        For a sufficiently large constant $A=A(k)$, if
        \[
            x=A\left(\sqrt{L_jL_{j+1}\log N}+L_{j+1}\log N\right),
        \]
        then $\Prob(\max_{t\le m}|M_t|\ge x)\le N^{-k-3}$.
    \end{claim}
    \begin{claimproof}[Proof of Claim~\ref{clm:short-fluctuations}]
    Removing a vertex $u\notin S$ destroys at most $L_{j+1}$ edges counted by $Y_S$, since each such edge contains $S\cup\{u\}$.
    At most $k$ vertices are removed in one step, so
    \[
        -kL_{j+1}\le\Delta Y_S(i)\le0.
    \]
    Equation~\eqref{eq:codegree-drift} and the bound on $r_i$ give $\E[-\Delta Y_S(i)\mid\mathcal F_i]=O_k(L_j/N)$ before $\tau$.
    Therefore
    \[
        \begin{aligned}
            \Var[\Delta Y_S(i)\mid\mathcal F_i]
            &\le\E[(\Delta Y_S(i))^2\mid\mathcal F_i]\\
            &\le kL_{j+1}\E[-\Delta Y_S(i)\mid\mathcal F_i]\\
            &=O_k(L_jL_{j+1}/N).
        \end{aligned}
    \]
    It follows that the increments of $M_t$ have absolute value at most $C_kL_{j+1}$, and their total conditional variance satisfies
    \[
        \sum_{i=0}^{m-1}\E[(M_{i+1}-M_i)^2\mid\mathcal F_i]
        \le C_kL_jL_{j+1}.
    \]
    We apply Lemma~\ref{lem:freedman} to the stopped martingale $M_t$, held constant after time $m$.
    The increment bound allows $B=C_kL_{j+1}$, and the conditional variance bound allows $v=C_kL_jL_{j+1}$, with $C_k$ large enough for both bounds above.
    Use the threshold $x$ in the claim, with $A\ge1$ to be fixed below.
    Since the total conditional variance is at most $v$ on every trajectory of the stopped process, the lemma gives
    \[
        \Prob\left(\max_{0\le t\le m}|M_t|\ge x\right)
        \le2\exp\left(-\frac{x^2}{2(v+Bx)}\right).
    \]
    The two terms in the definition of $x$ give
    \[
        \begin{aligned}
            \frac{x^2}{v}&\ge\frac{A^2}{C_k}\log N,\\
            \frac{x}{B}&\ge\frac{A}{C_k}\log N.
        \end{aligned}
    \]
    Consequently, since $A\ge1$,
    \[
        \frac{x^2}{2(v+Bx)}
        \ge\frac14\min\left\{\frac{x^2}{v},\frac{x}{B}\right\}
        \ge\frac{A}{4C_k}\log N.
    \]
    Choosing $A\ge4C_k(k+4)$ therefore yields
    \[
        \Prob\left(\max_{0\le t\le m}|M_t|\ge x\right)
        \le2N^{-k-4}\le N^{-k-3}.
    \]

    \end{claimproof}

    On the event in Theorem~\ref{thm:BB}, we have $\tau=m$ and $p_m=1/2$.
    Use~\eqref{eq:greedy-decay} to estimate $g_m$ in~\eqref{eq:codegree-decomposition}.
    Since $Y_S(0)=d_{G_0}(S)$, the drift error~\eqref{eq:codegree-drift-error} and Claim~\ref{clm:short-fluctuations} then give with high probability
    \[
        Y_S(m)=2^{-(k-j)}d_{G_0}(S)
        +O_k\left(L_j\sqrt{L_2\log N/D}
        +\sqrt{L_jL_{j+1}\log N}+L_{j+1}\log N\right).
    \]
    The error $O_k(L_j/N)$ from approximating $g_m$ is absorbed by the first error term: regularity gives $(k-1)D\le(N-1)L_2$.
    A union bound over all original $j$-sets, for $2\le j<k$, makes these estimates simultaneous with high probability.
    For every surviving $S$, its codegree in $F$ equals $Y_S(m)$, and $d_{G_0}(S)\le L_j$.
    This proves the required codegree bounds together with the degree estimate in Theorem~\ref{thm:BB}.
    This completes the proof.
\end{proof}


\subsection{Counting and iteration}\label{subsec:greedyiteration}

For a regular hypergraph $J$, write
\[
    R_J(\lambda)=\log Z_J(\lambda)-\frac{|V(J)|}{k}\bigl(\log(\lambda\deg(J))-k+1\bigr),
\]
where $\deg(J)$ denotes the degree of $J$.
The next lemma relates the remainder per vertex before and after one greedy stage.

For this lemma, let $G$ be a $D$-regular $k$-uniform hypergraph on $N$ vertices, where $D\ge1$ and $2k\mid N$.
Fix $L_j\ge\max\{1,\Delta_j(G)\}$ for $2\le j\le k$.
For a fixed $C_0>0$, let $\lambda>0$ and $\ell\ge1$ satisfy
\[
    \begin{gathered}
        \log N\le C_0\ell^2,\\
        |\log\lambda|+\log(D+1)\le C_0\ell,\\
        \ell\sqrt{L_2/D}=o(1).
    \end{gathered}
\]
Consider the first $m=N/(2k)$ steps of the random greedy process.
Call a trajectory \emph{good} if the degree estimate in Theorem~\ref{thm:BB} holds throughout these steps and its residual hypergraph satisfies the codegree bounds in Lemma~\ref{lem:short}.
Assume that good trajectories have total probability at least $1/2$, and let $\mathcal P$ be the family of prefixes they produce.
Thus every $P\in\mathcal P$ is an unordered matching of size $m$ admitting at least one good ordering.

Define
\[
    D'=\left\lceil2^{-(k-1)}D+C_1\sqrt{DL_2}\,\ell\right\rceil,
\]
where $C_1=C_1(k,C_0)$ is sufficiently large that every $F_P=G-V(P)$ has maximum degree at most $D'$.
For each $P\in\mathcal P$, choose a $D'$-regular completion $K_P$ from Lemma~\ref{lem:completion}, containing $a_P$ copies of $F_P$ and $w_P$ additional vertices, with $2k\mid |V(K_P)|$.
The codegree bounds in Lemma~\ref{lem:short}, increased to at least one, hold for every $K_P$.

\begin{lemma}[One greedy stage]\label{lem:stage}
    With the notation above,
    \[
        \frac{R_G(\lambda)}N\ge
        \min_{P\in\mathcal P}\left\{
        \left(\frac12+\frac{w_P}{a_PN}\right)
        \frac{R_{K_P}(\lambda)}{|V(K_P)|}\right\}
        -O_{k,C_0}\left(\ell^2\sqrt{L_2/D}\right),
    \]
    where $w_P/(a_PN)=O_{k,C_0}(\ell\sqrt{L_2/D})$ uniformly over $P\in\mathcal P$.
\end{lemma}
\begin{proof}[Proof of Lemma~\ref{lem:stage}]
    Use $p_j=1-kj/N$ as in Theorem~\ref{thm:BB}.
    We count the prefixes and their extensions separately.

    \proofparagraph{Counting good prefixes.}
    \begin{claim}\label{clm:greedy-prefixes}
    \[
        \log(|\mathcal P|\lambda^m) \ge\frac N{2k}\bigl(\log(\lambda D)+(k+1)\log2-k+1\bigr) -O_{k,C_0}(N\ell\sqrt{L_2/D}+\log N).
    \]

    \end{claim}
    \begin{claimproof}[Proof of Claim~\ref{clm:greedy-prefixes}]
    Along a good trajectory, at step $j<m$ the residual hypergraph has $Np_j$ vertices and every surviving degree is at least $Dp_j^{k-1}-O_{k,C_0}(\sqrt{DL_2}\,\ell)$.
    Since $p_j\ge1/2$, the number of available edges is at least
    \[
        \frac{ND}{k}p_j^k(1-O_{k,C_0}(\ell\sqrt{L_2/D})).
    \]
    Each good ordered trajectory has probability at most the reciprocal of the product of these lower bounds.
    Good trajectories have total probability at least $1/2$, and an unordered prefix has at most $m!$ orderings.
    Consequently,
    \begin{equation}\label{eq:greedy-prefix-product}
        |\mathcal P|\ge\frac1{2m!} \prod_{j=0}^{m-1}\frac{ND}{k}p_j^k(1-O_{k,C_0}(\ell\sqrt{L_2/D})).
    \end{equation}
    Stirling's formula and a Riemann sum give
    \begin{align*}
        \log m! & =m\log m-m+O(\log N),\\
        k\sum_{j<m}\log p_j & =N\int_{1/2}^1\log x\,\dd x+O_k(1) =\frac N2(\log2-1)+O_k(1).
    \end{align*}
    Taking logarithms in~\eqref{eq:greedy-prefix-product}, using these two estimates with $m=N/(2k)$, and adding $m\log\lambda$ gives the claimed bound on $\log(|\mathcal P|\lambda^m)$.
    \end{claimproof}

    \proofparagraph{Counting extensions and their multiplicity.}
    For each $P\in\mathcal P$, a matching in $F_P$ extends $P$ to a matching in $G$.
    Therefore
    \[
        \sum_{P\in\mathcal P}\lambda^mZ_{F_P}(\lambda) =\sum_{M\in\mathcal M(G)}\lambda^{|M|} |\{P\in\mathcal P:P\subseteq M\}| \le\binom{N/k}{m}Z_G(\lambda).
    \]
    In particular,
    \begin{equation}\label{eq:greedy-extension-count}
        \log Z_G(\lambda)\ge \log(|\mathcal P|\lambda^m)+\min_{P\in\mathcal P}\log Z_{F_P}(\lambda) -\log\binom{N/k}{m}.
    \end{equation}
    The factors $m!$ and $\binom{N/k}{m}$ remove different multiplicities.
    The former forgets the order of a single prefix; the latter accounts for the different prefixes contained in the same extended matching.

    \proofparagraph{Completion and cancellation of the main terms.}
    Every $F_P$ has $N/2$ vertices and degrees $2^{-(k-1)}D\pm O_{k,C_0}(\sqrt{DL_2}\,\ell)$.
    The choice of $D'$ therefore gives
    \[
        \sum_{v\in V(F_P)}(D'-d_{F_P}(v))=O_{k,C_0}(N\sqrt{DL_2}\,\ell)
    \]
    and
    \begin{equation}\label{eq:completion-vertex-cost}
        w_P/a_P=O_{k,C_0}(N\ell\sqrt{L_2/D}).
    \end{equation}
    Lemma~\ref{lem:completion} and $\log D'=\log D-(k-1)\log2+O_{k,C_0}(\ell\sqrt{L_2/D})$ imply
    \begin{align}
        \log Z_{F_P}(\lambda) & \ge \frac{|V(K_P)|}{ka_P}\bigl(\log(\lambda D')-k+1\bigr)+a_P^{-1}R_{K_P}(\lambda) -(w_P/a_P)\log(1+D'\max\{1,\lambda\})\notag\\
        & \ge\frac N{2k}\bigl(\log(\lambda D)-(k-1)\log2-k+1\bigr) +a_P^{-1}R_{K_P}(\lambda)-O_{k,C_0}(N\ell^2\sqrt{L_2/D}).
    \label{eq:completed-residual-count}
    \end{align}
    Here the contribution of the $w_P/a_P$ added vertices to the first term, and the cost of recording their matching edges, are both $O_{k,C_0}(N\ell^2\sqrt{L_2/D})$.

    The prefix main term plus the residual main term is exactly
    \begin{align*}
        & \frac N{2k}\bigl(\log(\lambda D)+(k+1)\log2-k+1\bigr) +\frac N{2k}\bigl(\log(\lambda D)-(k-1)\log2-k+1\bigr)\\
        & \hspace{35mm}=\frac Nk\bigl(\log(\lambda D)-k+1\bigr)+\frac Nk\log2,
    \end{align*}
    since the coefficient of $\log2$ is $(k+1)-(k-1)=2$.
    Finally,
    \[
        \log\binom{N/k}{N/(2k)}=\frac Nk\log2+O_k(\log N).
    \]
    Thus the binomial multiplicity cancels the remaining $\log2$ term.
    Since $(k-1)D\le(N-1)L_2$, the term $(\log N)/N$ is $O_{k,C_0}(\ell^2L_2/D)$ and is absorbed by $O_{k,C_0}(\ell^2\sqrt{L_2/D})$.
    Substituting the prefix count from Claim~\ref{clm:greedy-prefixes} and the residual count~\eqref{eq:completed-residual-count} into~\eqref{eq:greedy-extension-count} cancels the main terms as above.
    Subtracting $\frac Nk(\log(\lambda D)-k+1)$ and dividing by $N$ gives the stated inequality, since
    \[
        \frac{R_{K_P}(\lambda)}{a_PN}
        =\left(\frac12+\frac{w_P}{a_PN}\right)
        \frac{R_{K_P}(\lambda)}{|V(K_P)|}.
    \]
    Here $|V(K_P)|=a_PN/2+w_P$, and~\eqref{eq:completion-vertex-cost} bounds the correction to the coefficient $1/2$ uniformly over $P\in\mathcal P$.
    This completes the proof.
\end{proof}

The hypergraph $K_P$ is not a residual subhypergraph of the original hypergraph.
It contains several copies of $G-V(P)$, together with new vertices and edges that fill the degree deficits.
Its order may be larger than $N$.
The coefficient of its normalized remainder in Lemma~\ref{lem:stage} is approximately $1/2$.
This distinction is essential when the inequality is iterated.

\subsection{Proof of Theorem~\ref{thm:greedycount}}\label{subsec:proof-greedycount}
We now iterate the one-stage inequality, keeping track of its error in the original normalization.

\begin{proof}[Proof of Theorem~\ref{thm:greedycount}]
    We may assume that $2k\mid n$.
    Otherwise, apply the argument to the disjoint union of $2k$ copies of $H$ and divide the resulting logarithmic bound by $2k$.
    Degrees and codegrees are unchanged, partition functions multiply, and the logarithm of the order changes by only a bounded amount.
    Let $\ell=\log(n+2)$, where $n$ denotes the order of $H$, not that of a later completion.
    Simplicity gives
    \[
        d\le\binom{n-1}{k-1}\le n^{k-1},
    \]
    and hence $\log d=O_k(\ell)$.
    No polynomial lower bound on $d$ is required here; we retain the hypothesis $\Gamma(H)(\log n)^L=o(1)$.

    \proofparagraph{The stages and the stopping scale.}
    We will construct regular hypergraphs $H_i$, starting with $H_0=H$.
    Write $N_i=|V(H_i)|$ and let $D_i$ be their degrees, with $N_0=n$ and $D_0=d$.
    The deterministic quantity $L_j^{(i)}$ will bound the $j$-codegree at stage $i$.
    Initially $L_j^{(0)}=\Delta_j(H)$.
    Define $p_i=2^{-i}$ and, for sufficiently large constants depending only on $k$, set
    \begin{equation}\label{eq:greedy-degree-recurrence}
        D_{i+1}=\left\lceil 2^{-(k-1)}D_i+C_k\sqrt{D_iL_2^{(i)}}\,\ell\right\rceil
    \end{equation}
    and
    \begin{equation}\label{eq:greedy-codegree-recurrence}
        L_j^{(i+1)}=(2^{-(k-j)}+C_k\ell^{-2})L_j^{(i)}+C_k\ell^4L_{j+1}^{(i)}
    \end{equation}
    for $2\le j<k$, while $L_k^{(i)}=1$.
    These sequences are defined before the prefixes or completions are chosen.

    Choose a sufficiently large fixed $B=B(k)$ and let
    \[
        T=\left\lfloor\log_2\frac1{\Gamma(H)\ell^B}\right\rfloor.
    \]
    Taking $L>B$ ensures that $\Gamma(H)\ell^B\to0$, and
    \begin{equation}\label{eq:greedy-stopping-scale}
        \Gamma(H)\ell^B\le p_T<2\Gamma(H)\ell^B.
    \end{equation}
    Since $H$ is nonempty and simple, $\Delta_k(H)=1$, so
    \[
        \Gamma(H)\ge d^{-1/(k-1)}\ge n^{-1}.
    \]
    Thus $T\le\log_2 n=O_k(\ell)$.
    We first estimate the deterministic sequences and then verify that they bound every completed stage up to $T$.

    \proofparagraph{Degree and codegree estimates.}
    The next claim describes the decay of each codegree and the accumulated increase in the regular degree caused by completion.

    \begin{claim}\label{clm:greedy-sequences}
        For $0\le i\le T$ and $2\le j\le k$,
        \[
            L_j^{(i)}\le C_k\ell^{C_k}\sum_{b=j}^k\Delta_b(H)p_i^{k-b}.
        \]
        Also $D_i\ge dp_i^{k-1}$ and
        \[
            D_i=dp_i^{k-1}\exp\left(O_k\left(\ell\sum_{h<i}\sqrt{L_2^{(h)}/D_h}\right)\right).
        \]
    \end{claim}
    \begin{claimproof}[Proof of Claim~\ref{clm:greedy-sequences}]
        To iterate the codegree recurrence, write $q_j=2^{-(k-j)}+C_k\ell^{-2}$ for $j<k$ and $q_k=1$.
        Unrolling~\eqref{eq:greedy-codegree-recurrence} gives
        \[
            L_j^{(i)}\le q_j^i\Delta_j(H)
            +C_k\ell^4\sum_{h=0}^{i-1}q_j^{i-1-h}L_{j+1}^{(h)}.
        \]
        For $j<b$, the difference $q_b-q_j$ is bounded below by a positive constant depending only on $k$, for all sufficiently large $n$.
        Therefore
        \[
            \sum_{h=0}^{i-1}q_j^{i-1-h}q_b^h
            =\frac{q_b^i-q_j^i}{q_b-q_j}\le C_kq_b^i.
        \]
        Descending induction on $j$, starting with $L_k^{(i)}=1$, now gives
        \[
            L_j^{(i)}\le C_k\sum_{b=j}^k\ell^{4(b-j)}\Delta_b(H)q_b^i.
        \]
        Since $i=O_k(\ell)$, we have $q_b^i=p_i^{k-b}\exp(O_k(i\ell^{-2}))=O_k(p_i^{k-b})$ for $b<k$, and $q_k^i=1$.
        There are at most $k-2$ transitions between indices, so the powers of $\ell$ are bounded in terms of $k$ alone.
        This proves the codegree estimate.

        The degree recurrence directly gives $D_i\ge dp_i^{k-1}$.
        Since $D_i,L_2^{(i)},\ell\ge1$, the rounding error in~\eqref{eq:greedy-degree-recurrence} is absorbed by $\ell\sqrt{D_iL_2^{(i)}}$.
        Consequently,
        \[
            0\le\log\frac{D_{i+1}}{2^{-(k-1)}D_i}
            \le C_k\ell\sqrt{L_2^{(i)}/D_i}.
        \]
        Summing over $i$ and exponentiating proves the degree estimate.
    \end{claimproof}

    \begin{claim}\label{clm:greedy-small-errors}
        If $B$ is sufficiently large, then for every $i\le T$,
        \[
            \ell\sqrt{L_2^{(i)}/D_i}
            \le C_k\ell^{C_k}\sum_{j=2}^k
            \sqrt{\Delta_j(H)/d}\,p_i^{-(j-1)/2}
            \le\ell^{-6}.
        \]
        In particular,
        \[
            \ell\sum_{i<T}\sqrt{L_2^{(i)}/D_i}=o(1)
        \]
        and $D_i=(1+o(1))dp_i^{k-1}\le2d$ uniformly for $i\le T$.
    \end{claim}
    \begin{claimproof}[Proof of Claim~\ref{clm:greedy-small-errors}]
        Claim~\ref{clm:greedy-sequences} gives
        \[
            \frac{L_2^{(i)}}{D_i}
            \le C_k\ell^{C_k}\sum_{j=2}^k
            \frac{\Delta_j(H)}d\,p_i^{-(j-1)}.
        \]
        Taking square roots and using $\sqrt{\sum_j x_j}\le\sum_j\sqrt{x_j}$ proves the first inequality.
        The definition of $\Gamma(H)$ gives $\sqrt{\Delta_2(H)/d}\le\Gamma(H)$ and, for $j\ge3$, $\Delta_j(H)/d\le\Gamma(H)^{j-1}$.
        Since $p_i\ge p_T\ge\Gamma(H)\ell^B$,
        \[
            \sqrt{\Delta_2(H)/d}\,p_i^{-1/2}
            \le\Gamma(H)^{1/2}\ell^{-B/2}\le\ell^{-B/2},
        \]
        whereas, for $j\ge3$,
        \[
            \sqrt{\Delta_j(H)/d}\,p_i^{-(j-1)/2}
            \le\ell^{-B(j-1)/2}.
        \]
        Increasing $B$ makes their sum, including the fixed logarithmic factor, at most $\ell^{-6}$.
        Summing over $T=O_k(\ell)$ stages gives $O_k(\ell^{-5})=o(1)$, and the degree assertion follows from Claim~\ref{clm:greedy-sequences}.
    \end{claimproof}

    \proofparagraph{Validity of the completed stages.}
    We must check the short-trajectory hypotheses at the order of each auxiliary hypergraph.
    Although this order depends on the prefix, the following bound is uniform over the choices.

    \begin{claim}\label{clm:greedy-orders}
        Every chain of completions obtained from Lemma~\ref{lem:completion}, with degrees $D_i$ and with each residual on half the current vertices, satisfies
        \[
            \log N_i=O_k(\ell^2)
        \]
        for $i\le T$.
    \end{claim}
    \begin{claimproof}[Proof of Claim~\ref{clm:greedy-orders}]
        At stage $i$, the residual has $N_i/2$ vertices and its total deficit relative to $D_{i+1}$ is at most $D_{i+1}N_i/2$.
        If the completion uses $a_i$ copies and $w_i$ new vertices, Lemma~\ref{lem:completion} gives
        \[
            a_i=O_k(D_{i+1}^2)
        \]
        and $w_i/a_i\le(k-1)N_i/2$.
        Hence, by Claim~\ref{clm:greedy-small-errors},
        \[
            N_{i+1}=a_iN_i/2+w_i\le C_kd^2N_i.
        \]
        Iterating and taking logarithms yields
        \[
            \log N_i\le\log n+i\log(C_kd^2)=O_k(\ell^2),
        \]
        using $i=O_k(\ell)$ and $\log d=O_k(\ell)$.
    \end{claimproof}

    Suppose that $H_i$ has been constructed with the stated degree and codegree bounds.
    Claims~\ref{clm:greedy-small-errors} and~\ref{clm:greedy-orders} imply
    \[
        \frac{L_2^{(i)}(\log N_i)^5}{D_i}
        =O_k(\ell^{-14}\ell^{10})=o(1).
    \]
    Theorem~\ref{thm:BB} and Lemma~\ref{lem:short} therefore give good trajectories with probability at least $1/2$, uniformly over the possible $H_i$.
    Their residual degrees are
    \[
        2^{-(k-1)}D_i\pm O_k(\ell\sqrt{D_iL_2^{(i)}}).
    \]
    Choosing the constant in~\eqref{eq:greedy-degree-recurrence} sufficiently large makes $D_{i+1}$ at least every such residual degree.
    Lemma~\ref{lem:completion} fills the deficits to this prescribed integer, so $D_{i+1}$ is the same for all good completions at stage $i$.

    For the higher codegrees, Lemma~\ref{lem:short} and Claim~\ref{clm:greedy-orders} bound the three error terms by a constant multiple of
    \[
        \ell^{-6}L_j^{(i)}
        +\ell\sqrt{L_j^{(i)}L_{j+1}^{(i)}}
        +\ell^2L_{j+1}^{(i)}.
    \]
    Young's inequality gives
    \[
        \ell\sqrt{L_j^{(i)}L_{j+1}^{(i)}}
        \le\frac12\ell^{-2}L_j^{(i)}
        +\frac12\ell^4L_{j+1}^{(i)}.
    \]
    Thus~\eqref{eq:greedy-codegree-recurrence} bounds the residual codegrees, for a sufficiently large constant.
    It also bounds the completed codegrees, since $L_j^{(i+1)}\ge1$ and Lemma~\ref{lem:completion} preserves every upper bound of at least one on the maximum $j$-codegree.
    Starting with $H_0=H$, this verifies every stage through $T$.
    Moreover, $|\log\lambda|+\log(D_i+1)=O_{k,C_0}(\ell)$, so all hypotheses of Lemma~\ref{lem:stage} hold at every stage before $T$.

    \proofparagraph{Iteration of the normalized remainder.}
    At each stage choose a prefix attaining the minimum in Lemma~\ref{lem:stage}, and let $H_{i+1}$ be its completion.
    Write $\rho_i=R_{H_i}(\lambda)/N_i$.
    The next claim keeps the coefficients explicit; it is these coefficients that prevent the large auxiliary orders from enlarging the original error.

    \begin{claim}\label{clm:greedy-remainder-iteration}
        For $i<T$,
        \[
            \rho_i\ge-C_{k,C_0}\ell^2\sqrt{L_2^{(i)}/D_i}
            +\frac{N_{i+1}}{a_iN_i}\rho_{i+1}.
        \]
        For every $i\le T$, the accumulated coefficient is
        \[
            \begin{aligned}
                \prod_{h=0}^{i-1}\frac{N_{h+1}}{a_hN_h}
                &=\frac{N_i}{n\prod_{h=0}^{i-1}a_h}\\
                &=p_i\exp\left(O_{k,C_0}\left(
                \ell\sum_{h<i}\sqrt{L_2^{(h)}/D_h}\right)\right)
                =O_{k,C_0}(p_i).
            \end{aligned}
        \]
        Consequently,
        \[
            \begin{aligned}
                \rho_0\ge{}&-C_{k,C_0}\ell^2\sum_{i=0}^{T-1}
                \left(\prod_{h=0}^{i-1}\frac{N_{h+1}}{a_hN_h}\right)
                \sqrt{L_2^{(i)}/D_i}\\
                &+\left(\prod_{h=0}^{T-1}\frac{N_{h+1}}{a_hN_h}\right)\rho_T.
            \end{aligned}
        \]
        An empty product is interpreted as one.
    \end{claim}
    \begin{claimproof}[Proof of Claim~\ref{clm:greedy-remainder-iteration}]
        Since $N_{i+1}=a_iN_i/2+w_i$, Lemma~\ref{lem:stage} gives the first inequality and
        \[
            \frac{N_{i+1}}{a_iN_i}
            =\frac12+\frac{w_i}{a_iN_i}
            =\frac12\left(1+O_{k,C_0}(\ell\sqrt{L_2^{(i)}/D_i})\right).
        \]
        In the product, the intermediate orders $N_1,\ldots,N_{i-1}$ cancel.
        The remaining numerator is $N_i$, and the denominator is $n\prod_{h<i}a_h$.
        Taking logarithms of the ratios and summing gives the exponential expression.
        Claim~\ref{clm:greedy-small-errors} bounds its exponent by $o(1)$, uniformly in $i\le T$.
        Finally, all ratios are positive, so repeated substitution of the first inequality gives the asserted lower bound for $\rho_0$.
    \end{claimproof}

    By Claim~\ref{clm:greedy-remainder-iteration}, the stage-$i$ error contributes at most
    \[
        O_{k,C_0}(np_i\ell^2\sqrt{L_2^{(i)}/D_i})
    \]
    to the original logarithmic count.
    The factor is $np_i$, irrespective of the size of $N_i$.

    \proofparagraph{The final error bound.}
    The empty matching gives $Z_{H_T}(\lambda)\ge1$.
    Since Claim~\ref{clm:greedy-small-errors} gives $D_T\le2d$ and $|\log\lambda|=O_{C_0}(\ell)$,
    \[
        \rho_T\ge-\frac1k\bigl(\log(\lambda D_T)-k+1\bigr)
        \ge-O_{k,C_0}(\ell).
    \]
    Apply Claim~\ref{clm:greedy-remainder-iteration} and multiply by $n$ to obtain
    \begin{equation}\label{eq:greedy-total-error}
        R_H(\lambda)\ge
        -C_{k,C_0}n\ell^2\sum_{i<T}p_i\sqrt{L_2^{(i)}/D_i}
        -C_{k,C_0}np_T\ell.
    \end{equation}
    It remains to bound the sum on the right of~\eqref{eq:greedy-total-error}. Claim~\ref{clm:greedy-sequences}, with $j=2$, and $D_i\ge dp_i^{k-1}$ give
    \[
        p_i\sqrt{L_2^{(i)}/D_i}
        \le C_k\ell^{C_k}\sum_{j=2}^k
        \sqrt{\Delta_j(H)/d}\,p_i^{(3-j)/2}.
    \]
    Here $\sqrt{\sum_j x_j}\le\sum_j\sqrt{x_j}$, and the extra factor $p_i$ changes the exponent from $-(j-1)/2$ to $(3-j)/2$.
    Exchanging the two finite sums in~\eqref{eq:greedy-total-error}, the stage errors are therefore at most
    \[
        C_{k,C_0}n\ell^{C_k}\sum_{j=2}^k
        \sqrt{\Delta_j(H)/d}\sum_{i<T}p_i^{(3-j)/2}.
    \]
    For $j=2$,
    \[
        \sqrt{\Delta_2(H)/d}\sum_{i<T}p_i^{1/2}
        \le\frac{\Gamma(H)}{1-2^{-1/2}}.
    \]
    For $j=3$, the inner sum is $T$, so its contribution is at most $\Gamma(H)T=O_k(\Gamma(H)\ell)$.
    For $j\ge4$, the geometric sum satisfies
    \[
        \sum_{i<T}p_i^{(3-j)/2}
        =\frac{2^{T(j-3)/2}-1}{2^{(j-3)/2}-1}
        \le C_kp_T^{(3-j)/2}.
    \]
    Using~\eqref{eq:greedy-stopping-scale} and $\Delta_j(H)/d\le\Gamma(H)^{j-1}$, we obtain
    \[
        \begin{aligned}
            \sqrt{\Delta_j(H)/d}\,p_T^{(3-j)/2}
            &\le\Gamma(H)^{(j-1)/2}
            (\Gamma(H)\ell^B)^{(3-j)/2}\\
            &=\Gamma(H)\ell^{-B(j-3)/2}\le\Gamma(H).
        \end{aligned}
    \]
    The terminal contribution satisfies $np_T\ell=O_k(n\Gamma(H)\ell^{B+1})$.
    Substituting the stage and terminal estimates into~\eqref{eq:greedy-total-error}, and increasing the fixed $L=L(k)$ to absorb their logarithmic powers, gives
    \[
        R_H(\lambda)\ge-O_{k,C_0}\bigl(n\Gamma(H)(\log n)^L\bigr).
    \]
    Recalling the definition of the remainder,
    \[
        R_H(\lambda)=\log Z_H(\lambda)
        -\frac nk\bigl(\log(\lambda d)-k+1\bigr),
    \]
    we conclude that
    \[
        \log Z_H(\lambda)\ge\frac nk\bigl(\log(\lambda d)-k+1\bigr)
        -O_{k,C_0}\bigl(n\Gamma(H)(\log n)^L\bigr),
    \]
    uniformly over the stated activity range, as required.
    This completes the proof.
\end{proof}


\section{Intersection patterns}\label{sec:patterns}

We now develop the entropy bounds that convert the count in Theorem~\ref{thm:greedycount} into the correction in Theorem~\ref{thm:weightedcount}.
Throughout this section, $H$ satisfies the hypotheses of Theorem~\ref{thm:weightedcount}.
All $o(1)$ terms tend to zero uniformly over these hypergraphs as $n\to\infty$, with $k,\eps,C$ fixed; the estimates for random matchings are also uniform over their laws.
For the lower bound, delete each edge of a random source matching $M$ independently with probability $\theta$, leaving $J$.
The entropy chain rule gives
\[
    \HH(J)=\HH(M)+\E[|M|]h(\theta)-\HH(M\mid J).
\]
The term $\HH(M\mid J)$ accounts for different sources producing the same $J$; bounding it allows us to use the entropy gained by deletion.

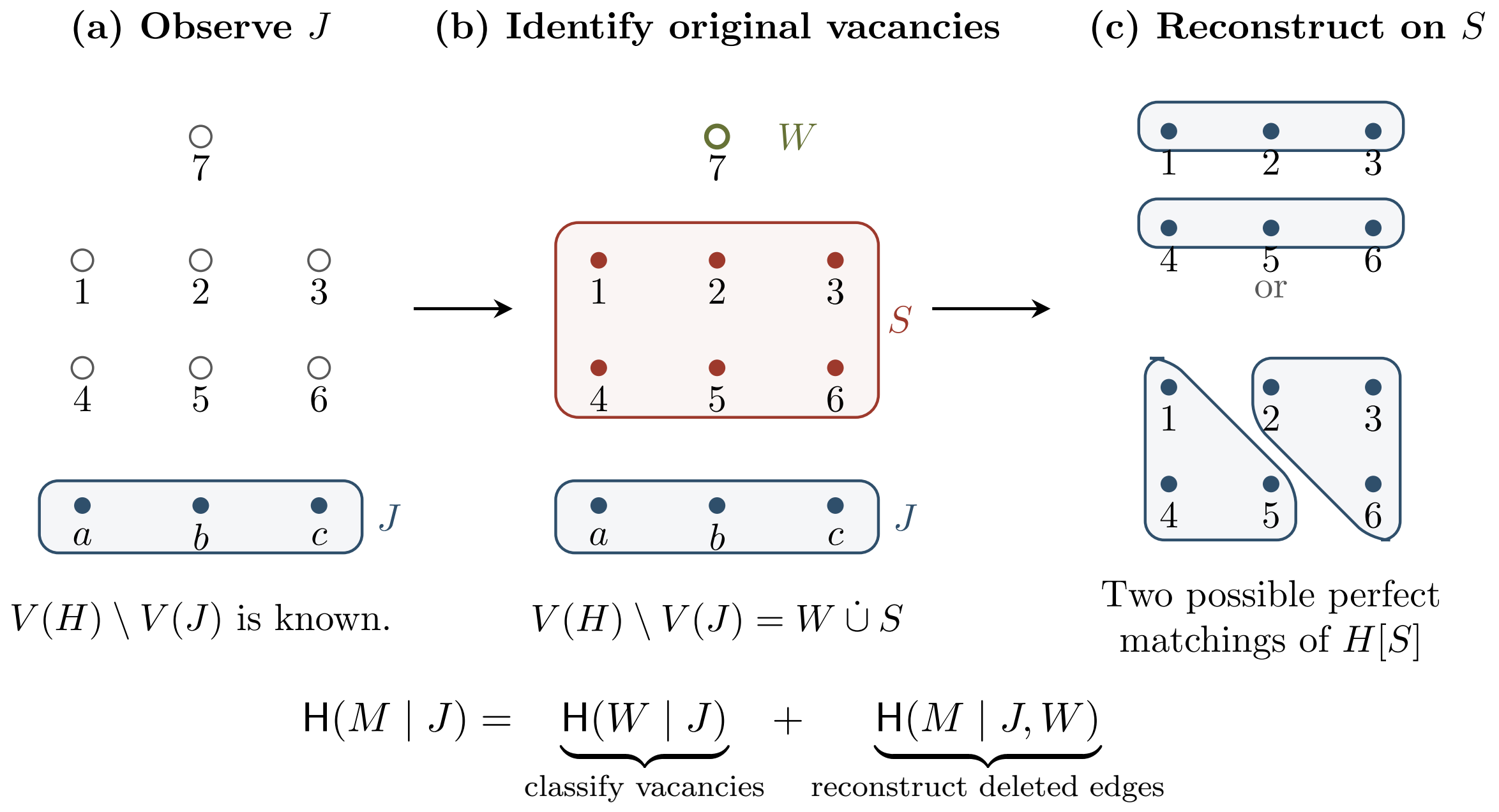
\begin{figure}
    \centering
    \singlespacing
    \matchingfigurestyles
    \resizebox{\linewidth}{!}{%
\begin{tikzpicture}[x=1cm,y=1cm]
\node[panel] at (1.600,4.150) {(a) Observe $J$};
\node[panel] at (6.400,4.150) {(b) Identify original vacancies};
\node[panel] at (11.700,4.150) {(c) Reconstruct on $S$};
\node[vacant] at (0.500,2.000) {};
\node[] at (0.500,1.710) {1};
\node[vacant] at (1.600,2.000) {};
\node[] at (1.600,1.710) {2};
\node[vacant] at (2.700,2.000) {};
\node[] at (2.700,1.710) {3};
\node[vacant] at (0.500,1.000) {};
\node[] at (0.500,0.710) {4};
\node[vacant] at (1.600,1.000) {};
\node[] at (1.600,0.710) {5};
\node[vacant] at (2.700,1.000) {};
\node[] at (2.700,0.710) {6};
\node[vacant] at (1.600,3.150) {};
\node[] at (1.600,2.860) {7};
\draw[source] (0.10,-.72) rectangle (3.10,-.05);
\node[vertex,fill=MBlue] at (0.500,-0.280) {};
\node[] at (0.500,-0.570) {$a$};
\node[vertex,fill=MBlue] at (1.600,-0.280) {};
\node[] at (1.600,-0.570) {$b$};
\node[vertex,fill=MBlue] at (2.700,-0.280) {};
\node[] at (2.700,-0.570) {$c$};
\node[text=MBlue] at (3.350,-0.390) {$J$};
\draw[draw=MRed,fill=MRed!5,line width=.75pt,rounded corners=6pt] (4.90,.54) rectangle (7.90,2.35);
\node[text=MRed] at (8.100,1.450) {$S$};
\node[vertex,fill=MRed] at (5.300,2.000) {};
\node[] at (5.300,1.710) {1};
\node[vertex,fill=MRed] at (6.400,2.000) {};
\node[] at (6.400,1.710) {2};
\node[vertex,fill=MRed] at (7.500,2.000) {};
\node[] at (7.500,1.710) {3};
\node[vertex,fill=MRed] at (5.300,1.000) {};
\node[] at (5.300,0.710) {4};
\node[vertex,fill=MRed] at (6.400,1.000) {};
\node[] at (6.400,0.710) {5};
\node[vertex,fill=MRed] at (7.500,1.000) {};
\node[] at (7.500,0.710) {6};
\node[vacant,draw=MOlive,line width=1.2pt] at (6.400,3.150) {};
\node[] at (6.400,2.860) {7};
\node[text=MOlive] at (7.150,3.150) {$W$};
\draw[source] (4.90,-.72) rectangle (7.90,-.05);
\node[vertex,fill=MBlue] at (5.300,-0.280) {};
\node[] at (5.300,-0.570) {$a$};
\node[vertex,fill=MBlue] at (6.400,-0.280) {};
\node[] at (6.400,-0.570) {$b$};
\node[vertex,fill=MBlue] at (7.500,-0.280) {};
\node[] at (7.500,-0.570) {$c$};
\node[text=MBlue] at (8.150,-0.390) {$J$};
\draw[arr] (3.58,1.55) -- (4.50,1.55);
\draw[arr] (8.40,1.55) -- (9.50,1.55);
\draw[source] (10.32,2.12) rectangle (12.78,2.57);
\draw[source] (10.32,3.02) rectangle (12.78,3.47);
\node[vertex,fill=MBlue] at (10.600,3.200) {};
\node[] at (10.600,2.910) {1};
\node[vertex,fill=MBlue] at (11.550,3.200) {};
\node[] at (11.550,2.910) {2};
\node[vertex,fill=MBlue] at (12.500,3.200) {};
\node[] at (12.500,2.910) {3};
\node[vertex,fill=MBlue] at (10.600,2.300) {};
\node[] at (10.600,2.010) {4};
\node[vertex,fill=MBlue] at (11.550,2.300) {};
\node[] at (11.550,2.010) {5};
\node[vertex,fill=MBlue] at (12.500,2.300) {};
\node[] at (12.500,2.010) {6};
\node[text=black!65] at (11.550,1.730) {or};
\draw[source] (10.38,1.09) -- (10.60,1.09) -- (11.78,-0.09) -- (11.78,-0.60) -- (10.38,-0.60) -- cycle;
\draw[source] (11.38,1.09) -- (12.75,1.09) -- (12.75,-0.60) -- (12.48,-0.60) -- (11.38,0.50) -- cycle;
\node[vertex,fill=MBlue] at (10.600,0.820) {};
\node[] at (10.600,0.530) {1};
\node[vertex,fill=MBlue] at (11.550,0.820) {};
\node[] at (11.550,0.530) {2};
\node[vertex,fill=MBlue] at (12.500,0.820) {};
\node[] at (12.500,0.530) {3};
\node[vertex,fill=MBlue] at (10.600,-0.080) {};
\node[] at (10.600,-0.370) {4};
\node[vertex,fill=MBlue] at (11.550,-0.080) {};
\node[] at (11.550,-0.370) {5};
\node[vertex,fill=MBlue] at (12.500,-0.080) {};
\node[] at (12.500,-0.370) {6};
\node[note] at (1.600,-1.350) {$V(H)\setminus V(J)$ is known.};
\node[note] at (6.400,-1.350) {$V(H)\setminus V(J)=W\mathbin{\dot\cup}S$};
\node[note,text width=4.6cm] at (11.550,-1.350) {Two possible perfect matchings of $H[S]$};
\node[font=\normalsize] at (6.550,-2.580) {$\mathsf H(M\mid J)=\underbrace{\mathsf H(W\mid J)}_{\text{classify vacancies}}+\underbrace{\mathsf H(M\mid J,W)}_{\text{reconstruct deleted edges}}$};
\end{tikzpicture}%
}
    \caption{The two reconstruction problems. The matching $J$ does not distinguish the original vacancy $W=\{7\}$ from the six vertices in $S$ made vacant by deletion. Once $W$ is specified, reconstruction amounts to choosing a perfect matching of $H[S]$; two possible choices are shown. The retained edge on $a,b,c$ is unchanged.}
    \label{fig:reconstruction}
\end{figure}

There are two reconstruction problems.
First, specify the set $W$ of vertices not covered by $J$ which were already uncovered in the source.
The other vertices form the deleted set $S$.
Once $S$ is known, reconstructing the source amounts to choosing a perfect matching on $H[S]$.
Intersection patterns control the choices in this second problem.
They record which whole source blocks must be deleted for an alternative edge to lie in $S$.
The vertices of a source edge are deleted together, so a vertex-by-vertex independence estimate would miss these correlations.
Figure~\ref{fig:reconstruction} illustrates the two reconstruction problems.

For the upper count, we first specify the uncovered set and then expose a perfect matching of the remaining vertices.
An alternative is available only if the matching blocks it meets have not yet been exposed.
The same pattern estimates therefore bound the entropy in both arguments.


\subsection{Coefficients and an integral estimate}\label{subsec:patterncoefficients}

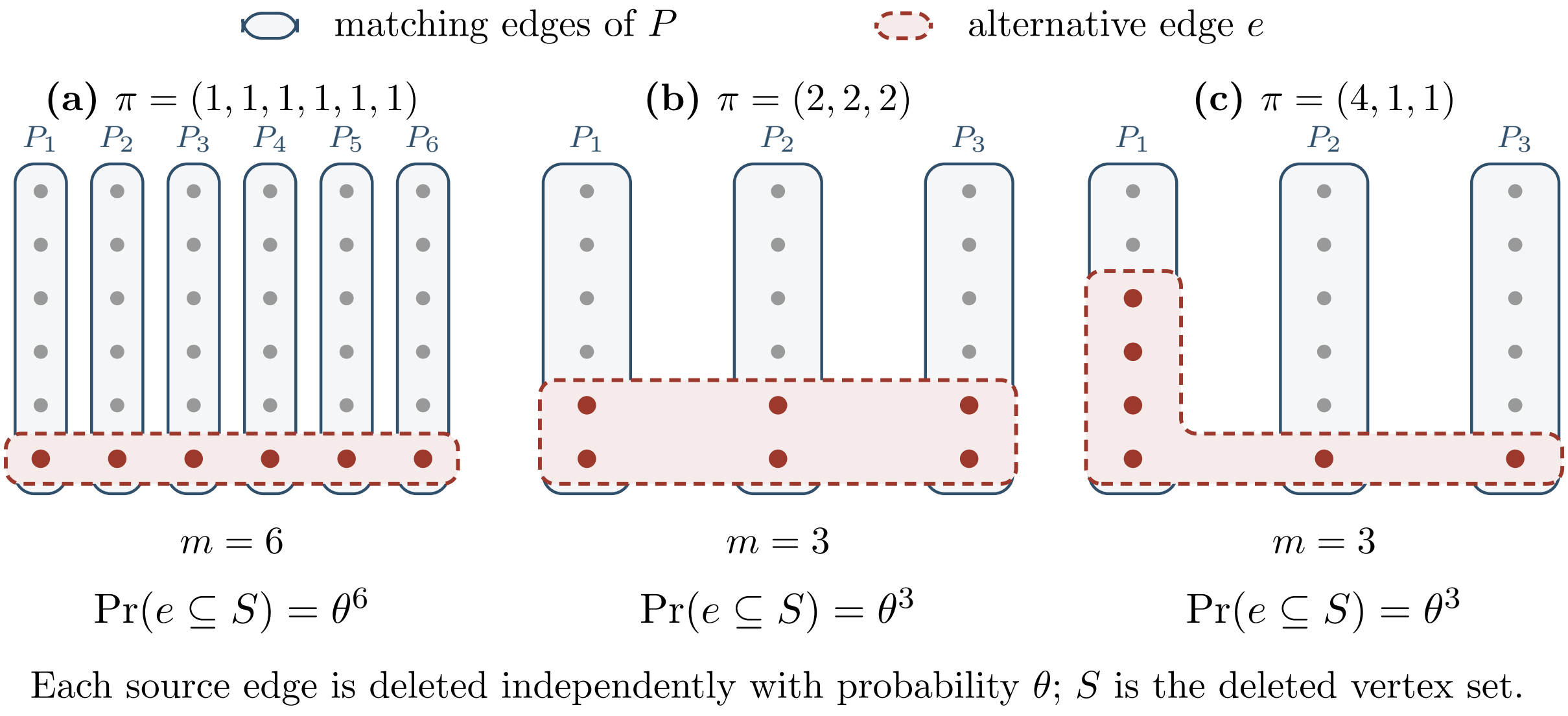
\begin{figure}
    \centering
    \singlespacing
    \matchingfigurestyles
    \resizebox{\linewidth}{!}{%
\begin{tikzpicture}[x=1cm,y=1cm]
\draw[source] (2.2,4.65) rectangle (2.7,4.88);
\node[anchor=west] at (2.9,4.76) {matching edges of $P$};
\draw[alternative] (8.0,4.65) rectangle (8.5,4.88);
\node[anchor=west] at (8.7,4.76) {alternative edge $e$};
\node[panel] at (2.100,4.080) {(a) $\pi=(1,1,1,1,1,1)$};
\draw[source] (0.115,.48) rectangle (0.585,3.50);
\node[text=MBlue,font=\scriptsize] at (0.350,3.730) {$P_{1}$};
\draw[source] (0.815,.48) rectangle (1.285,3.50);
\node[text=MBlue,font=\scriptsize] at (1.050,3.730) {$P_{2}$};
\draw[source] (1.515,.48) rectangle (1.985,3.50);
\node[text=MBlue,font=\scriptsize] at (1.750,3.730) {$P_{3}$};
\draw[source] (2.215,.48) rectangle (2.685,3.50);
\node[text=MBlue,font=\scriptsize] at (2.450,3.730) {$P_{4}$};
\draw[source] (2.915,.48) rectangle (3.385,3.50);
\node[text=MBlue,font=\scriptsize] at (3.150,3.730) {$P_{5}$};
\draw[source] (3.615,.48) rectangle (4.085,3.50);
\node[text=MBlue,font=\scriptsize] at (3.850,3.730) {$P_{6}$};
\draw[alternative] (0.030,.57) rectangle (4.170,1.03);
\node[circle,inner sep=1.7pt,fill=MRed] at (0.350,0.800) {};
\node[circle,inner sep=1.3pt,fill=black!40] at (0.350,1.290) {};
\node[circle,inner sep=1.3pt,fill=black!40] at (0.350,1.780) {};
\node[circle,inner sep=1.3pt,fill=black!40] at (0.350,2.270) {};
\node[circle,inner sep=1.3pt,fill=black!40] at (0.350,2.760) {};
\node[circle,inner sep=1.3pt,fill=black!40] at (0.350,3.250) {};
\node[circle,inner sep=1.7pt,fill=MRed] at (1.050,0.800) {};
\node[circle,inner sep=1.3pt,fill=black!40] at (1.050,1.290) {};
\node[circle,inner sep=1.3pt,fill=black!40] at (1.050,1.780) {};
\node[circle,inner sep=1.3pt,fill=black!40] at (1.050,2.270) {};
\node[circle,inner sep=1.3pt,fill=black!40] at (1.050,2.760) {};
\node[circle,inner sep=1.3pt,fill=black!40] at (1.050,3.250) {};
\node[circle,inner sep=1.7pt,fill=MRed] at (1.750,0.800) {};
\node[circle,inner sep=1.3pt,fill=black!40] at (1.750,1.290) {};
\node[circle,inner sep=1.3pt,fill=black!40] at (1.750,1.780) {};
\node[circle,inner sep=1.3pt,fill=black!40] at (1.750,2.270) {};
\node[circle,inner sep=1.3pt,fill=black!40] at (1.750,2.760) {};
\node[circle,inner sep=1.3pt,fill=black!40] at (1.750,3.250) {};
\node[circle,inner sep=1.7pt,fill=MRed] at (2.450,0.800) {};
\node[circle,inner sep=1.3pt,fill=black!40] at (2.450,1.290) {};
\node[circle,inner sep=1.3pt,fill=black!40] at (2.450,1.780) {};
\node[circle,inner sep=1.3pt,fill=black!40] at (2.450,2.270) {};
\node[circle,inner sep=1.3pt,fill=black!40] at (2.450,2.760) {};
\node[circle,inner sep=1.3pt,fill=black!40] at (2.450,3.250) {};
\node[circle,inner sep=1.7pt,fill=MRed] at (3.150,0.800) {};
\node[circle,inner sep=1.3pt,fill=black!40] at (3.150,1.290) {};
\node[circle,inner sep=1.3pt,fill=black!40] at (3.150,1.780) {};
\node[circle,inner sep=1.3pt,fill=black!40] at (3.150,2.270) {};
\node[circle,inner sep=1.3pt,fill=black!40] at (3.150,2.760) {};
\node[circle,inner sep=1.3pt,fill=black!40] at (3.150,3.250) {};
\node[circle,inner sep=1.7pt,fill=MRed] at (3.850,0.800) {};
\node[circle,inner sep=1.3pt,fill=black!40] at (3.850,1.290) {};
\node[circle,inner sep=1.3pt,fill=black!40] at (3.850,1.780) {};
\node[circle,inner sep=1.3pt,fill=black!40] at (3.850,2.270) {};
\node[circle,inner sep=1.3pt,fill=black!40] at (3.850,2.760) {};
\node[circle,inner sep=1.3pt,fill=black!40] at (3.850,3.250) {};
\node[] at (2.100,0.050) {$m=6$};
\node[font=\normalsize] at (2.100,-0.590) {$\Pr(e\subseteq S)=\theta^{6}$};
\node[panel] at (7.100,4.080) {(b) $\pi=(2,2,2)$};
\draw[source] (4.950,.48) rectangle (5.750,3.50);
\node[text=MBlue,font=\scriptsize] at (5.350,3.730) {$P_{1}$};
\draw[source] (6.700,.48) rectangle (7.500,3.50);
\node[text=MBlue,font=\scriptsize] at (7.100,3.730) {$P_{2}$};
\draw[source] (8.450,.48) rectangle (9.250,3.50);
\node[text=MBlue,font=\scriptsize] at (8.850,3.730) {$P_{3}$};
\draw[alternative] (4.920,.57) rectangle (9.280,1.52);
\node[circle,inner sep=1.7pt,fill=MRed] at (5.350,0.800) {};
\node[circle,inner sep=1.7pt,fill=MRed] at (5.350,1.290) {};
\node[circle,inner sep=1.3pt,fill=black!40] at (5.350,1.780) {};
\node[circle,inner sep=1.3pt,fill=black!40] at (5.350,2.270) {};
\node[circle,inner sep=1.3pt,fill=black!40] at (5.350,2.760) {};
\node[circle,inner sep=1.3pt,fill=black!40] at (5.350,3.250) {};
\node[circle,inner sep=1.7pt,fill=MRed] at (7.100,0.800) {};
\node[circle,inner sep=1.7pt,fill=MRed] at (7.100,1.290) {};
\node[circle,inner sep=1.3pt,fill=black!40] at (7.100,1.780) {};
\node[circle,inner sep=1.3pt,fill=black!40] at (7.100,2.270) {};
\node[circle,inner sep=1.3pt,fill=black!40] at (7.100,2.760) {};
\node[circle,inner sep=1.3pt,fill=black!40] at (7.100,3.250) {};
\node[circle,inner sep=1.7pt,fill=MRed] at (8.850,0.800) {};
\node[circle,inner sep=1.7pt,fill=MRed] at (8.850,1.290) {};
\node[circle,inner sep=1.3pt,fill=black!40] at (8.850,1.780) {};
\node[circle,inner sep=1.3pt,fill=black!40] at (8.850,2.270) {};
\node[circle,inner sep=1.3pt,fill=black!40] at (8.850,2.760) {};
\node[circle,inner sep=1.3pt,fill=black!40] at (8.850,3.250) {};
\node[] at (7.100,0.050) {$m=3$};
\node[font=\normalsize] at (7.100,-0.590) {$\Pr(e\subseteq S)=\theta^{3}$};
\node[panel] at (12.100,4.080) {(c) $\pi=(4,1,1)$};
\draw[source] (9.950,.48) rectangle (10.750,3.50);
\node[text=MBlue,font=\scriptsize] at (10.350,3.730) {$P_{1}$};
\draw[source] (11.700,.48) rectangle (12.500,3.50);
\node[text=MBlue,font=\scriptsize] at (12.100,3.730) {$P_{2}$};
\draw[source] (13.450,.48) rectangle (14.250,3.50);
\node[text=MBlue,font=\scriptsize] at (13.850,3.730) {$P_{3}$};
\draw[alternative] (9.920,.57) -- (14.280,.57) -- (14.280,1.03) -- (10.790,1.03) -- (10.790,2.52) -- (9.920,2.52) -- cycle;
\node[circle,inner sep=1.7pt,fill=MRed] at (10.350,0.800) {};
\node[circle,inner sep=1.7pt,fill=MRed] at (10.350,1.290) {};
\node[circle,inner sep=1.7pt,fill=MRed] at (10.350,1.780) {};
\node[circle,inner sep=1.7pt,fill=MRed] at (10.350,2.270) {};
\node[circle,inner sep=1.3pt,fill=black!40] at (10.350,2.760) {};
\node[circle,inner sep=1.3pt,fill=black!40] at (10.350,3.250) {};
\node[circle,inner sep=1.7pt,fill=MRed] at (12.100,0.800) {};
\node[circle,inner sep=1.3pt,fill=black!40] at (12.100,1.290) {};
\node[circle,inner sep=1.3pt,fill=black!40] at (12.100,1.780) {};
\node[circle,inner sep=1.3pt,fill=black!40] at (12.100,2.270) {};
\node[circle,inner sep=1.3pt,fill=black!40] at (12.100,2.760) {};
\node[circle,inner sep=1.3pt,fill=black!40] at (12.100,3.250) {};
\node[circle,inner sep=1.7pt,fill=MRed] at (13.850,0.800) {};
\node[circle,inner sep=1.3pt,fill=black!40] at (13.850,1.290) {};
\node[circle,inner sep=1.3pt,fill=black!40] at (13.850,1.780) {};
\node[circle,inner sep=1.3pt,fill=black!40] at (13.850,2.270) {};
\node[circle,inner sep=1.3pt,fill=black!40] at (13.850,2.760) {};
\node[circle,inner sep=1.3pt,fill=black!40] at (13.850,3.250) {};
\node[] at (12.100,0.050) {$m=3$};
\node[font=\normalsize] at (12.100,-0.590) {$\Pr(e\subseteq S)=\theta^{3}$};
\node[note] at (7.100,-1.300) {Each source edge is deleted independently with probability $\theta$; $S$ is the deleted vertex set.};
\end{tikzpicture}%
}
    \caption{Three intersection patterns in a $6$-uniform hypergraph. Each blue block is an edge of the source matching $P$; the red vertices form an alternative edge $e$. If $S$ is the union of the independently deleted source edges, then $\Prob(e\subseteq S)=\theta^{m(\pi)}$. The last two patterns meet the same number of source edges but have different intersection sizes.}
    \label{fig:intersectionpatterns}
\end{figure}

Consider a $6$-uniform source matching $P$ whose edges are independently deleted with probability $\theta$.
The probability that an alternative edge $e$ lies in the deleted set depends on its pattern, as shown in Figure~\ref{fig:intersectionpatterns}.
The last two patterns have availability probability $\theta^3$, rather than $\theta^6$, but their coefficients are bounded using different codegrees.

\begin{definition}\label{def:exposure-polynomial}
    Under the hypotheses of Theorem~\ref{thm:weightedcount}, define
    \begin{equation}\label{eq:pattern-polynomial}
        L_H(y)=1+dy^{k-1}+\sum_{\pi\in\Pi_k}b_\pi(H) y^{m(\pi)-1}
        \quad(0\le y\le1),
    \end{equation}
    where $b_\pi(H)$ is defined in Definition~\ref{def:coefficients}.
\end{definition}

The terms in $L_H$ correspond to the source edge, alternatives meeting $k$ source blocks, and alternatives of patterns $\pi\in\Pi_k$, respectively.
The polynomial depends only on $H$, not on the matching law or the deletion probability.

After independent deletion with probability $\theta$, an alternative of pattern $\pi$ lies in the deleted set $S$ with probability $\theta^{m(\pi)}$.
Dividing the expected incidences by $\E[|S|]=\theta\E[|V(M)|]$, when this is positive, bounds the average contribution by $b_\pi(H)\theta^{m(\pi)-1}$.

In the random-order exposure below, vertices receive independent uniform priorities and each matching edge is revealed at its highest-priority vertex.
Conditional on the current vertex having priority $u$, every vertex of another block must have priority less than $u$ for that block to remain unexposed; this has probability $u^k$.
Writing $y=u^k$, the resulting averaged choice bound is $L_H(y)$, or $L_H(\theta y)$ when the source is first thinned.
In particular, $dy^{k-1}$ bounds the contribution from alternatives meeting $k$ distinct source blocks before thinning.
The normalized contribution of a pattern with $m$ parts is
\[
    \frac{b_\pi(H) y^{m-1}}{dy^{k-1}}=\frac{b_\pi(H)}{d}y^{-(k-m)}.
\]
This ratio can be large near $y=0$, even when $b_\pi(H)/d$ is small.
The next lemma estimates the difference between $\int_0^1\log L_H(y)\dd y$ and $\log d-k+1$.

\begin{lemma}\label{lem:integralbound}
    Under the hypotheses of Theorem~\ref{thm:weightedcount},
    \begin{equation}\label{eq:pattern-integral}
        \int_0^1\log\frac{L_H(y)}{dy^{k-1}}\dd y=o(d^{-1/k}).
    \end{equation}
\end{lemma}
\begin{proof}[Proof of Lemma~\ref{lem:integralbound}]
    For $0<x\le1$ and an integer $q\ge1$, split the integral of $\log(1+xy^{-q})$ at $x^{1/q}$.
    Below this point, substitute $y=x^{1/q}z$; the resulting integral is
    \[
        x^{1/q}\int_0^1\log(1+z^{-q})\dd z=O_q(x^{1/q}),
    \]
    since $\log(1+z^{-q})\le\log2-q\log z$.
    Above this point, use $\log(1+xy^{-q})\le xy^{-q}$.
    Thus
    \begin{equation}\label{eq:pattern-log-integral}
        \int_0^1\log(1+xy^{-q})\dd y
        =\begin{cases}
            O(x\log(e/x)),&q=1,\\
            O_q(x^{1/q}),&q\ge2.
        \end{cases}
    \end{equation}
    To estimate the integral in~\eqref{eq:pattern-integral}, use $\log(1+\sum_i x_i)\le\sum_i\log(1+x_i)$ and apply~\eqref{eq:pattern-log-integral} to $1/(dy^{k-1})$ and each $(b_\pi(H)/d)y^{-(k-m(\pi))}$.
    The first of these contributes $O_k(d^{-1/(k-1)})=o(d^{-1/k})$, since $d\to\infty$.
    For a pattern with $m$ parts and $q=k-m\ge2$, the contribution is at most
    \[
        O_k\bigl((b_\pi(H)/d)^{1/q}\bigr)
        =O_{k,C}\bigl(d^{-1/k}n^{-\eps/q}\bigr).
    \]
    If $q=1$, monotonicity of $x\log(e/x)$ on $(0,1]$ and $b_\pi(H)/d=O(d^{-1/k}n^{-\eps})$ give a contribution
    \[
        O\bigl(d^{-1/k}n^{-\eps}\log n\bigr),
    \]
    where simplicity gives $\log d=O_k(\log n)$.
    Each pattern contribution is $o(d^{-1/k})$, and there are only finitely many patterns for fixed $k$.
    This proves~\eqref{eq:pattern-integral} and completes the proof.
\end{proof}


\subsection{Exposure and reconstruction}\label{subsec:patternbounds}

The following bounds apply to any random matching of $H$, including Gibbs matchings at activities that depend on $n$.
In the random-order exposure, we must count the possible matching edges through each vertex according to their intersection patterns.
The next definition and the random-order exposure lemma apply to any simple $k$-uniform hypergraph; neither requires regularity or the hypotheses of Theorem~\ref{thm:weightedcount}.

\begin{definition}\label{def:vertex-patterns}
    Let $H$ be a simple $k$-uniform hypergraph, let $P$ be a matching of $H$, and let $v\in V(P)$.
    For $\pi\in\Pi_k$, define $b_{v,\pi}(P)$ to be the number of edges $e\in E(H)\setminus P$ such that $v\in e\subseteq V(P)$ and the multiset of nonzero intersection sizes $|e\cap B|$, over $B\in P$, is $\pi$.
    Define $g_v(P)$ to be the number of edges $e\in E(H)\setminus P$ such that $v\in e\subseteq V(P)$ and $e$ meets exactly $k$ edges of $P$.
\end{definition}

Each edge counted by $N_\pi(P)$ contributes to $b_{v,\pi}(P)$ at each of its $k$ vertices, so Definition~\ref{def:coefficients} gives
\begin{equation}\label{eq:pattern-incidences}
    \sum_{v\in V(P)}b_{v,\pi}(P)=kN_\pi(P)\le |V(P)|b_\pi(H).
\end{equation}
This identity allows us to use the intersection coefficients $b_\pi(H)$ after summing the conditional entropy bounds over vertices.

\begin{lemma}[Random-order exposure]\label{lem:random-order-entropy}
    Let $H$ be a simple $k$-uniform hypergraph. Fix $U\subseteq V(H)$ and let $P$ be a random perfect matching of $H[U]$.
    With $g_v(P)$ and $b_{v,\pi}(P)$ as in Definition~\ref{def:vertex-patterns}, for every law of $P$ we have
    \[
        \HH(P)\le\frac1k\int_0^1\E_P\left[\sum_{v\in U}
        \log\left(1+g_v(P)y^{k-1}+\sum_{\pi\in\Pi_k}b_{v,\pi}(P)y^{m(\pi)-1}\right)\right]\dd y.
    \]
\end{lemma}
\begin{proof}[Proof of Lemma~\ref{lem:random-order-entropy}]
    We use the entropy argument of Radhakrishnan~\cite{Radhakrishnan}, in the hypergraph form used by Luria~\cite[proof of Theorem~3.1]{Luria}.
    For each $v\in U$, let $P_v$ be the edge of $P$ containing $v$.
    The vector $(P_v)_{v\in U}$ determines $P$.
    Independently of $P$, give the vertices independent uniform priorities $\alpha_v\in[0,1]$ and expose the variables $P_v$ in decreasing priority order.
    For a fixed priority vector $\alpha$, write $\mathcal F_v^\alpha=(P_w:\alpha_w>\alpha_v)$ for the information exposed before $v$.
    The chain rule in Proposition~\ref{prop:entropybasics}, averaged over these orders, gives
    \begin{equation}\label{eq:exposure-chain-rule}
        \HH(P)=\E_\alpha\left[\sum_{v\in U}\HH(P_v\mid\mathcal F_v^\alpha)\right].
    \end{equation}

    Let $A_v$ be the event that $v$ has the largest priority in $P_v$.
    If $A_v$ fails, an earlier variable already reveals $P_v$, so the corresponding conditional entropy is zero.
    Whether this has happened is determined by $\mathcal F_v^\alpha$: it is enough to check whether one of its revealed edges contains $v$.
    On $A_v$, let $N_v(P,\alpha)$ count edges of $H[U]$ through $v$ disjoint from all previously revealed matching edges.
    These include every possible value of $P_v$ given the history, even though some may not extend to a perfect matching.
    Set $N_v(P,\alpha)=1$ off $A_v$.
    Applying the conditional support bound~\eqref{eq:conditionalsupport} to each summand in~\eqref{eq:exposure-chain-rule} yields
    \begin{equation}\label{eq:exposure-support-sum}
        \HH(P)\le\sum_{v\in U}\E_{P,\alpha}
        \bigl[\mathbf1_{A_v}\log N_v(P,\alpha)\bigr].
    \end{equation}

    We estimate each expectation on the right of~\eqref{eq:exposure-support-sum} by first fixing $P$ and $\alpha_v=u$.
    The event $A_v$ requires the other $k-1$ vertices of $P_v$ to have priority less than $u$, and therefore has probability $u^{k-1}$.
    Conditional on $A_v$, the edge $P_v$ is always available and contributes one choice.
    If an alternative edge through $v$ meets $m$ blocks of $P$, it is available exactly when every vertex of its other $m-1$ blocks has priority less than $u$.
    Those $k(m-1)$ vertices lie outside $P_v$, so their priorities remain independent after conditioning on $A_v$.
    Hence
    \[
        \E_\alpha[N_v(P,\alpha)\mid P,\alpha_v=u,A_v]
        =1+g_v(P)u^{k(k-1)}
        +\sum_{\pi\in\Pi_k}b_{v,\pi}(P)u^{k(m(\pi)-1)}.
    \]
    Jensen's inequality now gives
    \[
        \begin{aligned}
            &\E_\alpha[\mathbf1_{A_v}\log N_v(P,\alpha)\mid P,\alpha_v=u]\\
            &\quad\le u^{k-1}\log\left(1+g_v(P)u^{k(k-1)}
            +\sum_{\pi\in\Pi_k}b_{v,\pi}(P)u^{k(m(\pi)-1)}\right).
        \end{aligned}
    \]
    Substitute this bound into~\eqref{eq:exposure-support-sum}, integrate over $u\in[0,1]$, and average over $P$.
    The change of variable $y=u^k$ satisfies $u^{k-1}\dd u=\dd y/k$ and gives the asserted bound.
\end{proof}

For a $d$-regular host, $g_v(P)\le d$, so $dy^{k-1}$ is a reference upper estimate for the contribution of alternatives meeting $k$ blocks.
Its integral supplies the leading term in the entropy upper bound:
\[
    \frac{|U|}{k}\int_0^1\log(dy^{k-1})\dd y
    =\frac{|U|}{k}(\log d-k+1).
\]
The subtraction of $k-1$ records the depletion of choices during exposure.
Under the hypotheses of Theorem~\ref{thm:weightedcount}, Lemma~\ref{lem:integralbound} shows that replacing $dy^{k-1}$ by $L_H(y)$ changes the integral of its logarithm by $o(d^{-1/k})$.
The next lemma also retains the loss when the available degrees are smaller than $d$.

\begin{lemma}[Exposure given the uncovered set]\label{lem:exposure}
    Under the hypotheses of Theorem~\ref{thm:weightedcount}, let $M$ be any random matching of $H$.
    Write $W=V(H)\setminus V(M)$, suppose $s=\E_M[|V(M)|]>0$, and let
    \[
        S_W=\sum_{v\notin W}(d-d_{H-W}(v)).
    \]
    Then
    \[
        \HH(M\mid W)\le\frac sk(\log d-k+1)-\frac{\E_M[S_W]}{kd}+o(sd^{-1/k}).
    \]
\end{lemma}

\begin{proof}[Proof of Lemma~\ref{lem:exposure}]
    We apply Lemma~\ref{lem:random-order-entropy} conditionally on $W$ and retain the reduction in available degrees measured by $S_W$.
    All expectations denoted by $\E_M$ below are with respect to the unconditional law of $M$; both $W$ and $S_W$ are determined by $M$.

    For each realized $W=w$, apply Lemma~\ref{lem:random-order-entropy} to the conditional law of $M$ on $H-w$.
    The expectation in that application is $\E_{M\mid W=w}$.
    Averaging over $W$ changes it to $\E_W\E_{M\mid W}=\E_M$, giving
    \begin{equation}\label{eq:conditional-exposure}
        \HH(M\mid W)\le\frac1k\int_0^1\E_M\left[\sum_{v\notin W}
        \log\left(1+g_v(M)y^{k-1}+\sum_{\pi\in\Pi_k}b_{v,\pi}(M)y^{m(\pi)-1}\right)\right]\dd y.
    \end{equation}
    To estimate the right side of~\eqref{eq:conditional-exposure}, apply the tangent inequality
    \[
        \log t\le\log L_H(y)+\frac{t-L_H(y)}{L_H(y)}
    \]
    to each logarithm.
    For every $\pi\in\Pi_k$, the terms involving $b_{v,\pi}(M)-b_\pi(H)$ in this tangent bound have a nonpositive expected sum, because
    \[
        \E_M\left[\sum_{v\notin W}(b_{v,\pi}(M)-b_\pi(H))\right]
        =k\E_M[N_\pi(M)]-s b_\pi(H)\le0.
    \]
    Here we used~\eqref{eq:pattern-incidences} and $s=k\E_M[|M|]$.
    The degree contribution satisfies, for each matching,
    \[
        \sum_{v\notin W}(g_v(M)-d)
        \le\sum_{v\notin W}(d_{H-W}(v)-d)=-S_W.
    \]
    Thus the tangent bound in~\eqref{eq:conditional-exposure} yields
    \begin{equation}\label{eq:exposure-integral-loss}
        \HH(M\mid W)\le\frac sk\int_0^1\log L_H(y)\dd y
        -\frac{\E_M[S_W]}k\int_0^1\frac{y^{k-1}}{L_H(y)}\dd y.
    \end{equation}

    We now evaluate the two integrals in~\eqref{eq:exposure-integral-loss}.
    Lemma~\ref{lem:integralbound} gives
    \[
        \int_0^1\log L_H(y)\dd y=\log d-k+1+o(d^{-1/k}).
    \]
    Since $u/(1+u)\le\log(1+u)$ for $u\ge0$,
    \[
        0\le1-\frac{dy^{k-1}}{L_H(y)}
        \le\log\frac{L_H(y)}{dy^{k-1}}.
    \]
    Integrating this inequality and applying Lemma~\ref{lem:integralbound} again gives
    \[
        \int_0^1\frac{y^{k-1}}{L_H(y)}\dd y\ge\frac1d-o(d^{-1-1/k}).
    \]
    Substitute these two integral estimates into~\eqref{eq:exposure-integral-loss}.
    Since $0\le\E_M[S_W]\le sd$, both resulting error terms are $o(sd^{-1/k})$, yielding the asserted bound for $\HH(M\mid W)$.
    This completes the proof.
\end{proof}

The next lemma bounds the conditional entropy needed to recover a matching $M$ after independent edge deletion leaves a matching $J$.
The term $nf(p,\theta(1-p))$ accounts for identifying the original uncovered set, and the remaining terms bound the entropy of the deleted perfect matching once that set is known.
Its error is uniform in the deletion probability $\theta$, allowing $\theta$ to tend to zero in the lower-bound argument of Section~\ref{sec:counting}.

\begin{lemma}[Reconstruction after thinning]\label{lem:reconstruction}
    Under the hypotheses of Theorem~\ref{thm:weightedcount}, let $M$ be any random matching of $H$, and let $p=\E_M[|V(H)\setminus V(M)|]/n$.
    Delete each edge of $M$ independently with probability $\theta\in(0,1)$, leaving $J$.
    Then
    \[
        \HH(M\mid J)\le nf(p,\theta(1-p))
        +\frac{n\theta(1-p)}k\bigl(\log(d\theta^{k-1})-k+1\bigr)+o(nd^{-1/k}).
    \]
    Here $f$ is the function in Definition~\ref{def:entropy-f}.
    The error is uniform over $\theta\in(0,1)$ and the law of $M$.
\end{lemma}
\begin{proof}[Proof of Lemma~\ref{lem:reconstruction}]
    Write $W=V(H)\setminus V(M)$.
    Let $S=V(M\setminus J)$.
    We write $\E_{M,J}$ for expectation over the law of $M$ and the subsequent independent deletions, and $\E_{J\mid M}$ for expectation over deletions with $M$ fixed.
    Since $W$ is determined by $M$, the chain rule gives
    \begin{equation}\label{eq:reconstruction-split}
        \HH(M\mid J)=\HH(W\mid J)+\HH(M\mid J,W).
    \end{equation}
    Let $\alpha_v=\Prob(v\in W)$ and $\beta_v=\Prob(v\in S)$.
    These two disjoint cases exhaust the vertices not covered by $J$.
    Since $\{v\in W\cup S\}$ is determined by $J$, conditioning reduces entropy and gives
    \[
        \HH(\mathbf1_{\{v\in W\}}\mid J)
        \le\HH(\mathbf1_{\{v\in W\}}\mid\mathbf1_{\{v\in W\cup S\}})
        =f(\alpha_v,\beta_v).
    \]
    Conditional subadditivity and Lemma~\ref{lem:entropyconcavity}(i) now give
    \begin{equation}\label{eq:vacancy-classification}
        \HH(W\mid J) \le\sum_v f(\alpha_v,\beta_v) \le nf(p,\theta(1-p)),
    \end{equation}
    because $\sum_v\alpha_v=np$ and $\sum_v\beta_v=n\theta(1-p)$.

    We now bound the second term of~\eqref{eq:reconstruction-split}, the entropy of the deleted matching on $S$.
    Conditional on $J,W$, the set $S=V(H)\setminus(V(J)\cup W)$ is fixed, and $M\setminus J$ is a random perfect matching of $H[S]$.
    Apply Lemma~\ref{lem:random-order-entropy} to this conditional law, whose expectation is $\E_{M\mid J,W}$.
    Averaging over $(J,W)$ changes this to $\E_{J,W}\E_{M\mid J,W}=\E_{M,J}$, since $W$ is determined by $M$.
    Thus
    \begin{equation}\label{eq:deleted-exposure}
        \begin{aligned}
            \HH(M\mid J,W)\le\frac1k\int_0^1\E_{M,J}\Biggl[\sum_{v\in S}
            \log\Biggl(1+g_v(M\setminus J)y^{k-1}+\sum_{\pi\in\Pi_k}b_{v,\pi}(M\setminus J)y^{m(\pi)-1}\Biggr)\Biggr]\dd y.
        \end{aligned}
    \end{equation}
    Empty $S$ contributes zero to this bound.

    To estimate the expected sum in~\eqref{eq:deleted-exposure}, fix the source matching $M$ before averaging over the deletions.
    An alternative edge of pattern $\pi$ belongs to $H[S]$ exactly when all $m(\pi)$ source blocks it meets are deleted, an event of probability $\theta^{m(\pi)}$.
    Consequently, with $M$ fixed,
    \[
        \E_{J\mid M}\left[\sum_{v\in S}b_{v,\pi}(M\setminus J)\right]
        =k\theta^{m(\pi)}N_\pi(M).
    \]
    Averaging this identity over $M$ and using~\eqref{eq:pattern-incidences} gives
    \[
        \E_{M,J}\left[\sum_{v\in S}b_{v,\pi}(M\setminus J)\right]
        =k\theta^{m(\pi)}\E_M[N_\pi(M)]
        \le n(1-p)b_\pi(H)\theta^{m(\pi)}.
    \]
    The same argument for alternatives meeting $k$ blocks gives
    \[
        \E_{M,J}\left[\sum_{v\in S}g_v(M\setminus J)\right]
        \le n(1-p)d\theta^k.
    \]
    Also $\E_{M,J}[|S|]=n\theta(1-p)$.
    If $p=1$, the reconstruction bound is immediate, so suppose $p<1$.
    Dividing the last two incidence bounds by $\E_{M,J}[|S|]$ gives $b_\pi(H)\theta^{m(\pi)-1}$ and $d\theta^{k-1}$, respectively.
    For each $y$, Jensen's inequality over the vertices of $S$, followed by Lemma~\ref{lem:entropyconcavity}(ii) with weight $|S|$, therefore bounds the expected sum of logarithms in~\eqref{eq:deleted-exposure} by
    \[
        \E_{M,J}[|S|]\log\left(1+d(\theta y)^{k-1}
        +\sum_{\pi\in\Pi_k}b_\pi(H)(\theta y)^{m(\pi)-1}\right)
        =n\theta(1-p)\log L_H(\theta y).
    \]
    Substituting into~\eqref{eq:deleted-exposure} yields
    \begin{equation}\label{eq:deleted-reconstruction-cost}
        \HH(M\mid J,W)\le\frac{n\theta(1-p)}k\int_0^1\log L_H(\theta y)\dd y.
    \end{equation}
    Equations~\eqref{eq:vacancy-classification} and~\eqref{eq:deleted-reconstruction-cost} bound the two terms of~\eqref{eq:reconstruction-split}. It remains to estimate the integral in~\eqref{eq:deleted-reconstruction-cost}.
    To this end, factoring $d(\theta y)^{k-1}$ and substituting $x=\theta y$ gives
    \begin{equation}\label{eq:reconstruction-integral}
        \begin{aligned}
            0&\le\theta\int_0^1\log L_H(\theta y)\dd y
            -\theta\bigl(\log(d\theta^{k-1})-k+1\bigr)\\
            &=\int_0^\theta\log\frac{L_H(x)}{dx^{k-1}}\dd x
            \le\int_0^1\log\frac{L_H(x)}{dx^{k-1}}\dd x
            =o(d^{-1/k}).
        \end{aligned}
    \end{equation}
    The inequality uses the nonnegative integrand, and the final estimate is Lemma~\ref{lem:integralbound}.
    Multiplying~\eqref{eq:reconstruction-integral} by $n(1-p)/k$ evaluates the right side of~\eqref{eq:deleted-reconstruction-cost}, with error $o(nd^{-1/k})$ uniformly over $\theta$ and $p$.
    Combining this estimate with~\eqref{eq:vacancy-classification} in~\eqref{eq:reconstruction-split} proves the lemma.
    This completes the proof.
\end{proof}


\section{Refined partition function estimates}\label{sec:counting}

We prove Theorem~\ref{thm:weightedcount} by estimating the correction to the leading term of $\log Z_H(\lambda)$.
In the lower- and upper-bound arguments, $H$ satisfies the hypotheses of Theorem~\ref{thm:weightedcount}, and the target activity $\lambda$ belongs to a fixed compact interval $I\subset(0,\infty)$. All error terms are uniform over these $H$ and $\lambda$.
To shorten the repeated expressions in this section, write
\[
    a_\lambda=(\lambda d)^{-1/k}
\]
and define
\begin{equation}\label{eq:counting-remainder}
    R_\lambda=\log Z_H(\lambda)-\frac nk\bigl(\log(\lambda d)-k+1\bigr).
\end{equation}
Thus $a_\lambda$ is the proposed leading term for the average vacancy probability, and our aim is to prove $R_\lambda=(1+o(1))na_\lambda$.
We use these abbreviations only in this section.

Theorem~\ref{thm:greedycount} supplies the leading-term lower bound with error $O(n\Gamma(H)(\log n)^L)$, including at the source activity $\lambda d$.
We now delete edges from a Gibbs matching at this larger activity and apply the reconstruction estimates of Section~\ref{sec:patterns} to obtain the positive correction $na_\lambda$.
For the upper bound, we estimate the entropy of the uncovered set and of a perfect matching on the remaining vertices.


\subsection{Lower bound}\label{subsec:countlower}

\paragraph{Generating partial matchings from the source.}
Sample $M$ at activity $\lambda T$ and delete its edges independently with probability $\theta$, leaving $J$.
Lemma~\ref{lem:entropy-counting} gives the lower bound
\[
    \log Z_H(\lambda)\ge\HH(J)+\E[|J|]\log\lambda.
\]
We use Lemma~\ref{lem:reconstruction} to bound the conditional entropy appearing when $\HH(J)$ is expanded.

To see the expected gain, consider a perfect source and temporarily omit the $o(nd^{-1/k})$ term in Lemma~\ref{lem:reconstruction}.
Deletion choices contribute $nh(\theta)/k$ in entropy.
The change in log weight subtracts $n\theta\log\lambda/k$. The upper bound for $\HH(M\mid J)$ in Lemma~\ref{lem:reconstruction}, apart from the omitted term, contributes the subtraction $n\theta\{\log(d\theta^{k-1})-(k-1)\}/k$.
For small $\theta$, the resulting main contribution is
\[
    \frac nk\left[ h(\theta)-\theta\log\lambda -\theta\{\log(d\theta^{k-1})-(k-1)\}\right] =n\theta\left(1-\log\frac{\theta}{a_\lambda}\right)+O(n\theta^2).
\]
It is maximized to first order at $\theta=a_\lambda$, with gain $na_\lambda$, as illustrated in Figure~\ref{fig:thinninggain}.
This calculation explains the correction combinatorially.
Perfect sources are only an illustration; their existence is not assumed.

\begin{figure}
    \centering
    \singlespacing
    \matchingfigurestyles
    \resizebox{\linewidth}{!}{%
\begin{tikzpicture}[x=1cm,y=1cm]
\node[font=\normalsize] at (6.5,4.75) {$G(\theta)=n\theta\bigl(1-\log(\theta/a_\lambda)\bigr)$};
\begin{scope}[x=2.9cm,y=2.7cm]
  \path[fill=MBlue!7] (0,0) -- plot[domain=.003:2.7182818,samples=100] (\x,{\x*(1-ln(\x))}) -- cycle;
  \draw[arr] (0,-.35) -- (0,1.40);
  \draw[arr] (-.03,0) -- (3.20,0);
  \draw[MBlue,line width=1.15pt] (0,0) -- plot[domain=.003:3.02,samples=120] (\x,{\x*(1-ln(\x))});
  \draw[black!45,densely dashed] (1,0) -- (1,1) -- (0,1);
  \fill[MRed] (1,1) circle[radius=1.7pt];
  \draw (1,.025) -- (1,-.025);
  \draw (2.7182818,.025) -- (2.7182818,-.025);
  \node[below=3pt] at (1,0) {$1$};
  \node[below=3pt] at (2.7182818,0) {$e$};
  \node[left=4pt] at (0,1) {$1$};
  \node[below left=2pt] at (0,0) {$0$};
\end{scope}
\node[anchor=west] at (.12,3.90) {$G(\theta)/(na_\lambda)$};
\node at (5.00,-.85) {$\theta/a_\lambda$};
\draw[arr,MRed] (5.9,3.46) -- (3.05,2.75);
\node[anchor=west,text=MRed] at (5.9,3.46) {maximum};
\node at (11.45,2.97) {$a_\lambda=(\lambda d)^{-1/k}$};
\node at (11.45,1.91) {Optimal deletion probability};
\node[font=\normalsize] at (11.45,1.36) {$\theta=a_\lambda$};
\node at (11.45,.45) {Counting gain};
\node[font=\normalsize] at (11.45,-.10) {$na_\lambda=n(\lambda d)^{-1/k}$};
\node[note] at (6.5,-1.50) {Main contribution for a perfect source; the error terms are omitted.};
\end{tikzpicture}%
}
    \caption{The main contribution in the perfect-source calculation, normalized by $na_\lambda$, as a function of $\theta/a_\lambda$. Its maximum occurs at $\theta=a_\lambda$ and gives the correction $na_\lambda$. The error terms are omitted.}
    \label{fig:thinninggain}
\end{figure}

For the actual source, we must also distinguish original vacancies from vertices uncovered by deletion.
Increasing the source activity gives a compensating term proportional to the number of original vacancies.
The next lemma keeps this term and removes the unknown source vacancy proportion.
The $o(nd^{-1/k})$ term below is inherited from Lemma~\ref{lem:reconstruction} and is uniform over the source activity and the deletion probability.

\begin{lemma}[Thinning]\label{lem:transfer}
    Under the hypotheses of Theorem~\ref{thm:weightedcount}, let $\lambda\in I$, $T>1$ and $0<\theta<1$, and define
    \[
        g=\frac{h(\theta)-\theta\log(\lambda d\theta^{k-1})+(k-1)\theta}{k}.
    \]
    If $g<(\log T)/k$, then
    \[
        R_\lambda\ge R_{\lambda T}+n\left[g+\theta\log(1-T^{-1/k}e^g)\right]-o(na_\lambda).
    \]
    The error is uniform over $T$ and $\theta$ satisfying these conditions.
\end{lemma}
\begin{proof}[Proof of Lemma~\ref{lem:transfer}]
    Let $M$ have the Gibbs measure at $\lambda T$, and obtain $J$ by deleting each edge independently with probability $\theta$.
    With $p=\E_{\lambda T}[|\calU_{\lambda T}(H)|]/n$, we have $\E[|M|]=n(1-p)/k$ and $\E[|J|]=(1-\theta)n(1-p)/k$.
    The law of $J$ need not be Gibbs at $\lambda$, but Lemma~\ref{lem:entropy-counting} gives
    \begin{equation}\label{eq:thinning-entropy-bound}
        \log Z_H(\lambda)\ge\HH(J)+\E[|J|]\log\lambda.
    \end{equation}

    Given $M$, the deletion choices have entropy $|M|h(\theta)$, so the chain rule gives
    \[
        \HH(J)=\HH(M)+\frac{n(1-p)}k h(\theta)-\HH(M\mid J).
    \]
    Lemma~\ref{lem:reconstruction} bounds the last term by
    \[
        \HH(M\mid J)\le nf(p,\theta(1-p))
        +\frac{n\theta(1-p)}k\bigl(\log(d\theta^{k-1})-k+1\bigr)+o(na_\lambda).
    \]
    Since $M$ has the Gibbs measure at activity $\lambda T$, equality holds in Lemma~\ref{lem:entropy-counting} at that activity.
    More explicitly, the definition of this measure gives
    \[
        \begin{aligned}
            \HH(M)
            &=\log Z_H(\lambda T)-\E_{\lambda T}[|M|]\log(\lambda T)\\
            &=\log Z_H(\lambda T)-\frac{n(1-p)}k\log(\lambda T).
        \end{aligned}
    \]
    Substitute the reconstruction bound into the preceding chain-rule expression for $\HH(J)$, and then use this entropy identity and $\E[|J|]=(1-\theta)n(1-p)/k$ on the right of~\eqref{eq:thinning-entropy-bound}.
    This yields
    \begin{equation}\label{eq:thinning-expanded-count}
        \begin{aligned}
            \log Z_H(\lambda)
            &\geq \HH(M)+\frac{n(1-p)}k h(\theta)-\HH(M\mid J)+\E[|J|]\log\lambda\\
            &\geq \log Z_H(\lambda T)-nf(p,\theta(1-p))-o(na_\lambda)\\
            &\quad+\frac{n(1-p)}k\Bigl[h(\theta)-\log(\lambda T)+(1-\theta)\log\lambda-\theta\bigl(\log(d\theta^{k-1})-k+1\bigr)\Bigr].
        \end{aligned}
    \end{equation}
    Since $-\log(\lambda T)+(1-\theta)\log\lambda=-\log T-\theta\log\lambda$, the bracket in~\eqref{eq:thinning-expanded-count} equals $-\log T+kg$, by the definition of $g$.
    Substituting this expression into~\eqref{eq:thinning-expanded-count} gives
    \[
        \begin{aligned}
            \log Z_H(\lambda)\ge{}&\log Z_H(\lambda T)-\frac{n(1-p)}k\log T\\
            &+n(1-p)g-nf(p,\theta(1-p))-o(na_\lambda).
        \end{aligned}
    \]
    By~\eqref{eq:counting-remainder},
    \[
        \log Z_H(\lambda T)-\frac nk\bigl(\log(\lambda d)-k+1\bigr)
        =R_{\lambda T}+\frac nk\log T.
    \]
    Thus subtracting $\frac nk(\log(\lambda d)-k+1)$ from both sides of the preceding bound leaves the term $(np/k)\log T$, and gives
    \begin{equation}\label{eq:thinning-before-cancellation}
        R_\lambda\ge R_{\lambda T}+n\left[(1-p)g+\frac pk\log T-f(p,\theta(1-p))\right]-o(na_\lambda).
    \end{equation}
    The term $(np/k)\log T$ arises because $M$ has $n(1-p)/k$ edges on average, rather than $n/k$.
    The next inequality uses this term to bound the entropy of the original vacancies without estimating $p$.

    \begin{claim}\label{clm:binaryentropy}
        For the function $f$ in Definition~\ref{def:entropy-f}, if $x,y\ge0$ and $u>0$, then $f(x,y)\le ux-y\log(1-e^{-u})$.
    \end{claim}
    \begin{claimproof}[Proof of Claim~\ref{clm:binaryentropy}]
        If $x+y=0$, the assertion is immediate.
        By concavity of $h$, its tangent line at $e^{-u}$ gives
        \[
            h(t)\le ut-(1-t)\log(1-e^{-u})
        \]
        for $0\le t\le1$.
        Apply this bound with $t=x/(x+y)$ and multiply by $x+y$.
    \end{claimproof}

    To bound the bracket on the right of~\eqref{eq:thinning-before-cancellation}, apply Claim~\ref{clm:binaryentropy} with $x=p$, $y=\theta(1-p)$ and $u=(\log T)/k-g>0$.
    Then
    \begin{align*}
        (1-p)g+\frac pk\log T-f(p,\theta(1-p))
        & \ge g+\theta(1-p)\log(1-T^{-1/k}e^g)\\
        & \ge g+\theta\log(1-T^{-1/k}e^g)
    \end{align*}
    since the logarithm is negative and $1-p\le1$.
    Substituting this lower bound into~\eqref{eq:thinning-before-cancellation} proves the lemma.
    This completes the proof.
\end{proof}

With $p$ as in the proof, the expected number of uncovered vertices in $J$ is $n\{p+\theta(1-p)\}$.
Thus, choosing $\theta=a_\lambda$ need not give $J$ the target vacancy scale; Lemma~\ref{lem:transfer} instead gives a lower bound for the partition function uniformly over $p$.

\begin{proposition}\label{prop:countrefinement}
    Under the hypotheses of Theorem~\ref{thm:weightedcount}, uniformly for $\lambda\in I$,
    \[
        R_\lambda\ge na_\lambda-o(na_\lambda).
    \]
\end{proposition}
\begin{proof}[Proof of Proposition~\ref{prop:countrefinement}]
    The expansion $h(\theta)=\theta\log(1/\theta)+\theta+O(\theta^2)$ in the definition of $g$ in Lemma~\ref{lem:transfer} gives
    \[
        g=\theta\bigl(1-\log(\theta/a_\lambda)\bigr)+O(\theta^2).
    \]
    Choose $\theta=a_\lambda$, which lies in $(0,1)$ for sufficiently large $n$.
    Then $g=a_\lambda+O(a_\lambda^2)$.
    Take $T=d$ in Lemma~\ref{lem:transfer}.
    Since $g=o(1)$ and $\log d\to\infty$, its hypothesis $g<(\log d)/k$ holds for sufficiently large $n$.
    Also $d^{-1/k}e^g=O(a_\lambda)$, so
    \[
        a_\lambda\log(1-d^{-1/k}e^g)=O(a_\lambda^2)=o(a_\lambda).
    \]
    These estimates for the two terms in the bracket in Lemma~\ref{lem:transfer} yield
    \begin{equation}\label{eq:thinning-lower-correction}
        R_\lambda\ge R_{\lambda d}+na_\lambda-o(na_\lambda).
    \end{equation}

    Condition~\eqref{eq:admissibility} gives $\Gamma(H)(\log n)^L=o(a_\lambda)$ for every fixed $L$.
    By simplicity, $\lambda d$ lies in a fixed polynomial activity range, so Theorem~\ref{thm:greedycount} gives
    \[
        R_{\lambda d}\ge-O\bigl(n\Gamma(H)(\log n)^L\bigr)=-o(na_\lambda).
    \]
    Substituting into~\eqref{eq:thinning-lower-correction} proves the asserted lower bound.
    This completes the proof.
\end{proof}


\subsection{Upper bound}\label{subsec:countupper}

We describe a matching by its uncovered set and a perfect matching on the remaining vertices.
Choosing the uncovered set contributes entropy, but each uncovered vertex also removes edge choices from the second part of this description.
Lemma~\ref{lem:exposure} bounds the conditional entropy of the matching given its uncovered set, retaining this degree loss.
Its $o(sd^{-1/k})$ term is $o(na_\lambda)$ because $s\le n$ and $\lambda\in I$.
We express the degree loss in terms of the number of edges meeting the uncovered set and apply Shearer's inequality.
The Gibbs variational inequality, with auxiliary weights chosen in terms of $a_\lambda$, then cancels the terms involving the unknown vacancy density.

\begin{proposition}\label{prop:deficit}
    Under the hypotheses of Theorem~\ref{thm:weightedcount}, uniformly for $\lambda\in I$,
    \[
        R_\lambda\le na_\lambda+o(na_\lambda).
    \]
\end{proposition}
\begin{proof}[Proof of Proposition~\ref{prop:deficit}]
    Let $M$ have the Gibbs measure at activity $\lambda$, write $W=\calU_\lambda(H)$, and let $p=\E_\lambda[|W|]/n$.
    We first bound the conditional entropy of $M$ given $W$.
    Since $W$ is determined by $M$, the chain rule and equality in Lemma~\ref{lem:entropy-counting} give
    \begin{equation}\label{eq:gibbs-uncovered-decomposition}
        \log Z_H(\lambda)=\HH(W)+\HH(M\mid W)+\frac{n(1-p)}k\log\lambda.
    \end{equation}
    Apply Lemma~\ref{lem:exposure} with $s=n(1-p)$ to the second term of~\eqref{eq:gibbs-uncovered-decomposition} and subtract $\frac nk(\log(\lambda d)-k+1)$, as in~\eqref{eq:counting-remainder}.
    We obtain
    \begin{equation}\label{eq:entropy-before-degree-loss}
        R_\lambda\le\HH(W)-\frac{\E_\lambda[S_W]}{kd}
        -\frac{np}{k}\bigl(\log(\lambda d)-k+1\bigr)+o(na_\lambda),
    \end{equation}
    where $S_W=\sum_{v\notin W}(d-d_{H-W}(v))$.

    The next claim bounds $\HH(W)-\E_\lambda[S_W]/(kd)$ in~\eqref{eq:entropy-before-degree-loss}.
    Its term proportional to $np$ cancels the corresponding term in that inequality.

    \begin{claim}\label{clm:vacancy-degree-loss}
        We have
        \[
            \HH(W)-\frac{\E_\lambda[S_W]}{kd}
            \le\frac{np}{k}\bigl(\log(\lambda d)-k+1\bigr)+na_\lambda+o(na_\lambda).
        \]
    \end{claim}
    \begin{claimproof}[Proof of Claim~\ref{clm:vacancy-degree-loss}]
        An edge $e$ meeting $W$ contributes $k-|e\cap W|$ to $S_W$, and an edge disjoint from $W$ contributes zero.
        Since $\sum_e|e\cap W|=d|W|$, we have
        \begin{equation}\label{eq:degree-loss-edges}
            S_W=k\sum_{e\in E(H)}\mathbf1_{\{e\cap W\ne\varnothing\}}-d|W|.
        \end{equation}
        Taking expectations and applying Lemma~\ref{lem:shearer} gives
        \begin{equation}\label{eq:vacancy-shearer-bound}
            \HH(W)-\frac{\E_\lambda[S_W]}{kd}
            \le\frac1d\sum_{e\in E(H)}
            \bigl(\HH(W\cap e)-\Prob_\lambda(W\cap e\ne\varnothing)\bigr)
            +\frac{np}{k}.
        \end{equation}

        Fix an edge $e$ and assign the following positive weights to its subsets, where $\mathrm e=\exp(1)$:
        \[
            w(A)=
            \begin{cases}
                1,&A=\varnothing,\\
                (\mathrm e a_\lambda)^{|A|}/\mathrm e,&A\ne\varnothing.
            \end{cases}
        \]
        Then
        \[
            \log w(A)=(1+\log a_\lambda)|A|-\mathbf1_{\{A\ne\varnothing\}},
        \]
        and the binomial theorem gives
        \[
            \sum_{A\subseteq e}w(A)
            =1+\frac{(1+\mathrm e a_\lambda)^k-1}{\mathrm e}
            =1+ka_\lambda+O_k(a_\lambda^2).
        \]
        Lemma~\ref{lem:entropy-counting}, applied to $W\cap e$ with these weights, therefore gives
        \begin{equation}\label{eq:edge-entropy-hit}
            \HH(W\cap e)-\Prob_\lambda(W\cap e\ne\varnothing)
            \le-(1+\log a_\lambda)\E_\lambda[|W\cap e|]+ka_\lambda+O_k(a_\lambda^2).
        \end{equation}
        This estimate holds for any law of $W\cap e$ as $a_\lambda\to0$.

        To sum~\eqref{eq:edge-entropy-hit}, use $|E(H)|=nd/k$ and
        \[
            \sum_{e\in E(H)}\E_\lambda[|W\cap e|]=d\E_\lambda[|W|]=dnp.
        \]
        Substituting into~\eqref{eq:vacancy-shearer-bound} gives
        \[
            \begin{aligned}
                \HH(W)-\frac{\E_\lambda[S_W]}{kd}
                &\le-np(1+\log a_\lambda)+na_\lambda+\frac{np}{k}+O_k(na_\lambda^2)\\
                &=\frac{np}{k}\bigl(\log(\lambda d)-k+1\bigr)+na_\lambda+O_k(na_\lambda^2).
            \end{aligned}
        \]
        The equality uses $\log(\lambda d)=-k\log a_\lambda$.
        Since $a_\lambda\to0$ uniformly for $\lambda\in I$, we have $na_\lambda^2=o(na_\lambda)$, proving the claim.
    \end{claimproof}

    Applying Claim~\ref{clm:vacancy-degree-loss} to~\eqref{eq:entropy-before-degree-loss} cancels the terms proportional to $np$ and gives
    \[
        R_\lambda\le na_\lambda+o(na_\lambda).
    \]
    This completes the proof.
\end{proof}

\begin{proof}[Proof of Theorem~\ref{thm:weightedcount}]
    Proposition~\ref{prop:deficit} and Proposition~\ref{prop:countrefinement} give
    \[
        R_\lambda=na_\lambda+o(na_\lambda)
    \]
    uniformly for $\lambda\in I$.
    By~\eqref{eq:counting-remainder}, this is~\eqref{eq:weightedcount}; taking $\lambda=1$ gives~\eqref{eq:totalcount}.
    This completes the proof.
\end{proof}


\section{Matching statistics}\label{sec:proofs}

We deduce Theorem~\ref{thm:statistics} from the counting estimate in Theorem~\ref{thm:weightedcount} using the exponential-moment identity~\eqref{eq:vacancy-mgf}.
Note that we have $d\le n^{k-1}$, so $n(\lambda d)^{-1/k}\to\infty$ uniformly on each fixed compact interval of positive activities.
All activity intervals below may be enlarged by a fixed amount when applying the counting estimate.


\subsection{Expectation and variance}\label{subsec:expectation}\label{sec:mean}\label{subsec:variance}\label{sec:variance}

We prove $L^2$ convergence of the normalized vacancy count, which gives the mean asymptotic and the variance bound together.
For a small tilt $s$, the leading term of the centered log moment-generating function is $n(\lambda d)^{-1/k}(e^s-1-s)$.
Its quadratic behavior near zero yields tail bounds that can be integrated to control the second moment.

\begin{proposition}\label{prop:variance}
    Suppose that $n(\lambda d)^{-1/k}\to\infty$ at a fixed activity $\lambda>0$ and that the counting estimate~\eqref{eq:weightedcount} holds uniformly on a fixed neighborhood of $\lambda$.
    Then
    \[
        \E_\lambda\left[\bigl(|\calU_\lambda(H)|-n(\lambda d)^{-1/k}\bigr)^2\right]=o(n^2(\lambda d)^{-2/k}).
    \]
    In particular, $|\calU_\lambda(H)|/(n(\lambda d)^{-1/k})\to1$ in $L^2$, and the same second-moment bound holds for $\Var_\lambda\bigl[|\calU_\lambda(H)|\bigr]$.
    For every fixed $\eta>0$,
    \[
        \Prob_\lambda(\bigl\lvert\,|\calU_\lambda(H)|-n(\lambda d)^{-1/k}\,\bigr\rvert>\eta n(\lambda d)^{-1/k})\le2e^{-c_\eta n(\lambda d)^{-1/k}}
    \]
    for all sufficiently large $n$.
\end{proposition}
\begin{proof}[Proof of Proposition~\ref{prop:variance}]
    Choose a fixed $s_0>0$ such that $\lambda e^{-ks}$ stays in the given neighborhood for $|s|\le s_0$.
    By~\eqref{eq:vacancy-mgf} and~\eqref{eq:weightedcount}, uniformly over this interval,
    \begin{align}
        \log\E_\lambda\bigl[e^{s(|\calU_\lambda(H)|-n(\lambda d)^{-1/k})}\bigr]
        &=sn-sn(\lambda d)^{-1/k}+\log Z_H(\lambda e^{-ks})-\log Z_H(\lambda)\notag\\
        &=n(\lambda d)^{-1/k}(e^s-1-s)+o(n(\lambda d)^{-1/k})\notag\\
        &\le Cn(\lambda d)^{-1/k}s^2+o(n(\lambda d)^{-1/k}),
    \label{eq:centered-vacancy-mgf}
    \end{align}
    where the last inequality follows from $e^s-1-s=O(s^2)$ for $|s|\le s_0$, by Taylor expansion at zero.
    \begin{claim}\label{clm:vacancy-tails}
        Uniformly for $x\ge0$,
        \[
            \begin{aligned}
                &\Prob_\lambda\bigl(\bigl\lvert\,|\calU_\lambda(H)|-n(\lambda d)^{-1/k}\,\bigr\rvert\ge x\bigr)\\
                &\quad\le\min\left\{1,\,2\exp\left(o(n(\lambda d)^{-1/k})
                -c\min\left\{\frac{x^2}{n(\lambda d)^{-1/k}},x\right\}\right)\right\},
            \end{aligned}
        \]
        where the $o(n(\lambda d)^{-1/k})$ error does not depend on $x$.
    \end{claim}
    \begin{claimproof}[Proof of Claim~\ref{clm:vacancy-tails}]
        For $x>0$, choose
        \[
            s=\min\left\{s_0,\frac{x}{2Cn(\lambda d)^{-1/k}}\right\}.
        \]
        On the upper-tail event, $e^{s(|\calU_\lambda(H)|-n(\lambda d)^{-1/k})}\ge e^{sx}$.
        Markov's inequality and~\eqref{eq:centered-vacancy-mgf} therefore give
        \[
            \begin{aligned}
                &\Prob_\lambda\bigl(|\calU_\lambda(H)|-n(\lambda d)^{-1/k}\ge x\bigr)\\
                &\quad\le e^{-sx}\E_\lambda\bigl[e^{s(|\calU_\lambda(H)|-n(\lambda d)^{-1/k})}\bigr]\\
                &\quad\le\exp\bigl(-sx+Cn(\lambda d)^{-1/k}s^2+o(n(\lambda d)^{-1/k})\bigr).
            \end{aligned}
        \]
        The choice $s\le x/(2Cn(\lambda d)^{-1/k})$ ensures that
        \[
            Cn(\lambda d)^{-1/k}s^2\le\frac{sx}{2}.
        \]
        Moreover,
        \[
            \begin{aligned}
                \frac{sx}{2}
                &=\min\left\{\frac{s_0x}{2},\frac{x^2}{4Cn(\lambda d)^{-1/k}}\right\}\\
                &\ge\min\left\{\frac{s_0}{2},\frac1{4C}\right\}
                \min\left\{x,\frac{x^2}{n(\lambda d)^{-1/k}}\right\}.
            \end{aligned}
        \]
        Thus the upper tail is at most the exponential in the claim, without the factor two, with $c=\min\{s_0/2,1/(4C)\}$.

        For the lower tail, $|\calU_\lambda(H)|-n(\lambda d)^{-1/k}\le-x$ implies
        $e^{-s(|\calU_\lambda(H)|-n(\lambda d)^{-1/k})}\ge e^{sx}$.
        Applying Markov's inequality with $-s$ gives
        \[
            \begin{aligned}
                &\Prob_\lambda\bigl(|\calU_\lambda(H)|-n(\lambda d)^{-1/k}\le-x\bigr)\\
                &\quad\le e^{-sx}\E_\lambda\bigl[e^{-s(|\calU_\lambda(H)|-n(\lambda d)^{-1/k})}\bigr]\\
                &\quad\le\exp\bigl(-sx+Cn(\lambda d)^{-1/k}s^2+o(n(\lambda d)^{-1/k})\bigr),
            \end{aligned}
        \]
        so the same calculation bounds this tail.
        The error in~\eqref{eq:centered-vacancy-mgf} has a common nonnegative $o(n(\lambda d)^{-1/k})$ upper bound for all $|s|\le s_0$.
        We use this bound for both signs of $s$, so it remains independent of $x$, even though $s$ depends on $x$.
        Adding the two tail bounds and using that a probability is at most one proves the claim for $x>0$.
        The case $x=0$ follows from the same trivial bound.
    \end{claimproof}

    For each fixed $\eta>0$, apply Claim~\ref{clm:vacancy-tails} with $x=\eta n(\lambda d)^{-1/k}$.
    The asserted tail probability is at most
    \[
        \begin{aligned}
            &2\exp\bigl(-c\min\{\eta^2,\eta\}n(\lambda d)^{-1/k}+o(n(\lambda d)^{-1/k})\bigr)\\
            &\quad\le2\exp\left(-\frac c2\min\{\eta^2,\eta\}n(\lambda d)^{-1/k}\right)
        \end{aligned}
    \]
    for all sufficiently large $n$.
    This proves the exponential tail bound with $c_\eta=(c/2)\min\{\eta^2,\eta\}$.

    To estimate the second moment, fix $0<\eta<1$.
    For $x\ge\eta n(\lambda d)^{-1/k}$, we have
    \[
        \min\left\{\frac{x^2}{n(\lambda d)^{-1/k}},x\right\}\ge\eta x.
    \]
    Since the error in Claim~\ref{clm:vacancy-tails} is independent of $x$, for all sufficiently large $n$ it satisfies
    \[
        o(n(\lambda d)^{-1/k})
        \le\frac{c\eta^2}{2}n(\lambda d)^{-1/k}
        \le\frac{c\eta x}{2}
    \]
    throughout this range.
    Substituting these inequalities into Claim~\ref{clm:vacancy-tails} gives
    \begin{equation}\label{eq:vacancy-integrable-tail}
        \Prob_\lambda\bigl(\bigl\lvert\,|\calU_\lambda(H)|-n(\lambda d)^{-1/k}\,\bigr\rvert>x\bigr)
        \le2e^{-c\eta x/2}
    \end{equation}
    whenever $x\ge\eta n(\lambda d)^{-1/k}$.

    For each outcome, the squared deviation is the integral of $2x$ from zero to the absolute deviation.
    Taking expectations and interchanging the nonnegative integral with expectation gives the first equality below.
    We then split the integral at $\eta n(\lambda d)^{-1/k}$, use the bound one for the probability below this point, and apply~\eqref{eq:vacancy-integrable-tail} above it:
    \[
        \begin{aligned}
            &\E_\lambda\left[\bigl(|\calU_\lambda(H)|-n(\lambda d)^{-1/k}\bigr)^2\right]\\
            &\quad=2\int_0^\infty x\Prob_\lambda\bigl(\bigl\lvert\,|\calU_\lambda(H)|-n(\lambda d)^{-1/k}\,\bigr\rvert>x\bigr)\dd x\\
            &\quad\le2\int_0^{\eta n(\lambda d)^{-1/k}}x\dd x
            +4\int_{\eta n(\lambda d)^{-1/k}}^\infty xe^{-c\eta x/2}\dd x\\
            &\quad=\eta^2n^2(\lambda d)^{-2/k}
            +4\int_{\eta n(\lambda d)^{-1/k}}^\infty xe^{-c\eta x/2}\dd x.
        \end{aligned}
    \]
    Integration by parts evaluates the last integral:
    \[
        \begin{aligned}
            4\int_{\eta n(\lambda d)^{-1/k}}^\infty xe^{-c\eta x/2}\dd x
            &=\left[-\left(\frac{8x}{c\eta}+\frac{16}{c^2\eta^2}\right)e^{-c\eta x/2}\right]_{\eta n(\lambda d)^{-1/k}}^\infty\\
            &=\left(\frac{8n(\lambda d)^{-1/k}}c+\frac{16}{c^2\eta^2}\right)
            e^{-c\eta^2n(\lambda d)^{-1/k}/2}.
        \end{aligned}
    \]
    For fixed $\eta$, this is $o(n^2(\lambda d)^{-2/k})$ because $n(\lambda d)^{-1/k}\to\infty$.
    Dividing the second-moment bound by $n^2(\lambda d)^{-2/k}$ and taking the upper limit as $n\to\infty$, we obtain
    \[
        \limsup_{n\to\infty}
        \frac{\E_\lambda\left[\bigl(|\calU_\lambda(H)|-n(\lambda d)^{-1/k}\bigr)^2\right]}
        {n^2(\lambda d)^{-2/k}}\le\eta^2.
    \]
    Letting $\eta\to0$ proves the stated $L^2$ convergence.
    The variance bound follows from
    \[
        \begin{aligned}
            \Var_\lambda\bigl[|\calU_\lambda(H)|\bigr]
            &=\E_\lambda\left[\bigl(|\calU_\lambda(H)|-n(\lambda d)^{-1/k}\bigr)^2\right]\\
            &\quad-\bigl(\E_\lambda[|\calU_\lambda(H)|]-n(\lambda d)^{-1/k}\bigr)^2\\
            &\le\E_\lambda\left[\bigl(|\calU_\lambda(H)|-n(\lambda d)^{-1/k}\bigr)^2\right].
        \end{aligned}
    \]
    When the counting estimate is uniform on an enlarged compact activity interval, the choices and bounds above are uniform on the original interval as well.
    This completes the proof.
\end{proof}

Proposition~\ref{prop:variance} proves Theorem~\ref{thm:statistics}(a).
To obtain the expected vacancy, apply Cauchy--Schwarz to the deviation and the constant random variable one:
\[
    \begin{aligned}
        \bigl|\E_\lambda[|\calU_\lambda(H)|]-n(\lambda d)^{-1/k}\bigr|
        &=\bigl|\E_\lambda[|\calU_\lambda(H)|-n(\lambda d)^{-1/k}]\bigr|\\
        &\le\left(\E_\lambda\left[\bigl(|\calU_\lambda(H)|-n(\lambda d)^{-1/k}\bigr)^2\right]\right)^{1/2}\\
        &=o(n(\lambda d)^{-1/k}).
    \end{aligned}
\]
Dividing by $n(\lambda d)^{-1/k}$ gives
\begin{equation}\label{eq:mean-vacancy}
    \E_\lambda[|\calU_\lambda(H)|]=(1+o(1))n(\lambda d)^{-1/k}.
\end{equation}


\subsection{Large deviations}\label{subsec:largedeviations}

For deviations by a fixed factor, it suffices to change the activity by a fixed factor.
The counting estimate identifies the exponential-moment limit, and the variance bound shows that the tilted law is concentrated near the desired vacancy.
This gives both sides of the tail asymptotic.

\begin{proof}[Proof of Theorem~\ref{thm:statistics}(b)]
    By~\eqref{eq:vacancy-mgf} and~\eqref{eq:weightedcount}, for every fixed real $z$,
    \begin{equation}\label{eq:vacancy-ld-mgf}
        \frac{1}{nd^{-1/k}}\log\E\bigl[e^{z|\calU(H)|}\bigr]=e^z-1+o(1).
    \end{equation}
    For $x>1$, Chernoff's inequality and~\eqref{eq:vacancy-ld-mgf}, with $z=\log x>0$, give
    \[
        \log\Prob(|\calU(H)|\ge xnd^{-1/k})
        \le nd^{-1/k}(x-1-x\log x+o(1)).
    \]
    For $0<x<1$, use $z=\log x<0$ and the event $|\calU(H)|\le xnd^{-1/k}$ to obtain the same bound for the lower tail.

    For the lower bound, we change the activity so that the desired vacancy is typical, and then compare the probabilities of the same matchings under the two measures.
    Fix $y>0$ and choose activity $y^{-k}$, since
    \[
        n(y^{-k}d)^{-1/k}=ynd^{-1/k}.
    \]
    For fixed $\eta>0$, let
    \[
        E=\left\{M\in\calM(H):
        \bigl|\,n-k|M|-ynd^{-1/k}\,\bigr|\le\eta nd^{-1/k}\right\}.
    \]
    Thus $E$ consists of matchings whose number of uncovered vertices is within $\eta nd^{-1/k}$ of $ynd^{-1/k}$.
    Proposition~\ref{prop:variance} and Markov's inequality applied to the squared deviation give
    \[
        \begin{aligned}
            \Prob_{y^{-k}}(E^c)
            &\le\frac{\E_{y^{-k}}\left[\bigl(|\calU_{y^{-k}}(H)|-ynd^{-1/k}\bigr)^2\right]}
            {\eta^2n^2d^{-2/k}}\\
            &=o(1).
        \end{aligned}
    \]
    Here $y$ and $\eta$ are fixed as $n\to\infty$.
    In particular, $\Prob_{y^{-k}}(E)\to1$.

    For each fixed matching $M$, the definition of the Gibbs measure gives
    \[
        \frac{\Prob_1(M)}{\Prob_{y^{-k}}(M)}
        =\frac{1/Z_H(1)}{y^{-k|M|}/Z_H(y^{-k})}
        =y^{k|M|}\frac{Z_H(y^{-k})}{Z_H(1)}.
    \]
    By the counting estimate~\eqref{eq:weightedcount},
    \[
        \begin{aligned}
            \log Z_H(y^{-k})-\log Z_H(1)
            &=\frac nk\log(y^{-k})+n(y^{-k}d)^{-1/k}-nd^{-1/k}+o(nd^{-1/k})\\
            &=-n\log y+(y-1)nd^{-1/k}+o(nd^{-1/k}).
        \end{aligned}
    \]
    Taking the logarithm of the probability ratio and substituting this difference gives
    \begin{equation}\label{eq:vacancy-probability-comparison}
        \begin{aligned}
            \log\frac{\Prob_1(M)}{\Prob_{y^{-k}}(M)}
            &=k|M|\log y+\log Z_H(y^{-k})-\log Z_H(1)\\
            &=-(n-k|M|)\log y+(y-1)nd^{-1/k}+o(nd^{-1/k})\\
            &=-(y\log y-y+1)nd^{-1/k}\\
            &\quad-\bigl(n-k|M|-ynd^{-1/k}\bigr)\log y+o(nd^{-1/k}).
        \end{aligned}
    \end{equation}
    For $M\in E$, the definition of $E$ gives
    \[
        -\bigl(n-k|M|-ynd^{-1/k}\bigr)\log y
        \ge-\eta nd^{-1/k}|\log y|.
    \]
    Substituting this bound into~\eqref{eq:vacancy-probability-comparison}, we obtain
    \[
        \log\frac{\Prob_1(M)}{\Prob_{y^{-k}}(M)}
        \ge-nd^{-1/k}\bigl(y\log y-y+1+\eta|\log y|+o(1)\bigr).
    \]
    The error comes from the two partition functions and is independent of $M$.
    Exponentiating this bound and summing over $M\in E$, we obtain
    \[
        \begin{aligned}
            \Prob_1(E)
            &=\sum_{M\in E}\frac{\Prob_1(M)}{\Prob_{y^{-k}}(M)}\Prob_{y^{-k}}(M)\\
            &\ge e^{-nd^{-1/k}(y\log y-y+1+\eta|\log y|+o(1))}
            \Prob_{y^{-k}}(E).
        \end{aligned}
    \]
    Since $\Prob_{y^{-k}}(E)\to1$, its logarithm is $o(1)$.
    Taking logarithms in the preceding bound, dividing by $nd^{-1/k}\to\infty$, and then taking the lower limit gives
    \begin{equation}\label{eq:vacancy-event-lower}
        \liminf_{n\to\infty}\frac{\log\Prob_1(E)}{nd^{-1/k}}
        \ge-\bigl(y\log y-y+1+\eta|\log y|\bigr).
    \end{equation}

    For the upper tail, fix $x>1$, choose $y>x$, and then choose $0<\eta<y-x$.
    Every $M\in E$ satisfies
    \[
        n-k|M|\ge(y-\eta)nd^{-1/k}>xnd^{-1/k}.
    \]
    Thus $E$ is contained in the upper-tail event, and~\eqref{eq:vacancy-event-lower} yields
    \[
        \liminf_{n\to\infty}
        \frac{\log\Prob(|\calU(H)|\ge xnd^{-1/k})}{nd^{-1/k}}
        \ge-\bigl(y\log y-y+1+\eta|\log y|\bigr).
    \]
    First let $\eta\downarrow0$ with $y$ fixed, and then let $y\downarrow x$.
    Continuity of $y\log y-y+1$ gives the lower bound $-(x\log x-x+1)$, which matches the upper bound already proved.

    For the lower tail, fix $0<x<1$, choose $0<y<x$, and then choose $0<\eta<x-y$.
    Every $M\in E$ now satisfies
    \[
        n-k|M|\le(y+\eta)nd^{-1/k}<xnd^{-1/k}.
    \]
    Consequently,~\eqref{eq:vacancy-event-lower} gives
    \[
        \liminf_{n\to\infty}
        \frac{\log\Prob(|\calU(H)|\le xnd^{-1/k})}{nd^{-1/k}}
        \ge-\bigl(y\log y-y+1+\eta|\log y|\bigr).
    \]
    Letting $\eta\downarrow0$ and then $y\uparrow x$ gives the matching lower bound for this tail.
    Choosing $y$ strictly beyond $x$ in the relevant direction ensures that the whole event $E$ lies in the desired tail; concentration at $y=x$ alone would not give this inclusion.
    This completes the proof.
\end{proof}


\section{Proofs of the applications}\label{sec:applicationproofs}

We now prove the statements of Section~\ref{sec:applications}, using the counting and statistical results established above.
We retain $N,K,D$ for the order, uniformity and degree of each auxiliary hypergraph. When verifying~\eqref{eq:admissibility}, its parameters $n,k,d$ are replaced by $N,K,D$.

\subsection{Latin squares and toroidal queens}\label{subsec:latinqueensproofs}

\begin{proof}[Proof of Corollary~\ref{cor:partiallatin}]
    Let $R,C,S$ be disjoint sets of $n$ row, column and symbol labels.
    The auxiliary hypergraph $\calH$ has vertex set $(R\times C)\cup(R\times S)\cup(C\times S)$.
    Each triple $(i,j,a)\in R\times C\times S$ gives the edge $\{(i,j),(i,a),(j,a)\}$, representing symbol $a$ in cell $(i,j)$.
    Two such edges meet exactly when their entries repeat a cell, a symbol in a row, or a symbol in a column.
    Thus matchings in $\calH$ are in bijection with $\mathcal{P}_n$.
    The hypergraph is linear, $3$-uniform and $n$-regular on $3n^2$ vertices.
    Its degree is a fixed positive power of its order, and
    \[
        |\calU(\calH)|=3(n^2-|P|).
    \]
    Theorem~\ref{thm:linearcount} at activity one gives the count, and Theorem~\ref{thm:statistics} gives the $L^2$ conclusion.
    The expectation and variance estimates follow from this convergence.
    Row and column permutations act transitively on the cells and preserve the uniform measure, so each cell is empty with probability $\E[n^2-|P|]/n^2\sim n^{-1/3}$.
    This completes the proof.
\end{proof}

\begin{proof}[Proof of Corollary~\ref{cor:latin}]
    Use three disjoint sets of $n$ vertices for the rows, columns and symbols, and associate to each cell $(i,j)$ the edge consisting of its row, column and symbol $L(i,j)$.
    The resulting auxiliary hypergraph $\calH$ is $n$-regular, linear and $3$-uniform on $3n$ vertices.
    Its matchings are in bijection with $\calT(L)$, and $|\calU(\calH)|=3(n-|T|)$.
    Theorem~\ref{thm:linearcount} at activity one gives the count, and Theorem~\ref{thm:statistics} gives the $L^2$ conclusion.
    For a group multiplication table, the maps $(x,y,xy)\mapsto(ax,yb,axyb)$ act transitively on each of the three vertex classes.
    Every matching covers equally many vertices in the three classes, so all their vacancy probabilities agree and equal $\E[|\calU(\calH)|]/(3n)\sim n^{-1/3}$.
    This completes the proof.
\end{proof}

\begin{proof}[Proof of Corollary~\ref{cor:queens}]
    Use four disjoint copies of $\mathbb Z_n$ for the rows, columns and two types of diagonals.
    Each cell gives an edge through its four indices.
    The resulting auxiliary hypergraph $\calH$ is $n$-regular and $4$-uniform on $4n$ vertices.
    Its maximum pair codegree is at most two: specifying the two diagonal indices leaves at most two solutions for a cell, and every other pair of compatible indices determines at most one cell.
    For all sufficiently large $n$, the pair-codegree bound satisfies $\Delta_2\le n^{1/8}$.
    The degree is polynomial in the order, so Corollary~\ref{cor:pairrange} with $\gamma=1/8$ gives the counting estimate and the statistical conclusion.
    Since $|\calU(\calH)|=4(n-Q)$, the result follows.
    Translations of the board act transitively on each row, column and diagonal class.
    Since each matching occupies equally many vertices in each class, every vacancy probability equals $\E[|\calU(\calH)|]/(4n)\sim n^{-1/4}$, for odd and even $n$ alike.
    This completes the proof.
\end{proof}

\subsection{Steiner systems and partial designs}\label{subsec:designproofs}

\begin{proof}[Proof of Corollary~\ref{cor:steiner}]
    We have $d=\Theta(n^{t-1})$, by~\eqref{eq:steiner-degree} for a Steiner system and by assumption in the more general case, and $t-1<k/2$.
    Counting $t$-sets extending a fixed $r$-set gives $\Delta_r(H)=O(n^{t-r})$ for $2\le r<t$ and $\Delta_r(H)=O(1)$ for $r\ge t$.
    Hence
    \[
        \Gamma(H)=O(n^{-1/2}+d^{-1/(k-1)})=d^{-1/k}n^{-\Omega(1)}.
    \]
    For $\pi=(a_1,\ldots,a_m)$, choose the smallest $h$ for which $a_1+\cdots+a_h\ge t$.
    Specifying $h$ matching edges and their intersection vertices leaves $O(1)$ candidate edges, since their union contains at least $t$ specified vertices.
    Thus $N_\pi(P)=O(|P|^h)$, and hence $b_\pi(H)=O(n^{h-1})$.
    Since the parts are in decreasing order,
    \[
        \frac{(h-1)k}{m}\le a_1+\cdots+a_{h-1}\le t-1,
    \]
    where the sum is zero if $h=1$.
    The assumption $\gcd(k,t-1)=1$ and the inequality $m<k$ imply that $(t-1)m/k$ is not an integer.
    Therefore $h-1\le((t-1)m-1)/k$, and
    \[
        b_\pi(H)=O(d^{m/k}n^{-1/k}).
    \]
    This estimate for $b_\pi(H)$, together with the preceding bound $\Gamma(H)=d^{-1/k}n^{-\Omega(1)}$, verifies both conditions in~\eqref{eq:admissibility} for a sufficiently small fixed $\eps>0$.
    Applying Theorem~\ref{thm:weightedcount} and Theorem~\ref{thm:statistics} gives the conclusion.
    This completes the proof.
\end{proof}

\begin{proof}[Proof of Theorem~\ref{thm:partialdesigns}]
    Let $\calH$ have the $t$-subsets of $[n]$ as vertices, with an edge consisting of the $t$-shadow of each $k$-subset.
    Its order, uniformity and degree are $N,K,D$, and its matchings correspond bijectively to the partial systems.
    For $k=t+1$, it is linear with degree polynomial in its order, so Theorem~\ref{thm:linearcount} and Theorem~\ref{thm:statistics} apply.
    We check Theorem~\ref{thm:counting} for $k=t+2,t+3$.

    For this purpose, fix a partial system $\calP$ and a candidate $k$-set $B\notin\calP$.
    A block $C\in\calP$ contributes $\binom{|B\cap C|}{t}$ vertices to the intersection pattern of the edge corresponding to $B$.
    Two blocks of $\calP$ intersect in at most $t-1$ points.
    We need only consider candidates whose entire $t$-shadow is covered by $\calP$.

    \begin{claim}\label{clm:set-design-two}
        For $k=t+2$, the auxiliary hypergraph satisfies the hypotheses of Theorem~\ref{thm:counting} with a fixed positive polynomial margin.
    \end{claim}
    \begin{claimproof}[Proof of Claim~\ref{clm:set-design-two}]
    Every nonsingleton intersection has size $t+1$, arising from $|B\cap C|=t+1$.
    There is at most one such intersection, since two $(t+1)$-subsets of $B$ intersect in at least $t$ points.
    Choosing $C\in\calP$, then $B\cap C$, and then the remaining point of $B$ gives
    \[
        b_{(t+1,1,\ldots,1)}(\calH)=O(n).
    \]
    Also $\Delta_j(\calH)=O(n)$ for $2\le j\le t+1$, and $\Delta_j(\calH)=O(1)$ for $j>t+1$: more than $t+1$ distinct $t$-subsets span at least $t+2$ points.
    Therefore
    \[
        \Gamma(\calH)=O\bigl(n^{-1/2}+n^{-1/t}+n^{-2/(K-1)}\bigr).
    \]
    Each exponent is strictly larger than $2/K$, since $K=\binom{t+2}2>2t$ and $K>4$.
    The only pattern in $\Pi_K$ that occurs has $m=K-t$ parts, and
    \[
        \frac{b_\pi(\calH)}{D^{m/K}}=O(n^{-1+2t/K})=n^{-\Omega(1)}.
    \]

    Since $D=\Theta(n^2)$ and $N=\Theta(n^t)$, these strict exponent inequalities verify~\eqref{eq:admissibility}, as well as the polynomial degree condition, for a sufficiently small fixed exponent in $N$.
    \end{claimproof}

    \begin{claim}\label{clm:set-design-three}
        For $k=t+3$, the auxiliary hypergraph also satisfies the hypotheses of Theorem~\ref{thm:counting} with a fixed positive polynomial margin.
    \end{claim}
    \begin{claimproof}[Proof of Claim~\ref{clm:set-design-three}]
    We have $K=\binom{t+3}3$; write $c=\binom{t+2}t$ for the larger nonsingleton intersection size.
    A nonsingleton part has size $t+1$ or $c$.
    A part of size $c$ is the only nonsingleton part: a $(t+2)$-set and a $(t+1)$-set in $B$ meet in at least $t$ points.
    Fixing the matching block which gives the part of size $c$ shows that its coefficient is $O(n)$.

    If there are $h$ parts of size $t+1$, the corresponding $(t+1)$-subsets of $B$ have pairwise disjoint two-point complements.
    Thus $h\le\lfloor(t+3)/2\rfloor$.
    For $h=1$, the coefficient is $O(n^2)$.
    For $h\ge2$, we obtain the sharper bound $O(n)$ as follows.
    Fix a block $C_1\in\calP$ and its $(t+1)$-set $I_1\subseteq B$.
    A second such intersection $I_2$ meets $I_1$ in exactly $t-1$ points.
    There are only constantly many choices for this intersection inside $C_1$.
    For a fixed $(t-1)$-set $U$, each block through $U$ contains $k-t+1$ of its $n-t+1$ extensions to a $t$-set.
    These extensions are disjoint between blocks, so at most $(n-t+1)/(k-t+1)=O(n)$ blocks of $\calP$ contain $U$.
    Choose $C_2$ and $I_2\subseteq C_2$; then $B=I_1\cup I_2$ is determined.
    There are therefore $O(|\calP|n)$ candidates.

    If $j$ distinct $t$-sets span $u$ points, their common codegree is $\binom{n-u}{k-u}$, interpreted as zero when $u>k$.
    More than $t+1$ such sets span at least $t+2$ points, and more than $c$ span at least $t+3$ points.
    Thus
    \[
        \Delta_j(\calH)\le
        \begin{cases}
            O(n^2),&2\le j\le t+1,\\
            O(n),&t+1<j\le c,\\
            O(1),&j>c.
        \end{cases}
    \]
    Substituting these codegree bounds into~\eqref{eq:Gamma} gives
    \begin{equation}\label{eq:set-design-three-gamma}
        \Gamma(\calH)=O\bigl(n^{-1/2}+n^{-1/t}+n^{-2/(c-1)}+n^{-3/(K-1)}\bigr).
    \end{equation}
    Every exponent of $n^{-1}$ in~\eqref{eq:set-design-three-gamma} is strictly larger than $3/K$.
    The coefficient estimates above give
    \begin{equation}\label{eq:set-design-three-patterns}
        \frac{b_\pi(\calH)}{D^{m(\pi)/K}}=
        \begin{cases}
            O\bigl(n^{-2+3(c-1)/K}\bigr),&\text{one part of size }c,\\
            O\bigl(n^{-1+3t/K}\bigr),&\text{one part of size }t+1,\\
            O\bigl(n^{-2+3ht/K}\bigr),&\text{$h\ge2$ parts of size }t+1.
        \end{cases}
    \end{equation}
    The exponents of $n$ in~\eqref{eq:set-design-three-patterns} are negative: $K>6$, $K>3t$, $2K>3(c-1)$, and
    \[
        2K>3t\lfloor(t+3)/2\rfloor.
    \]
    The last inequality follows from
    \[
        2(t+2)(t+1)-9t=2t^2-3t+4>0
    \]
    and $\lfloor(t+3)/2\rfloor\le(t+3)/2$.

    Since $D=\Theta(n^3)$ and $N=\Theta(n^t)$, estimates~\eqref{eq:set-design-three-gamma} and~\eqref{eq:set-design-three-patterns} verify~\eqref{eq:admissibility} with a fixed margin, and the degree is polynomial in $N$.
    \end{claimproof}

    Claims~\ref{clm:set-design-two} and~\ref{clm:set-design-three} establish the hypotheses in both nonlinear cases.
    Apply Theorem~\ref{thm:statistics} to obtain the $L^2$ conclusion.
    The identity $|\calU(\calH)|=N-K|\calS|$ translates this to the number of blocks, and the counting estimate at activity one gives the total count.
    This completes the proof.
\end{proof}

\begin{proof}[Proof of Theorem~\ref{thm:subspaces}]
    Take the $t$-dimensional subspaces as vertices, and for each $(t+s)$-dimensional subspace take its collection of $t$-dimensional subspaces as an edge.
    The resulting hypergraph $\calH$ is $K$-uniform and $D$-regular on $N$ vertices.
    Its matchings correspond bijectively to $\mathfrak D_q(h,t+s,t)$.
    Since $q,t,s$ are fixed, the auxiliary hypergraph has $N=\Theta(q^{ht})$ vertices and degree $D=\Theta(q^{hs})$.
    For $s=1$, distinct blocks intersect in dimension at most $t$, so $\calH$ is linear and Theorem~\ref{thm:linearcount} and Theorem~\ref{thm:statistics} apply.
    We verify Theorem~\ref{thm:counting} for the remaining cases.

    Fix a matching $\calP$ and a candidate block $B\notin\calP$ whose $t$-dimensional subspaces are all covered by $\calP$.
    Distinct blocks of $\calP$ intersect in dimension at most $t-1$.
    A block $C\in\calP$ contributes exactly $\qbinom{\dim(B\cap C)}t$ vertices to the intersection pattern.
    We will also use the following codegree observation: specified $t$-dimensional subspaces with span of dimension $\ell$ lie in exactly $\qbinom{h-\ell}{t+s-\ell}$ candidate blocks when $\ell\le t+s$, and in none otherwise.

    \begin{claim}\label{clm:subspace-two}
        For $s=2$ and $t\ge1$, the auxiliary hypergraph satisfies the hypotheses of Theorem~\ref{thm:counting} with a fixed positive polynomial margin.
    \end{claim}
    \begin{claimproof}[Proof of Claim~\ref{clm:subspace-two}]
    Let $a=\qbinom{t+1}t$.
    A nonsingleton part has size $a$, and there is at most one such part: two $(t+1)$-dimensional subspaces of $B$ intersect in dimension at least $t$.
    Fix the matching block and its $(t+1)$-dimensional intersection with $B$.
    There are $O(q^h)$ choices for $B$, giving $b_{(a,1,\ldots,1)}(\calH)=O(q^h)$.
    Furthermore, $\Delta_j(\calH)=O(q^h)$ for $2\le j\le a$, while $\Delta_j(\calH)=O(1)$ for $j>a$.
    Consequently,
    \[
        \Gamma(\calH)=O\bigl(q^{-h/2}+q^{-h/(a-1)}+q^{-2h/(K-1)}\bigr).
    \]
    The Gaussian binomial formula gives
    \[
        \frac Ka=\frac{q^{t+2}-1}{q^2-1}>q\ge2.
    \]
    Also $a\ge q+1\ge3$, so $K>4$ and $K>2(a-1)$, including when $t=1$.
    Every exponent of $q^{-h}$ in the codegree bound is therefore strictly greater than $2/K$.
    The only pattern in $\Pi_K$ that occurs has $m=K-a+1$ parts, and
    \[
        \frac{b_\pi(\calH)}{D^{m/K}}=O(q^{h(-1+2(a-1)/K)})=q^{-\Omega(h)}.
    \]

    Since $N=\Theta(q^{ht})$ and $D=\Theta(q^{2h})$, these bounds give~\eqref{eq:admissibility} with a fixed margin in $N$ and a polynomial degree lower bound.
    \end{claimproof}

    \begin{claim}\label{clm:subspace-three}
        For $s=3$ and $t\ge2$, the auxiliary hypergraph satisfies the hypotheses of Theorem~\ref{thm:counting} with a fixed positive polynomial margin.
    \end{claim}
    \begin{claimproof}[Proof of Claim~\ref{clm:subspace-three}]
    The two possible nonsingleton intersection sizes are $a=\qbinom{t+1}t$ and $c=\qbinom{t+2}t$.
    A $c$-part is the only nonsingleton part, by the dimension formula, and has coefficient $O(q^h)$.
    A single $a$-part has coefficient $O(q^{2h})$.
    If there are $j\ge2$ parts of size $a$, the corresponding $(t+1)$-dimensional intersections have pairwise intersection dimension exactly $t-1$ and span $B$ in pairs.
    Their annihilators in the dual of $B$ are therefore two-dimensional subspaces with pairwise trivial intersections.
    Counting their nonzero vectors gives
    \[
        j\le J:=\left\lfloor\frac{q^{t+3}-1}{q^2-1}\right\rfloor.
    \]
    To bound their coefficient, fix the first matching block and its intersection with $B$, and then choose the $(t-1)$-dimensional intersection with a second such space.
    Only $O(q^h)$ blocks of $\calP$ contain this $(t-1)$-space: each uses a fixed positive number of the $O(q^h)$ $t$-spaces through it, and these $t$-spaces are distinct between matching blocks.
    After choosing the second block and its intersection with $B$, the candidate $B$ is determined as their span.
    Thus every pattern with $j\ge2$ parts of size $a$ has coefficient $O(q^h)$.

    More than $a$ distinct $t$-spaces span dimension at least $t+2$, and more than $c$ span dimension at least $t+3$.
    The codegree observation therefore gives
    \[
        \Delta_j(\calH)\le
        \begin{cases}
            O(q^{2h}),&2\le j\le a,\\
            O(q^h),&a<j\le c,\\
            O(1),&j>c.
        \end{cases}
    \]
    Consequently,
    \begin{equation}\label{eq:subspace-three-gamma}
        \Gamma(\calH)=O\bigl(q^{-h/2}+q^{-h/(a-1)}+q^{-2h/(c-1)}+q^{-3h/(K-1)}\bigr).
    \end{equation}
    The coefficient estimates give
    \begin{equation}\label{eq:subspace-three-patterns}
        \frac{b_\pi(\calH)}{D^{m(\pi)/K}}=
        \begin{cases}
            O\bigl(q^{h(-2+3(c-1)/K)}\bigr),&\text{one part of size }c,\\
            O\bigl(q^{h(-1+3(a-1)/K)}\bigr),&\text{one part of size }a,\\
            O\bigl(q^{h(-2+3j(a-1)/K)}\bigr),&\text{$j\ge2$ parts of size }a.
        \end{cases}
    \end{equation}
    To verify~\eqref{eq:admissibility}, we check that every exponent of $q^{-h}$ in~\eqref{eq:subspace-three-gamma} exceeds $3/K$ and every exponent of $q^h$ in~\eqref{eq:subspace-three-patterns} is negative.
    Since $t\ge2$, the Gaussian binomial formula gives
    \[
        \frac Kc=\frac{q^{t+3}-1}{q^3-1}>q^2
    \]
    and
    \[
        \frac ca=\frac{q^{t+2}-1}{q^2-1}>q^2.
    \]
    Hence $K>6$, $K>3(a-1)$ and $2K>3(c-1)$.
    For the patterns with several $a$-parts, the bound on $J$ gives
    \[
        \frac{2K}{J(a-1)}\ge\frac{2(q^{t+2}-1)(q^{t+1}-1)}{q(q^3-1)(q^t-1)}>2q\ge4>3.
    \]
    For the strict inequality above, use $(q^{t+2}-1)/(q^t-1)>q^2$ and $(q^{t+1}-1)/(q^3-1)\ge1$.
    Since $N=\Theta(q^{ht})$ and $D=\Theta(q^{3h})$, these strict margins in~\eqref{eq:subspace-three-gamma} and~\eqref{eq:subspace-three-patterns} verify~\eqref{eq:admissibility} and the degree condition.
    \end{claimproof}

    Claims~\ref{clm:subspace-two} and~\ref{clm:subspace-three} complete the verification for the nonlinear cases.
    Theorem~\ref{thm:statistics} gives the $L^2$ conclusion, and $|\calU(\calH)|=N-K|\calS|$ translates it to the number of blocks.
    The counting estimate at activity one gives the total count.
    This completes the proof.
\end{proof}

\subsection{Matroids of fixed rank}\label{subsec:matroidproofs}

For the paving upper bound, deleting the vertices covered by large hyperplanes leaves a linear hypergraph that need not be regular.
We use the following upper bound.

\begin{proposition}[A maximum-degree upper bound]\label{prop:maxdegree}
    For each fixed $k\ge3$, as the positive integer $d\to\infty$, every $n$-vertex linear $k$-uniform hypergraph $F$ with maximum degree at most $d$ satisfies
    \[
        \log|\calM(F)|\le\frac nk(\log d-k+1)+nd^{-1/k}
        +O_k\bigl(nd^{-1/(k-1)}\bigr).
    \]
    The estimate is uniform over $n$ and $F$.
\end{proposition}
All the tools needed for the proof have been developed in Sections~\ref{sec:patterns} and~\ref{sec:counting}.
The proof follows the same entropy argument, keeping track of vertices whose degrees are less than $d$.
We give the details in Appendix~\ref{sec:maxdegreeproof}.

\begin{proof}[Proof of Theorem~\ref{thm:matroids}]
    The auxiliary hypergraph $\calH$ for partial $S(n,r,r-1)$ systems has $N=\binom n{r-1}$ vertices and degree $D=n-r+1$.
    Applying Theorem~\ref{thm:partialdesigns} with $k=r$ and $t=r-1$ gives
    \begin{equation}\label{eq:sparse-paving-count}
        \log s(n,r)=\frac Nr(\log D-r+1)+(1+o(1))ND^{-1/r}.
    \end{equation}

    To bound $p(n,r)$, fix the family $\mathcal L$ of hyperplanes of size at least $r+1$.
    Their $(r-1)$-shadows are pairwise disjoint, so they occupy
    \begin{equation}\label{eq:large-hyperplane-shadow}
        w=\sum_{L\in\mathcal L}\binom{|L|}{r-1}
    \end{equation}
    vertices of the auxiliary hypergraph.
    Deleting these vertices leaves a linear hypergraph $F$ on $N-w$ vertices with maximum degree at most $D$.
    The remaining hyperplanes of size $r$ form a matching in $F$.
    Together with $\mathcal L$, this matching determines the paving matroid: every $(r-1)$-set not contained in a member of $\mathcal L$ or a block of the matching is itself a hyperplane.
    Figure~\ref{fig:matroidshadow} illustrates this deletion and the possible loss of regularity.

    \begin{figure}
        \centering
        \singlespacing
        \matchingfigurestyles
        \resizebox{\linewidth}{!}{%
\begin{tikzpicture}[x=1cm,y=1cm]
\node[panel] at (1.350,4.600) {(a) A large hyperplane};
\node[panel] at (6.150,4.600) {(b) Delete its shadow};
\node[panel] at (11.750,4.600) {(c) Match the remaining vertices};
\node[] at (1.350,3.920) {$r=4,\ n=7$};
\draw[draw=MOlive,fill=MOlive!7,line width=.9pt,rounded corners=5pt]
 (1.35,3.22) -- (.13,2.33) -- (.60,.63) -- (2.10,.63) -- (2.57,2.33) -- cycle;
\node[vertex,fill=MOlive] at (1.350,3.000) {};
\node[] at (1.350,2.710) {1};
\node[vertex,fill=MOlive] at (0.400,2.270) {};
\node[] at (0.400,1.980) {2};
\node[vertex,fill=MOlive] at (0.770,1.180) {};
\node[] at (0.770,0.890) {3};
\node[vertex,fill=MOlive] at (1.930,1.180) {};
\node[] at (1.930,0.890) {4};
\node[vertex,fill=MOlive] at (2.300,2.270) {};
\node[] at (2.300,1.980) {5};
\node[] at (1.350,0.160) {$L=\{1,2,3,4,5\}$};
\draw[arr] (2.88,2.05) -- (3.90,2.05);
\node[] at (6.150,3.920) {$\partial_3L=\binom L3$};
\draw[draw=MRed,fill=MRed!4,line width=.8pt,rounded corners=6pt] (4.20,.42) rectangle (8.10,3.43);
\node[circle,draw=MRed,fill=white,line width=.7pt,inner sep=1.8pt] at (4.800,3.080) {};
\node[text=MRed] at (5.270,3.080) {123};
\node[circle,draw=MRed,fill=white,line width=.7pt,inner sep=1.8pt] at (4.800,2.530) {};
\node[text=MRed] at (5.270,2.530) {124};
\node[circle,draw=MRed,fill=white,line width=.7pt,inner sep=1.8pt] at (4.800,1.980) {};
\node[text=MRed] at (5.270,1.980) {125};
\node[circle,draw=MRed,fill=white,line width=.7pt,inner sep=1.8pt] at (4.800,1.430) {};
\node[text=MRed] at (5.270,1.430) {134};
\node[circle,draw=MRed,fill=white,line width=.7pt,inner sep=1.8pt] at (4.800,0.880) {};
\node[text=MRed] at (5.270,0.880) {135};
\node[circle,draw=MRed,fill=white,line width=.7pt,inner sep=1.8pt] at (6.700,3.080) {};
\node[text=MRed] at (7.170,3.080) {145};
\node[circle,draw=MRed,fill=white,line width=.7pt,inner sep=1.8pt] at (6.700,2.530) {};
\node[text=MRed] at (7.170,2.530) {234};
\node[circle,draw=MRed,fill=white,line width=.7pt,inner sep=1.8pt] at (6.700,1.980) {};
\node[text=MRed] at (7.170,1.980) {235};
\node[circle,draw=MRed,fill=white,line width=.7pt,inner sep=1.8pt] at (6.700,1.430) {};
\node[text=MRed] at (7.170,1.430) {245};
\node[circle,draw=MRed,fill=white,line width=.7pt,inner sep=1.8pt] at (6.700,0.880) {};
\node[text=MRed] at (7.170,0.880) {345};
\node[] at (6.150,-0.120) {$w=\binom53=10$};
\draw[arr] (8.40,2.05) -- (9.20,2.05);
\draw[source] (9.48,1.30) rectangle (13.65,3.01);
\node[vertex,fill=MBlue] at (10.250,2.620) {};
\node[] at (10.250,2.330) {126};
\node[vertex,fill=MBlue] at (12.850,2.620) {};
\node[] at (12.850,2.330) {127};
\node[vertex,fill=MBlue] at (10.250,1.870) {};
\node[] at (10.250,1.580) {167};
\node[vertex,fill=MBlue] at (12.850,1.870) {};
\node[] at (12.850,1.580) {267};
\draw[source] (9.48,-0.77) rectangle (13.65,0.94);
\node[vertex,fill=MBlue] at (10.250,0.550) {};
\node[] at (10.250,0.260) {346};
\node[vertex,fill=MBlue] at (12.850,0.550) {};
\node[] at (12.850,0.260) {347};
\node[vertex,fill=MBlue] at (10.250,-0.200) {};
\node[] at (10.250,-0.490) {367};
\node[vertex,fill=MBlue] at (12.850,-0.200) {};
\node[] at (12.850,-0.490) {467};
\node[note] at (11.600,3.600) {Two surviving $4$-edges};
\node[text=black!65] at (11.600,-1.220) {plus $17$ further vertices};
\node[font=\normalsize] at (2.200,-1.920) {$F=\mathcal H-\partial_3L$};
\node[font=\normalsize] at (6.300,-1.920) {$|V(F)|=35-10=25$};
\node[] at (11.750,-1.920) {$d_F(126)=1$};
\node[] at (11.750,-2.470) {$d_F(167)=4$};
\node[font=\normalsize] at (4.400,-2.650) {$\Delta(F)\le4$, but $F$ is no longer regular.};
\end{tikzpicture}%
}
        \caption{Deleting the shadow of a large hyperplane for $r=4$ and $n=7$. The hyperplane $L=\{1,2,3,4,5\}$ removes its ten triples from $\calH$. The two displayed edges of $F$ are the shadows of $\{1,2,6,7\}$ and $\{3,4,6,7\}$; triple labels abbreviate sets. The residual degrees $d_F(126)=1$ and $d_F(167)=4$ show why a maximum-degree bound is needed.}
        \label{fig:matroidshadow}
    \end{figure}
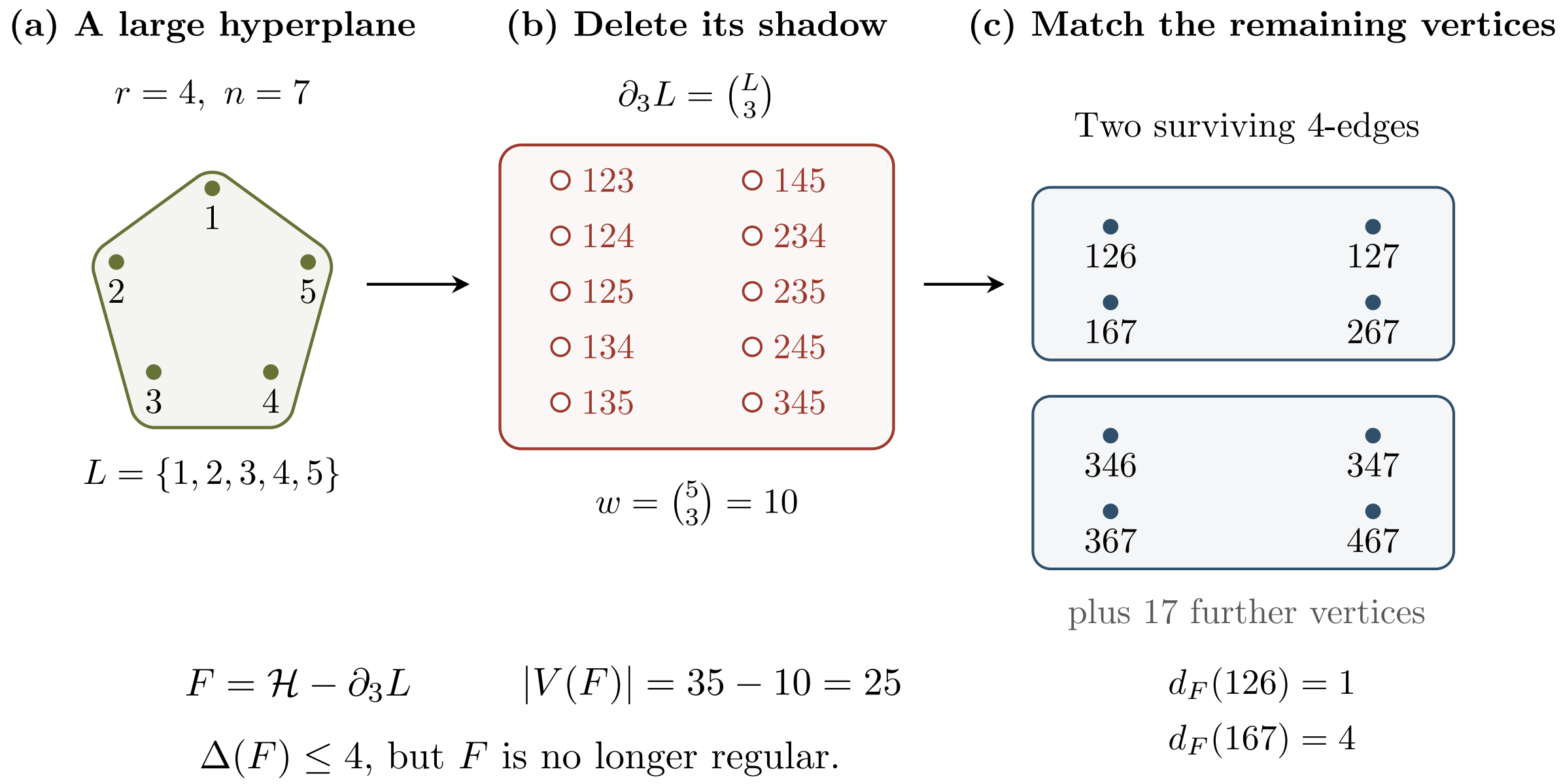

    For the residual hypergraph $F$ on $N-w$ vertices, with $w$ as in~\eqref{eq:large-hyperplane-shadow}, Proposition~\ref{prop:maxdegree} gives
    \begin{equation}\label{eq:paving-residual-count}
        \log|\calM(F)|\le\frac{N-w}{r}(\log D-r+1)
        +ND^{-1/r}+O_r\bigl(ND^{-1/(r-1)}\bigr).
    \end{equation}
    Exponentiating~\eqref{eq:paving-residual-count} and expanding $w$ by~\eqref{eq:large-hyperplane-shadow} gives a factor $(e^{r-1}/D)^{\binom{|L|}{r-1}/r}$ for each $L\in\mathcal L$.
    Summing over $\mathcal L$ and omitting the compatibility restrictions between its members, we obtain
    \begin{align}
        \log p(n,r)
        &\le\frac Nr(\log D-r+1)+ND^{-1/r}+O_r\bigl(ND^{-1/(r-1)}\bigr)\notag\\
        &\quad+\sum_{b=r+1}^{n}\binom nb
        \log\left(1+\left(\frac{e^{r-1}}D\right)^{\binom b{r-1}/r}\right)\notag\\
        &\le\frac Nr(\log D-r+1)+ND^{-1/r}+O_r\bigl(ND^{-1/(r-1)}\bigr)
        +\sum_{b=r+1}^{n}\binom nb\left(\frac{e^{r-1}}D\right)^{\binom b{r-1}/r}.
    \label{eq:paving-count-upper}
    \end{align}
    The following claim bounds the final sum in~\eqref{eq:paving-count-upper}, which accounts for the additional choices of large hyperplanes.
    \begin{claim}\label{clm:large-hyperplane-cost}
        For fixed $r\ge4$ and $D=n-r+1$,
        \[
            \sum_{b=r+1}^{n}\binom nb
            \left(\frac{e^{r-1}}D\right)^{\binom b{r-1}/r}
            =O_r(n^{(r+1)/2}).
        \]
    \end{claim}
    \begin{claimproof}[Proof of Claim~\ref{clm:large-hyperplane-cost}]
        The first term has order
        \[
            n^{r+1-\binom{r+1}{r-1}/r}=n^{(r+1)/2}.
        \]
        For $r+1\le b<n$, the ratio of the next term to the current term is
        \[
            \begin{aligned}
                \frac{n-b}{b+1}\left(\frac{e^{r-1}}D\right)^{\binom b{r-2}/r}
                &\le n\left(\frac{e^{r-1}}D\right)^{(r^2-1)/6}\\
                &=O_r(n^{1-(r^2-1)/6})=o(1).
            \end{aligned}
        \]
        Here $e^{r-1}/D<1$ for sufficiently large $n$, and $\binom b{r-2}/r\ge(r^2-1)/6$.
        The ratios are uniformly at most $1/2$ for large $n$, so the sum is at most twice its first term.
    \end{claimproof}

    Since $r\ge4$,
    \[
        r-1-\frac1r-\frac{r+1}{2}=\frac{r-3}{2}-\frac1r>0.
    \]
    Thus $n^{(r+1)/2}=o(ND^{-1/r})$, and $D^{-1/(r-1)}=o(D^{-1/r})$ controls the other error.
    Claim~\ref{clm:large-hyperplane-cost} in~\eqref{eq:paving-count-upper}, together with $p(n,r)\ge s(n,r)$ and~\eqref{eq:sparse-paving-count}, gives the formula for paving matroids.

    Finally, Pendavingh and van der Pol~\cite[Theorem~3]{PvdP} proved that
    \[
        \log p(n,r)\le\log m(n,r)
        \le\left(1+\frac{r+o(1)}{n-r+1}\right)\log p(n,r).
    \]
    The additional error is $O_r(n^{r-2}\log n)=o(ND^{-1/r})$, so the same formula holds for $m(n,r)$.
    Expanding $N$ and $D$, we have
    \[
        \frac Nr(\log D-r+1)
        =\frac{n^{r-1}}{r!}(\log n-r+1)+O_r(n^{r-2}\log n)
    \]
    and
    \[
        ND^{-1/r}=\frac{n^{r-1-1/r}}{(r-1)!}(1+O_r(n^{-1})).
    \]
    Substitute these expansions into~\eqref{eq:sparse-paving-count}. Since $n^{r-2}\log n=o(n^{r-1-1/r})$, this gives the stated formula for $s(n,r)$; the comparisons with $p(n,r)$ and $m(n,r)$ above have errors of the same smaller order and give the formula for all three classes.
    This completes the proof.
\end{proof}


\section{Concluding remarks}\label{sec:concluding}

We determined the correction $nd^{-1/k}$ to the logarithmic matching count under a criterion on higher codegrees and intersection patterns.
The weighted estimate yields the sharp mean vacancy, concentration and large deviations, and gives refined counts in the applications, including matroids of each fixed rank at least four.
The main remaining questions concern sharper fluctuations, weaker hypotheses and the mixed hyperplane sizes in rank three.

\paragraph{Sharp variance and fluctuations.}
Theorem~\ref{thm:statistics} gives concentration on the deficit scale:
\[
    \Var_\lambda\bigl[|\calU_\lambda(H)|\bigr]=o\bigl(n^2(\lambda d)^{-2/k}\bigr).
\]
This does not determine the variance on the expected scale $n(\lambda d)^{-1/k}$.
The natural target is
\begin{equation}\label{eq:varianceasymptotic}
    \begin{gathered}
        \Var_\lambda\bigl[|\calU_\lambda(H)|\bigr]=(1+o(1))n(\lambda d)^{-1/k},\\
        \Var_\lambda[|M|]=(1+o(1))\frac{n(\lambda d)^{-1/k}}{k^2}.
    \end{gathered}
\end{equation}
At activity one in the linear case, this is the asymptotic formula proposed by Kahn and Kim~\cite{KK}.

\begin{question}\label{ques:variance}
    For fixed $k\ge3$ and $\eps>0$, is~\eqref{eq:varianceasymptotic} valid uniformly over all $n$-vertex $d$-regular linear $k$-uniform hypergraphs with $d\ge n^\eps$, and over $\lambda$ in each fixed compact interval of $(0,\infty)$?
\end{question}

Even $\Var_\lambda\bigl[|\calU_\lambda(H)|\bigr]=O(n(\lambda d)^{-1/k})$ would improve our estimate.
One can also ask whether
\[
    \frac{|\calU_\lambda(H)|-\E_\lambda[|\calU_\lambda(H)|]}{\sqrt{n(\lambda d)^{-1/k}}}
    \ \xrightarrow{\mathrm d}\ \mathcal N(0,1).
\]
Centering at the actual mean is essential: the present expectation estimate does not locate it to within $o(\sqrt{n(\lambda d)^{-1/k}})$.
The large-deviation rate alone does not imply this limit, a Poisson approximation, or the sharp variance asymptotic.

\paragraph{Pointwise vacancy.}
For the regular linear hypergraphs covered here, does
\[
    \max_{v\in V(H)}|d^{1/k}p_H(v)-1| \longrightarrow 0
\]
hold under $d\ge n^\eps$?
Vertex transitivity makes this a consequence of the averaged result, but regularity alone does not make the vertices interchangeable.
The author's linear counterexamples~\cite{Lee} have orders very large compared with their degrees.
Appendix~\ref{sec:counterexample} disproves the broader small-codegree formulation even at polynomial degree, but does not settle this restricted linear question.

\paragraph{The graph case.}
For a simple regular graph, $\Delta_2=1$, so the pair-codegree contribution $\sqrt{\Delta_2/d}=d^{-1/2}$ is already of the same order as the average vacancy probability.
An error of order $n\sqrt{\Delta_2/d}$ therefore does not resolve the correction to the logarithmic count.
There is a separate sharp upper bound: Davies, Jenssen, Perkins and Roberts~\cite{DJPR} give $\log Z_G(\lambda)\le n\log Z_{K_{d,d}}(\lambda)/(2d)$ for every $n$-vertex $d$-regular graph.
As $d\to\infty$, the explicit matching polynomial of $K_{d,d}$ and Stirling's formula yield, uniformly on compact intervals of positive activities,
\[
    \log Z_G(\lambda)\le\frac n2\bigl(\log(\lambda d)-1\bigr)
    +n(\lambda d)^{-1/2}+O(n\log d/d).
\]
This illustrates why the graph case is better treated with its available graph-specific estimates.

\paragraph{The relation between degree and order.}
The refined counting formula retains $d\ge n^\eps$ and fixed polynomial margins to absorb the logarithmic error from Theorem~\ref{thm:greedycount} below the correction scale.
This is a sufficient hypothesis, not a proposed necessary threshold.
The more basic question is whether the averaged conclusion needs any relation between degree and order.
Without such a relation, Gupta~\cite{Gupta} bounds the average vacancy by $O_k(\log\log d/\log d)$.

\begin{question}\label{ques:averagedvacancy}
    For each fixed $k\ge3$, is
    \[
        \E[|\calU(H)|]=(1+o_d(1))nd^{-1/k}
    \]
    valid uniformly over all $n$-vertex $d$-regular linear $k$-uniform hypergraphs as $d\to\infty$, with no restriction on $n$?
\end{question}

The known pointwise counterexamples do not determine this average.
Even a uniform upper bound $\E[|\calU(H)|]=O_k(nd^{-1/k})$ would substantially strengthen the available estimate.
The corresponding count requires precision below the leading exponential scale, since its correction is only $nd^{-1/k}$.

\paragraph{The pair-codegree range.}
Corollary~\ref{cor:pairrange} gives the sharp count and mean vacancy under polynomial degree and $\Delta_2\le d^\gamma$ for every fixed $\gamma<1/k$.
The counterexamples in Appendix~\ref{sec:counterexample} show that polynomial degree and a polynomial gap between $d$ and $\Delta_2(H)$ do not suffice for the lower estimate $\E[|\calU(H)|]\ge(1-o(1))nd^{-1/k}$.
Determining how far the sufficient pair-codegree range can be extended remains open.
The intersection coefficients distinguish hypergraphs with the same maximum pair codegree, so the exponent in Corollary~\ref{cor:pairrange} does not describe the full scope of Theorem~\ref{thm:weightedcount}.

These counterexamples leave open the one-sided upper bound $\E[|\calU(H)|]=O(nd^{-1/k})$ under weaker codegree assumptions.
Appendix~\ref{sec:conditionallifting} gives a precise extension that would improve maximum-matching bounds even for regular linear hypergraphs of polynomial degree.
This reduction retains a polynomial relation between degree and order.

\paragraph{Perfect matchings.}\label{par:perfect-matchings}
Combining~\eqref{eq:perfect-probability} with Theorem~\ref{thm:weightedcount} at $\lambda=1$ shows that, under its hypotheses, the number $\operatorname{pm}(H)$ of perfect matchings satisfies
\[
    \operatorname{pm}(H)\le
    \left(\frac d{e^{k-1}}\right)^{n/k}\exp\bigl(o(nd^{-1/k})\bigr).
\]
A positive lower bound requires additional assumptions, since perfect matchings may not exist, for example when $k\nmid n$.

\paragraph{Steiner systems outside the stated range.}
Corollary~\ref{cor:steiner} assumes $t<k/2+1$ and $\gcd(k,t-1)=1$.
The first condition provides a polynomial margin in the error in the greedy counting estimate, and the second provides one in the pattern counts.
Removing either restriction requires estimates beyond the uniform bounds used here.
For instance, in $S(n,4,3)$,
\[
    \begin{aligned}
        d&=\frac{(n-1)(n-2)}6,\\
        \Delta_2&=\frac{n-2}2,
    \end{aligned}
\]
and
\[
    nd^{-1/4}\sim6^{1/4}\sqrt n.
\]
The error in the greedy counting estimate $n\Gamma(H)(\log n)^L=\Theta(\sqrt n(\log n)^L)$ already exceeds the desired correction, and the elementary bound for the pattern $(2,2)$ has no polynomial margin.
It remains natural to ask for the same sharp count and typical deficit throughout the larger Steiner range.

\paragraph{Further partial designs.}
Which partial $S(n,k,t)$ systems with $k\ge t+4$, and which partial subspace designs with $s\ge4$ or $(t,s)=(1,3)$, have auxiliary hypergraphs satisfying~\eqref{eq:admissibility}?
These are the parameter ranges not covered by Theorems~\ref{thm:partialdesigns} and~\ref{thm:subspaces}.

\paragraph{Rank-three matroids and linear spaces.}\label{par:rankthree}
Theorem~\ref{thm:matroids} leaves open the corresponding refinement for rank three.
Kwan, Sah and Sawhney~\cite[Theorems~1.1 and~1.2]{KSSMatroids} proved
\[
    \log f(n,3)=\frac{n^2}{6}(\log n+\log c)+o(n^2)
\]
for $f\in\{p,m\}$, where $c=\frac{1+\sqrt3}{2}e^{\sqrt3/2-3}$.
Here paving matroids correspond to linear spaces, or equivalently to partitions of the edges of $K_n$ into cliques.
The optimizing profile has both three-point and four-point lines in quantities of order $n^2$~\cite[Lemma~2.5]{KSSMatroids}.
Thus the auxiliary matching problem has edges of sizes three and six already at leading order, and the contribution of the larger hyperplanes cannot be absorbed as in rank at least four.

\begin{question}\label{ques:rankthree}
    Is there a constant $C>0$ such that
    \[
        \log p(n,3)=\frac{n^2}{6}(\log n+\log c)+Cn^{5/3}+o(n^{5/3})?
    \]
    If so, determine $C$.
\end{question}

A formal stationary-point calculation in the bounded-line variational expression of~\cite[Lemma~2.2]{KSSMatroids} suggests the scale $n^{5/3}$.
Indeed, the trial line counts have the form $\sigma_j=\binom nj z^{\binom j2}$ with $z\asymp n^{-1/3}$, so both $\sigma_2$ and $\sigma_5$ have order $n^{5/3}$.
This calculation does not control the counting error on that scale; a refinement must also account for five-point lines and the dependence between line sizes.
The same question for $m(n,3)$ is equivalent, since~\cite[Theorem~3]{PvdP} gives $\log m(n,3)-\log p(n,3)=O(n\log n)=o(n^{5/3})$.


\section*{Statement of AI use}

This research involved a human--AI collaboration using OpenAI's
ChatGPT 5.5 Pro, ChatGPT 5.6 Pro and ChatGPT 6 Pro.
The main research questions and the statements of the main results
were formulated entirely by the author, who also developed the
overall mathematical strategy.
The specific contributions of the AI models were as follows.

\begin{enumerate}[(i)]
    \item \emph{Counterexample.}
        When asked by the author to prove or disprove Conjecture~\ref{conj:Lee}, ChatGPT 5.5 Pro produced
        the counterexample presented in Appendix~\ref{sec:counterexample}.
        This prompted the author to focus on averaged vacancy
        probabilities under suitable degree and codegree conditions.

    \item \emph{Random greedy counting.}
        ChatGPT 5.6 Pro and ChatGPT 6 Pro made substantial
        contributions to the development and proof of the random
        greedy counting argument in Section~\ref{sec:greedy},
        through repeated exchanges guided by the author.
        In particular, all technical calculations in that section were carried out by the models.

    \item \emph{Entropy method.}
        The entropy method based on intersection patterns in
        Section~\ref{sec:patterns} was developed jointly by
        the author, ChatGPT 5.6 Pro and ChatGPT 6 Pro,
        through repeated exchanges under the author's
        mathematical guidance.

    \item \emph{Matroid application.}
        ChatGPT 6 Pro identified the application to counting
        matroids of fixed rank in Section~\ref{subsec:matroids}.

    \item \emph{Literature and presentation.}
        The models assisted with literature searches, writing
        and proofreading.
        ChatGPT 6 Pro generated the TikZ code for all figures
        in the paper.
\end{enumerate}

The author has independently verified all AI-generated mathematical
results and arguments included in this paper, reviewed all other
AI-assisted content, and takes full responsibility for the entire
content of the paper.


\section*{Acknowledgements}
Part of this work was done during the author's visit to the University of Vienna. He gratefully acknowledges the University of Vienna for its hospitality and support.
The author thanks his advisor, Jaehoon Kim, for his valuable comments on the initial draft of this manuscript. He also thanks Donggyu Kim for very helpful discussions on matroids.


\providecommand{\MR}[1]{}
\providecommand{\MRhref}[2]{\href{http://www.ams.org/mathscinet-getitem?mr=#1}{#2}}
\bibliographystyle{amsplain_initials_nobysame}
\bibliography{bibfile}

@article{AKS,
  author = {Alon, Noga and Kim, Jeong Han and Spencer, Joel},
  title = {Nearly perfect matchings in regular simple hypergraphs},
  journal = {Israel J. Math.}, volume = {100}, year = {1997}, pages = {171--187},
  doi = {10.1007/BF02773639}
}

@article{AsratianKuzjurin,
  author = {Asratian, Armen S. and Kuzjurin, Nikolai N.},
  title = {On the number of nearly perfect matchings in almost regular uniform hypergraphs},
  journal = {Discrete Math.}, volume = {207}, year = {1999}, pages = {1--8},
  doi = {10.1016/S0012-365X(99)00127-2}
}

@article{AKPartial,
  author = {Asratian, Armen S. and Kuzjurin, Nikolai N.},
  title = {On the number of partial {Steiner} systems},
  journal = {J. Combin. Des.}, volume = {8}, number = {5},
  year = {2000}, pages = {347--352},
  doi = {10.1002/1520-6610(2000)8:5<347::AID-JCD4>3.0.CO;2-8}
}

@article{BB,
  author = {Bennett, Patrick and Bohman, Tom},
  title = {A natural barrier in random greedy hypergraph matching},
  journal = {Combin. Probab. Comput.}, volume = {28}, number = {6},
  year = {2019}, pages = {816--825}, doi = {10.1017/S0963548319000051},
  eprint = {1210.3581}, archivePrefix = {arXiv}
}

@unpublished{BF,
  author = {Bennett, Patrick and Frieze, Alan},
  title = {Binomial random matroids}, year = {2026},
  note = {Preprint, \href{https://arxiv.org/abs/2603.10293}{arXiv:2603.10293}},
  eprint = {2603.10293}, archivePrefix = {arXiv}
}

@book{BLM,
  author = {Boucheron, St{\'e}phane and Lugosi, G{\'a}bor and Massart, Pascal},
  title = {Concentration Inequalities: A Nonasymptotic Theory of Independence},
  publisher = {Oxford University Press}, year = {2013},
  doi = {10.1093/acprof:oso/9780199535255.001.0001}
}

@article{CGT,
  author = {Carroll, Teena and Galvin, David and Tetali, Prasad},
  title = {Matchings and independent sets of a fixed size in regular graphs},
  journal = {J. Combin. Theory Ser. A}, volume = {116}, number = {7},
  year = {2009}, pages = {1219--1227}, doi = {10.1016/j.jcta.2008.12.008}
}

@article{CFGS,
  author = {Chung, Fan R. K. and Graham, Ronald L. and Frankl, Peter and Shearer, James B.},
  title = {Some intersection theorems for ordered sets and graphs},
  journal = {J. Combin. Theory Ser. A}, volume = {43}, year = {1986}, pages = {23--37},
  doi = {10.1016/0097-3165(86)90019-1}
}

@article{CPST,
  author = {Cohen, Emma and Perkins, Will and Sarantis, Michail and Tetali, Prasad},
  title = {On the number of independent sets in uniform, regular, linear hypergraphs},
  journal = {European J. Combin.}, volume = {99}, year = {2022}, pages = {103401},
  doi = {10.1016/j.ejc.2021.103401}
}

@book{CoverThomas,
  author = {Cover, Thomas M. and Thomas, Joy A.},
  title = {Elements of Information Theory}, edition = {2nd},
  publisher = {John Wiley \& Sons}, year = {2006},
  doi = {10.1002/047174882X}
}

@article{DDK,
  author = {Dai, Tantan and Divoux, Alexander and Kelly, Tom},
  title = {Entropy bounds for perfect matchings in bipartite hypergraphs},
  journal = {Electron. J. Combin.}, volume = {33}, number = {2},
  year = {2026}, pages = {Paper No. P2.20}, doi = {10.37236/14339},
  eprint = {2506.17652}, archivePrefix = {arXiv}
}

@article{DJPR,
  author = {Davies, Ewan and Jenssen, Matthew and Perkins, Will and Roberts, Barnaby},
  title = {Independent sets, matchings, and occupancy fractions},
  journal = {J. Lond. Math. Soc.}, volume = {96}, number = {1}, year = {2017}, pages = {47--66},
  doi = {10.1112/jlms.12056}
}

@article{EH,
  author = {Erd{\H{o}}s, Paul and Hanani, Haim},
  title = {On a limit theorem in combinatorical analysis},
  journal = {Publ. Math. Debrecen}, volume = {10}, year = {1963}, pages = {10--13},
  doi = {10.5486/PMD.1963.10.1-4.02}
}

@article{FR,
  author = {Frankl, Peter and R{\"o}dl, Vojt{\v{e}}ch},
  title = {Near perfect coverings in graphs and hypergraphs},
  journal = {European J. Combin.}, volume = {6}, year = {1985}, pages = {317--326},
  doi = {10.1016/S0195-6698(85)80045-7}
}

@article{Freedman,
  author = {Freedman, David A.},
  title = {On tail probabilities for martingales},
  journal = {Ann. Probab.}, volume = {3}, number = {1},
  year = {1975}, pages = {100--118}, doi = {10.1214/aop/1176996452}
}

@article{GJKKL,
  author = {Glock, Stefan and Joos, Felix and Kim, Jaehoon and K{\"u}hn, Marcus and Lichev, Lyuben},
  title = {Conflict-free hypergraph matchings},
  journal = {J. Lond. Math. Soc.}, volume = {109}, number = {5},
  year = {2024}, pages = {e12899}, doi = {10.1112/jlms.12899},
  eprint = {2205.05564}, archivePrefix = {arXiv}
}

@article{GKLO,
  author = {Glock, Stefan and K{\"u}hn, Daniela and Lo, Allan and Osthus, Deryk},
  title = {The existence of designs via iterative absorption: hypergraph {$F$}-designs for arbitrary {$F$}},
  journal = {Mem. Amer. Math. Soc.}, volume = {284}, number = {1406},
  year = {2023}, pages = {v+131}, doi = {10.1090/memo/1406}
}

@unpublished{GK,
  author = {Gould, Stephen and Kelly, Tom},
  title = {Advancing the {R{\"o}dl} nibble: New bounds on matchings and the list chromatic index of hypergraphs},
  year = {2025},
  note = {Preprint, \href{https://arxiv.org/abs/2511.11375}{arXiv:2511.11375}},
  eprint = {2511.11375}, archivePrefix = {arXiv}
}

@unpublished{Gupta,
  author = {Gupta, Anish},
  title = {A rate for the average vacancy of uniformly random matchings in linear hypergraphs},
  year = {2026},
  note = {Preprint, \href{https://arxiv.org/abs/2609.26133}{arXiv:2609.26133}},
  eprint = {2609.26133}, archivePrefix = {arXiv}
}

@article{HL,
  author = {Heilmann, Ole J. and Lieb, Elliott H.},
  title = {Theory of monomer-dimer systems},
  journal = {Comm. Math. Phys.}, volume = {25}, year = {1972}, pages = {190--232},
  doi = {10.1007/BF01877590}
}

@article{HPP,
  author = {van der Hofstad, Remco and Pendavingh, Rudi and van der Pol, Jorn},
  title = {The number of partial {Steiner} systems and {$d$}-partitions},
  journal = {Adv. Comb.}, year = {2022}, pages = {Paper No. 2, 23},
  doi = {10.19086/aic.32563}, eprint = {1811.11810}, archivePrefix = {arXiv}
}

@article{IK,
  author = {Ilinca, Liviu and Kahn, Jeff},
  title = {Asymptotics of the upper matching conjecture},
  journal = {J. Combin. Theory Ser. A}, volume = {120}, number = {5},
  year = {2013}, pages = {976--983}, doi = {10.1016/j.jcta.2013.01.013},
  eprint = {1205.4342}, archivePrefix = {arXiv}
}

@inproceedings{Kahn95,
  author = {Kahn, Jeff},
  title = {Asymptotics of hypergraph matching, covering and coloring problems},
  booktitle = {Proceedings of the International Congress of Mathematicians, Vol. 1, 2 ({Z{\"u}rich}, 1994)},
  publisher = {Birkh{\"a}user}, address = {Basel}, year = {1995}, pages = {1353--1362},
  doi = {10.1007/978-3-0348-9078-6_130}
}

@incollection{Kahn97,
  author = {Kahn, Jeff},
  title = {On some hypergraph problems of {Paul Erd{\H{o}}s} and the asymptotics of matchings, covers and colorings},
  booktitle = {The Mathematics of {Paul Erd{\H{o}}s}, {I}},
  series = {Algorithms Combin.}, volume = {13},
  publisher = {Springer}, address = {Berlin}, year = {1997}, pages = {345--371},
  doi = {10.1007/978-3-642-60408-9_26}
}

@article{KK,
  author = {Kahn, Jeff and Kim, Jeong Han},
  title = {Random matchings in regular graphs},
  journal = {Combinatorica}, volume = {18}, number = {2}, year = {1998}, pages = {201--226},
  doi = {10.1007/PL00009817}
}

@article{KKMO,
  author = {Kang, Dong Yeap and K{\"u}hn, Daniela and Methuku, Abhishek and Osthus, Deryk},
  title = {New bounds on the size of nearly perfect matchings in almost regular hypergraphs},
  journal = {J. Lond. Math. Soc.}, volume = {108}, number = {5}, year = {2023}, pages = {1701--1746},
  doi = {10.1112/jlms.12792}, eprint = {2010.04183}, archivePrefix = {arXiv}
}

@unpublished{Keevash,
  author = {Keevash, Peter}, title = {The existence of designs},
  year = {2014},
  note = {Preprint, revised 2024, \href{https://arxiv.org/abs/1401.3665}{arXiv:1401.3665}},
  eprint = {1401.3665}, archivePrefix = {arXiv}
}

@article{KeevashCount,
  author = {Keevash, Peter}, title = {Counting designs},
  journal = {J. Eur. Math. Soc.}, volume = {20}, number = {4},
  year = {2018}, pages = {903--927}, doi = {10.4171/JEMS/779},
  eprint = {1504.02909}, archivePrefix = {arXiv}
}

@unpublished{KeevashII,
  author = {Keevash, Peter}, title = {The existence of designs {II}},
  year = {2018},
  note = {Preprint, \href{https://arxiv.org/abs/1802.05900}{arXiv:1802.05900}},
  eprint = {1802.05900}, archivePrefix = {arXiv}
}

@article{KSW,
  author = {Kwan, Matthew and Safavi, Roodabeh and Wang, Yiting},
  title = {Counting perfect matchings in {Dirac} hypergraphs},
  journal = {Combinatorica}, volume = {46},
  year = {2026}, pages = {Paper No. 5, 32}, doi = {10.1007/s00493-025-00194-8},
  eprint = {2408.09589}, archivePrefix = {arXiv}
}

@article{KSSMatroids,
  author = {Kwan, Matthew and Sah, Ashwin and Sawhney, Mehtaab},
  title = {Enumerating matroids and linear spaces},
  journal = {C. R. Math. Acad. Sci. Paris}, volume = {361},
  year = {2023}, pages = {565--575}, doi = {10.5802/crmath.423},
  eprint = {2112.03788}, archivePrefix = {arXiv}
}

@article{KPSY,
  author = {Keevash, Peter and Pokrovskiy, Alexey and Sudakov, Benny and Yepremyan, Liana},
  title = {New bounds for {Ryser}'s conjecture and related problems},
  journal = {Trans. Amer. Math. Soc. Ser. B}, volume = {9}, year = {2022}, pages = {288--321},
  doi = {10.1090/btran/92}
}

@article{Lee,
  author = {Lee, Hyunwoo}, title = {Random matchings in linear hypergraphs},
  journal = {Forum Math. Sigma}, volume = {14}, year = {2026}, pages = {e117},
  doi = {10.1017/fms.2026.10253}, eprint = {2406.06421}, archivePrefix = {arXiv}
}

@article{LinialLuria,
  author = {Linial, Nathan and Luria, Zur},
  title = {An upper bound on the number of {Steiner} triple systems},
  journal = {Random Structures Algorithms}, volume = {43}, number = {4},
  year = {2013}, pages = {399--406}, doi = {10.1002/rsa.20487}
}

@unpublished{Luria,
  author = {Luria, Zur}, title = {New bounds on the number of {$n$}-queens configurations},
  year = {2017},
  note = {Preprint, \href{https://arxiv.org/abs/1705.05225}{arXiv:1705.05225}},
  eprint = {1705.05225}, archivePrefix = {arXiv}
}

@unpublished{Montgomery,
  author = {Montgomery, Richard},
  title = {A proof of the {Ryser--Brualdi--Stein} conjecture for large even {$n$}},
  year = {2023},
  note = {Preprint, \href{https://arxiv.org/abs/2310.19779}{arXiv:2310.19779}},
  eprint = {2310.19779}, archivePrefix = {arXiv}
}

@article{McKayTian,
  author = {McKay, Brendan D. and Tian, Fang},
  title = {Asymptotic enumeration of linear hypergraphs with given number of vertices and edges},
  journal = {Adv. Appl. Math.}, volume = {115},
  year = {2020}, pages = {102000}, doi = {10.1016/j.aam.2020.102000},
  eprint = {1908.06333}, archivePrefix = {arXiv}
}

@article{PvdP,
  author = {Pendavingh, Rudi and van der Pol, Jorn},
  title = {Enumerating matroids of fixed rank},
  journal = {Electron. J. Combin.}, volume = {24}, number = {1},
  year = {2017}, pages = {Paper No. P1.8, 28}, doi = {10.37236/5894},
  eprint = {1512.06655}, archivePrefix = {arXiv}
}

@article{PS,
  author = {Pippenger, Nicholas and Spencer, Joel},
  title = {Asymptotic behavior of the chromatic index for hypergraphs},
  journal = {J. Combin. Theory Ser. A}, volume = {51}, year = {1989}, pages = {24--42},
  doi = {10.1016/0097-3165(89)90074-5}
}

@article{Radhakrishnan,
  author = {Radhakrishnan, Jaikumar}, title = {An entropy proof of {Bregman}'s theorem},
  journal = {J. Combin. Theory Ser. A}, volume = {77}, year = {1997}, pages = {161--164},
  doi = {10.1006/jcta.1996.2727}
}

@article{Rodl,
  author = {R{\"o}dl, Vojt{\v{e}}ch}, title = {On a packing and covering problem},
  journal = {European J. Combin.}, volume = {6}, year = {1985}, pages = {69--78},
  doi = {10.1016/S0195-6698(85)80023-8}
}

@article{Schrijver,
  author = {Schrijver, Alexander},
  title = {Counting 1-factors in regular bipartite graphs},
  journal = {J. Combin. Theory Ser. B}, volume = {72}, number = {1},
  year = {1998}, pages = {122--135}, doi = {10.1006/jctb.1997.1798}
}

@incollection{Ryser,
  author = {Ryser, Herbert J.}, title = {Neuere {Probleme} der {Kombinatorik}},
  booktitle = {Vortr{\"a}ge {\"u}ber {Kombinatorik} ({Oberwolfach}, 1967)},
  publisher = {Mathematisches Forschungsinstitut Oberwolfach},
  year = {1967}, pages = {69--91}
}

@article{Stein,
  author = {Stein, Sherman K.}, title = {Transversals of {Latin} squares and their generalizations},
  journal = {Pacific J. Math.}, volume = {59}, year = {1975}, pages = {567--575},
  doi = {10.2140/pjm.1975.59.567}
}

@unpublished{TZ,
  author = {Tang, Yantao and Zhao, Yi},
  title = {Counting near-spanning matchings in {Latin} squares and {Steiner} triple systems},
  year = {2026},
  note = {Preprint, \href{https://arxiv.org/abs/2609.11006}{arXiv:2609.11006}},
  eprint = {2609.11006}, archivePrefix = {arXiv}
}

@unpublished{BK,
  author = {Bowtell, Candida and Keevash, Peter},
  title = {The {$n$}-queens problem}, year = {2021},
  note = {Preprint, \href{https://arxiv.org/abs/2109.08083}{arXiv:2109.08083}},
  eprint = {2109.08083}, archivePrefix = {arXiv}
}

@article{KSS,
  author = {Keevash, Peter and Sah, Ashwin and Sawhney, Mehtaab},
  title = {The existence of subspace designs},
  journal = {Proc. Lond. Math. Soc.}, volume = {131}, number = {1},
  year = {2025}, pages = {e70071}, doi = {10.1112/plms.70071},
  eprint = {2212.00870}, archivePrefix = {arXiv}
}

@article{EKOO,
  author = {Etzion, Tuvi and Kurz, Sascha and Otal, Kamil and {\"O}zbudak, Ferruh},
  title = {Subspace packings: constructions and bounds},
  journal = {Des. Codes Cryptogr.}, volume = {88}, number = {9},
  year = {2020}, pages = {1781--1810}, doi = {10.1007/s10623-020-00732-z},
  eprint = {1909.06081}, archivePrefix = {arXiv}
}

@unpublished{JJP,
  author = {Jany, Benjamin and Jurrius, Relinde and Pendavingh, Rudi},
  title = {Counting {$q$}-matroids}, year = {2026},
  note = {Preprint, \href{https://arxiv.org/abs/2606.20348}{arXiv:2606.20348}},
  eprint = {2606.20348}, archivePrefix = {arXiv}
}

@book{Oxley,
  author = {Oxley, James}, title = {Matroid Theory},
  edition = {Second}, series = {Oxford Graduate Texts in Mathematics}, volume = {21},
  publisher = {Oxford University Press}, address = {Oxford}, year = {2011},
  doi = {10.1093/acprof:oso/9780198566946.001.0001}
}

@article{DP,
  author = {Delcourt, Michelle and Postle, Luke},
  title = {Finding an almost perfect matching in a hypergraph avoiding forbidden submatchings},
  journal = {J. Lond. Math. Soc.}, volume = {114}, number = {2},
  year = {2026}, pages = {e70660}, doi = {10.1112/jlms.70660},
  eprint = {2204.08981}, archivePrefix = {arXiv}
}

@article{GL,
  author = {Glebov, Roman and Luria, Zur},
  title = {On the maximum number of {Latin} transversals},
  journal = {J. Combin. Theory Ser. A}, volume = {141},
  year = {2016}, pages = {136--146}, doi = {10.1016/j.jcta.2016.02.007},
  eprint = {1506.00983}, archivePrefix = {arXiv}
}

@article{Kwan,
  author = {Kwan, Matthew},
  title = {Almost all {Steiner} triple systems have perfect matchings},
  journal = {Proc. Lond. Math. Soc.}, volume = {121}, number = {6},
  year = {2020}, pages = {1468--1495}, doi = {10.1112/plms.12373},
  eprint = {1611.02246}, archivePrefix = {arXiv}
}

@article{GHM,
  author = {Ghandehari, Mahya and Hatami, Hamed and Mahmoodian, Ebadollah S.},
  title = {On the size of the minimum critical set of a {Latin} square},
  journal = {Discrete Math.}, volume = {293}, number = {1--3},
  year = {2005}, pages = {121--127}, doi = {10.1016/j.disc.2004.08.026},
  eprint = {math/0701015}, archivePrefix = {arXiv}
}

@article{FalconStones,
  author = {Falc{\'o}n, Ra{\'u}l M. and Stones, Rebecca J.},
  title = {Enumerating partial {Latin} rectangles},
  journal = {Electron. J. Combin.}, volume = {27}, number = {2},
  year = {2020}, pages = {Paper No. P2.47}, doi = {10.37236/9093},
  eprint = {1908.10610}, archivePrefix = {arXiv}
}

@book{vanLintWilson,
  author = {van Lint, J. H. and Wilson, R. M.},
  title = {A Course in Combinatorics}, edition = {2},
  publisher = {Cambridge University Press}, year = {2001},
  doi = {10.1017/CBO9780511987045}
}

@article{McKayWanless,
  author = {McKay, Brendan D. and Wanless, Ian M.},
  title = {On the number of {Latin} squares},
  journal = {Ann. Comb.}, volume = {9}, number = {3},
  year = {2005}, pages = {335--344},
  doi = {10.1007/s00026-005-0261-7}
}

@article{Simkin,
  author = {Simkin, Michael},
  title = {The number of $n$-queens configurations},
  journal = {Adv. Math.}, volume = {427},
  year = {2023}, pages = {109127},
  doi = {10.1016/j.aim.2023.109127}
}

@article{Knuth,
  author = {Knuth, Donald E.},
  title = {The asymptotic number of geometries},
  journal = {J. Combin. Theory Ser. A}, volume = {16}, number = {3},
  year = {1974}, pages = {398--400},
  doi = {10.1016/0097-3165(74)90063-6}
}

@article{BPvdP,
  author = {Bansal, Nikhil and Pendavingh, Rudi and van der Pol, Jorn},
  title = {On the number of matroids},
  journal = {Combinatorica}, volume = {35}, number = {3},
  year = {2015}, pages = {253--277},
  doi = {10.1007/s00493-014-3029-z}
}

@article{MNWW,
  author = {Mayhew, Dillon and Newman, Mike and Welsh, Dominic and Whittle, Geoff},
  title = {On the asymptotic proportion of connected matroids},
  journal = {European J. Combin.}, volume = {32}, number = {6},
  year = {2011}, pages = {882--890},
  doi = {10.1016/j.ejc.2011.01.016}
}

@article{WainwrightJordan,
  author = {Wainwright, Martin J. and Jordan, Michael I.},
  title = {Graphical models, exponential families, and variational inference},
  journal = {Found. Trends Mach. Learn.}, volume = {1}, number = {1--2},
  year = {2008}, pages = {1--305},
  doi = {10.1561/2200000001}
}

@article{Csikvari,
  author = {Csikv{\'a}ri, P{\'e}ter},
  title = {Lower matching conjecture, and a new proof of {Schrijver}'s and {Gurvits}'s theorems},
  journal = {J. Eur. Math. Soc.}, volume = {19}, number = {6},
  year = {2017}, pages = {1811--1844}, doi = {10.4171/JEMS/706},
  eprint = {1406.0766}, archivePrefix = {arXiv}
}

@inproceedings{Lelarge,
  author = {Lelarge, Marc},
  title = {Counting matchings in irregular bipartite graphs and random lifts},
  booktitle = {Proceedings of the Twenty-Eighth Annual {ACM--SIAM} Symposium on Discrete Algorithms},
  publisher = {SIAM}, year = {2017}, pages = {2230--2237},
  doi = {10.1137/1.9781611974782.146},
  eprint = {1507.04739}, archivePrefix = {arXiv}
}

@article{DJPRStability,
  author = {Davies, Ewan and Jenssen, Matthew and Perkins, Will and Roberts, Barnaby},
  title = {Tight bounds on the coefficients of partition functions via stability},
  journal = {J. Combin. Theory Ser. A}, volume = {160},
  year = {2018}, pages = {1--30}, doi = {10.1016/j.jcta.2018.06.005},
  eprint = {1704.07784}, archivePrefix = {arXiv}
}

@article{SpencerPacking,
  author = {Spencer, Joel},
  title = {Asymptotic packing via a branching process},
  journal = {Random Structures Algorithms}, volume = {7}, number = {2},
  year = {1995}, pages = {167--172}, doi = {10.1002/rsa.3240070206}
}

@article{RodlThoma,
  author = {R{\"o}dl, Vojt{\v{e}}ch and Thoma, Lubo{\v{s}}},
  title = {Asymptotic packing and the random greedy algorithm},
  journal = {Random Structures Algorithms}, volume = {8}, number = {3},
  year = {1996}, pages = {161--177},
  doi = {10.1002/(SICI)1098-2418(199605)8:3<161::AID-RSA1>3.0.CO;2-W}
}

@article{LinialLuriaPermutations,
  author = {Linial, Nathan and Luria, Zur},
  title = {An upper bound on the number of high-dimensional permutations},
  journal = {Combinatorica}, volume = {34}, number = {4},
  year = {2014}, pages = {471--486}, doi = {10.1007/s00493-011-2842-8},
  eprint = {1106.0649}, archivePrefix = {arXiv}
}

@article{BFL,
  author = {Bohman, Tom and Frieze, Alan and Lubetzky, Eyal},
  title = {Random triangle removal},
  journal = {Adv. Math.}, volume = {280},
  year = {2015}, pages = {379--438}, doi = {10.1016/j.aim.2015.04.015},
  eprint = {1203.4223}, archivePrefix = {arXiv}
}

@article{MadimanTetali,
  author = {Madiman, Mokshay and Tetali, Prasad},
  title = {Information inequalities for joint distributions, with interpretations and applications},
  journal = {IEEE Trans. Inform. Theory}, volume = {56}, number = {6},
  year = {2010}, pages = {2699--2713}, doi = {10.1109/TIT.2010.2046253}
}

@article{PvdPComparison,
  author = {Pendavingh, Rudi and van der Pol, Jorn},
  title = {On the number of matroids compared to the number of sparse paving matroids},
  journal = {Electron. J. Combin.}, volume = {22}, number = {2},
  year = {2015}, pages = {Paper No. P2.51}, doi = {10.37236/4899}
}

@article{BEOVW,
  author = {Braun, Michael and Etzion, Tuvi and {\"O}sterg{\aa}rd, Patric R. J. and Vardy, Alexander and Wassermann, Alfred},
  title = {Existence of {$q$}-analogs of {Steiner} systems},
  journal = {Forum Math. Pi}, volume = {4},
  year = {2016}, pages = {e7}, doi = {10.1017/fmp.2016.5},
  eprint = {1304.1462}, archivePrefix = {arXiv}
}

\makeatletter
\@setaddresses
\global\let\addresses\@empty
\makeatother


\clearpage
\appendix


\section{A counterexample to the small-codegree conjecture}\label{sec:counterexample}

The following exactly regular hypergraphs disprove Conjecture~\ref{conj:Lee} even after averaging over vertices.

\begin{theorem}\label{thm:LeeCounterexample}
    For every fixed $k\ge3$, there is a sequence of simple $n_q$-vertex $d_q$-regular $k$-uniform hypergraphs $H_q$ with $n_q=kq^3$ and
    \[
        d_q=q^3\binom{(k-1)q-1}{k-2}=\Theta_k(q^{k+1}).
    \]
    They satisfy $\Delta_2(H_q)\le d_q/q$, and every $v\in V(H_q)$ satisfies
    \[
        \Prob(v\in\calU(H_q))\le q^{-3/2}=o(d_q^{-1/k}).
    \]
    In particular, $\E[|\calU(H_q)|]=o(n_qd_q^{-1/k})$.
    For all sufficiently large $q$, these hypergraphs satisfy the hypotheses of the conjecture with $\delta=1/4$ and $\eps=1/2$.
\end{theorem}

\begin{proof}[Proof of Theorem~\ref{thm:LeeCounterexample}]
    We first record an identity for uniformly chosen matchings.
    \begin{claim}\label{clm:insertion}
        For a uniformly chosen matching $M$ of a finite hypergraph $H$,
        \[
            \E[|M|]=\E\bigl[e(H[\calU(H)])\bigr].
        \]
    \end{claim}
    \begin{claimproof}[Proof of Claim~\ref{clm:insertion}]
        Deleting a marked edge is a bijection between pairs $(M,e)$ with $e\in M$ and pairs $(J,e)$ with $e\in E(H)$ contained in $V(H)\setminus V(J)$.
        Counting both sets of pairs and dividing by $|\calM(H)|$ proves the identity.
    \end{claimproof}

    Take disjoint sets $A_1,\ldots,A_{q^2}$ of size $(k-1)q$ and a set $B$ of size $q^3$.
    Let the edges of $H_q$ be the sets $S\cup\{b\}$ with $S\in\binom{A_i}{k-1}$ for some $i$ and $b\in B$.
    Then $H_q$ is simple and has $n_q=kq^3$ vertices.
    A vertex in $A_i$ has degree
    \begin{equation}\label{eq:counterexample-degree}
        d_q=q^3\binom{(k-1)q-1}{k-2},
    \end{equation}
    while a vertex in $B$ has degree $q^2\binom{(k-1)q}{k-1}=d_q$.
    Thus $H_q$ is exactly regular.

    Two vertices in the same $A_i$ have codegree $q^3\binom{(k-1)q-2}{k-3}$, and a pair with one vertex in $A_i$ and one in $B$ has codegree $\binom{(k-1)q-1}{k-2}$.
    All other pair codegrees are zero, so
    \begin{equation}\label{eq:counterexample-codegree}
        \frac{\Delta_2(H_q)}{d_q}
        =\max\left\{\frac{k-2}{(k-1)q-1},\frac1{q^3}\right\}
        \le\frac1q.
    \end{equation}

    For a uniformly chosen matching $M$, let $(k-1)u_i$ be the number of uncovered vertices in $A_i$.
    Each $u_i$ is a nonnegative integer, since every matching edge uses $k-1$ vertices from one $A_i$.
    The number of uncovered vertices in $B$ is
    \[
        T=\sum_{i=1}^{q^2}u_i=q^3-|M|,
    \]
    and hence $|\calU(H_q)|=kT$.
    Moreover,
    \[
        e(H_q[\calU(H_q)])
        =T\sum_{i=1}^{q^2}\binom{(k-1)u_i}{k-1}
        \ge T^2.
    \]
    Indeed, partitioning the $(k-1)u_i$ uncovered vertices of $A_i$ into $(k-1)$-sets gives $\binom{(k-1)u_i}{k-1}\ge u_i$.
    Claim~\ref{clm:insertion} therefore gives
    \[
        \E[T^2]\le\E[|M|]=q^3-\E[T]\le q^3.
    \]
    By Cauchy--Schwarz,
    \begin{equation}\label{eq:counterexample-mean-vacancy}
        \E[|\calU(H_q)|]=k\E[T]\le kq^{3/2}.
    \end{equation}

    The automorphism group is transitive on $\bigcup_i A_i$ and separately on $B$.
    Every matching leaves $(k-1)T$ vertices uncovered in the first class and $T$ in the second, whose sizes are $(k-1)q^3$ and $q^3$.
    Thus~\eqref{eq:counterexample-mean-vacancy} gives, for every vertex $v$,
    \begin{equation}\label{eq:counterexample-pointwise-vacancy}
        \Prob(v\in\calU(H_q))=\frac{\E[T]}{q^3}\le q^{-3/2}.
    \end{equation}
    Since $d_q=\Theta_k(q^{k+1})$ and $k\ge3$,
    \[
        \frac{q^{-3/2}}{d_q^{-1/k}}
        =O_k(q^{-1/2+1/k})\longrightarrow0.
    \]
    Finally, for $q\ge k$,~\eqref{eq:counterexample-codegree} gives $\Delta_2(H_q)/d_q\le1/q\le n_q^{-1/4}$, and $d_q>n_q^{1/2}$ for all sufficiently large $q$.
    Exact regularity verifies the remaining degree condition in Conjecture~\ref{conj:Lee}.
    This completes the proof.
\end{proof}


\section{A conditional codegree extension and large matchings}\label{sec:conditionallifting}

We show that extending the averaged vacancy upper bound to larger pair codegrees would improve maximum-matching bounds for regular linear hypergraphs of polynomial degree.
Write $\nu(F)$ for the maximum number of edges in a matching of $F$.

\begin{lemma}[A regular lifting]\label{lem:regularlifting}
    Let $F$ be a simple $m$-vertex $D$-regular $s$-uniform hypergraph, where $s\ge2$ and $D\ge1$, and let $b\ge1$ be an integer.
    There is a simple $Db$-regular $(s+1)$-uniform hypergraph $F^+$ on $\frac{s+1}{s}mDb$ vertices satisfying
    \[
        \Delta_2(F^+)=b\Delta_2(F)
    \]
    and $\nu(F^+)\le Db\nu(F)$.
    Consequently,
    \begin{equation}\label{eq:liftingdeficit}
        1-\frac{(s+1)\nu(F^+)}{|V(F^+)|}
        \ge1-\frac{s\nu(F)}m.
    \end{equation}
\end{lemma}
\begin{proof}[Proof of Lemma~\ref{lem:regularlifting}]
    Take $Db$ vertex-disjoint copies $F^{(1)},\ldots,F^{(Db)}$ of $F$.
    For each $e\in E(F)$, introduce $b$ new vertices $x_{e,1},\ldots,x_{e,b}$, all distinct and outside these copies.
    For $i\in[Db]$ and $j\in[b]$, include the edge
    \[
        e^{(i)}\cup\{x_{e,j}\},
    \]
    where $e^{(i)}$ is the copy of $e$ in $F^{(i)}$.
    These are all edges of $F^+$, so $F^+$ is simple and $(s+1)$-uniform.
    Each old vertex has $b$ extensions of each of its $D$ incident edges, and each $x_{e,j}$ belongs to one edge from each of the $Db$ copies.
    Thus $F^+$ is $Db$-regular, and $|E(F)|=mD/s$ gives the asserted order.

    A pair of old vertices in the same copy has codegree $b$ times its codegree in $F$.
    A pair with one old and one new vertex has codegree at most one, and all other pairs have codegree zero.
    Since $\Delta_2(F)\ge1$, this proves the codegree formula.
    Finally, the edges of a matching in $F^+$ that use $F^{(i)}$ project to a matching in $F^{(i)}$.
    There are at most $\nu(F)$ such edges for each $i$, so $\nu(F^+)\le Db\nu(F)$.
    Dividing by the order formula proves~\eqref{eq:liftingdeficit}.
    This completes the proof.
\end{proof}

\begin{proposition}[A conditional consequence]\label{prop:conditionallifting}
    For fixed integers $3\le r<k$, suppose that for every fixed $\alpha>0$ there is a constant $C=C(k,\alpha)$ such that every sufficiently large simple $N$-vertex $d$-regular $k$-uniform hypergraph $H$ with $d\ge N^\alpha$ and $\Delta_2(H)\le d^{1-1/(k-2)}$ satisfies $\E[|\calU(H)|]\le CNd^{-1/k}$ for a uniformly chosen matching.
    Then, for every fixed $\eps>0$, every $n$-vertex $D$-regular linear $r$-uniform hypergraph $G$ with $D\ge n^\eps$ satisfies
    \begin{equation}\label{eq:conditionalmaximum}
        n-r\nu(G)=O_{k,r,\eps}\bigl(nD^{-1+2/k}\bigr).
    \end{equation}
\end{proposition}
\begin{proof}[Proof of Proposition~\ref{prop:conditionallifting}]
    Let $b=D^{k-3}$ and $d=Db=D^{k-2}$.
    Apply Lemma~\ref{lem:regularlifting} to $G$ with this value of $b$, and then apply it $k-r-1$ more times with $b=1$.
    The resulting simple $k$-uniform hypergraph $H$ is $d$-regular on $N=\frac kr nd^{k-r}$ vertices, and
    \begin{equation}\label{eq:lifted-codegree}
        \Delta_2(H)=D^{k-3}=d^{1-1/(k-2)}.
    \end{equation}
    Since $D\ge n^\eps$, we have
    \[
        N\le\frac kr d^{k-r+1/(\eps(k-2))}.
    \]
    Choose a fixed $\alpha>0$ with
    $\alpha(k-r+1/(\eps(k-2)))<1$.
    Then $d\ge N^\alpha$ for all sufficiently large $n$, so the assumed vacancy bound applies to $H$.
    Iterating~\eqref{eq:liftingdeficit} and using $\nu(H)\ge\E[|M_H|]$ gives
    \[
        \frac{n-r\nu(G)}n
        \le\frac{N-k\nu(H)}N
        \le\frac{\E[|\calU(H)|]}N
        \le Cd^{-1/k}=CD^{-1+2/k}.
    \]
    Multiplying by $n$ gives the deficit bound~\eqref{eq:conditionalmaximum}.
    This completes the proof.
\end{proof}

For $r=3$ and $k=5$, Proposition~\ref{prop:conditionallifting} would give a matching leaving only $O(nD^{-3/5})$ vertices uncovered in every regular linear $3$-uniform hypergraph of polynomial degree.
This would improve the $O(nD^{-1/2}(\log D)^{3/2})$ bound of Alon, Kim and Spencer~\cite{AKS} in that range.
See also~\cite{KKMO,GK} for further quantitative matching bounds.


\section{The maximum-degree upper bound}\label{sec:maxdegreeproof}

\begin{proof}[Proof of Proposition~\ref{prop:maxdegree}]
    Let $M$ be a uniformly chosen matching of $F$, write $W=V(F)\setminus V(M)$, and let $p=\E[|W|]/n$.
    We first bound the conditional entropy of $M$ given $W$.
    Lemma~\ref{lem:random-order-entropy} applies to every simple $k$-uniform hypergraph.
    Apply it to the conditional law of $M$ given $W$ and average over $W$.
    Since $F$ is linear, all coefficients $b_\pi(F)$ vanish; the tangent calculation in the proof of Lemma~\ref{lem:exposure} therefore gives~\eqref{eq:exposure-integral-loss} with $H$ replaced by $F$, $s=n(1-p)$ and $L_F(y)=1+dy^{k-1}$.
    This calculation uses only $g_v(M)\le d_{F-W}(v)\le d$, so it does not require regularity.
    Equation~\eqref{eq:pattern-log-integral} with $x=1/d$ and $q=k-1$ gives
    \[
        \int_0^1\log\left(1+\frac1{dy^{k-1}}\right)\dd y=O_k(d^{-1/(k-1)}).
    \]
    Using this estimate in the two integrals in~\eqref{eq:exposure-integral-loss} gives
    \begin{equation}\label{eq:linear-exposure}
        \HH(M\mid W)\le\frac{n(1-p)}k(\log d-k+1)
        -\frac{\E[S_W]}{kd}+O_k(nd^{-1/(k-1)}),
    \end{equation}
    where $S_W=\sum_{v\notin W}(d-d_{F-W}(v))$.
    Indeed, the first integral is $\log d-k+1+O_k(d^{-1/(k-1)})$, the second is at least $d^{-1}(1-O_k(d^{-1/(k-1)}))$, and $0\le S_W\le nd$.

    Write
    \[
        R=\log|\calM(F)|-\frac nk(\log d-k+1).
    \]
    Since $M$ is uniform and $W$ is determined by $M$, the entropy chain rule gives $\log|\calM(F)|=\HH(W)+\HH(M\mid W)$.
    Substituting~\eqref{eq:linear-exposure} and subtracting $\frac nk(\log d-k+1)$ gives
    \begin{equation}\label{eq:maxdegree-before-degree-loss}
        R\le\HH(W)-\frac{\E[S_W]}{kd}
        -\frac{np}{k}(\log d-k+1)+O_k(nd^{-1/(k-1)}).
    \end{equation}

    The next claim bounds $\HH(W)-\E[S_W]/(kd)$ in~\eqref{eq:maxdegree-before-degree-loss}.
    Its term proportional to $np$ cancels the corresponding term in that inequality.

    \begin{claim}\label{clm:maxdegree-loss}
        We have
        \[
            \HH(W)-\frac{\E[S_W]}{kd}
            \le\frac{np}{k}(\log d-k+1)+nd^{-1/k}+O_k(nd^{-2/k}).
        \]
    \end{claim}
    \begin{claimproof}[Proof of Claim~\ref{clm:maxdegree-loss}]
        Write $a_v=d-d_F(v)$ and $T=\sum_v a_v$.
        Counting the original deficits at covered vertices and the incidences lost through edges meeting $W$ gives
        \[
            \begin{aligned}
                S_W
                &=\sum_{v\notin W}a_v
                +\sum_{e\in E(F)}(k-|e\cap W|)\mathbf1_{\{e\cap W\ne\varnothing\}}\\
                &=T+k\sum_{e\in E(F)}\mathbf1_{\{e\cap W\ne\varnothing\}}-d|W|.
            \end{aligned}
        \]
        The last equality uses $a_v+d_F(v)=d$.
        Add $a_v$ copies of $\{v\}$ to the edge family, so every coordinate belongs to exactly $d$ members.
        Write $p_v=\Prob(v\in W)$.
        Taking expectations in the identity for $S_W$ and applying Lemma~\ref{lem:shearer} gives
        \begin{equation}\label{eq:maxdegree-shearer-bound}
            \begin{aligned}
                \HH(W)-\frac{\E[S_W]}{kd}
                &\le\frac1d\sum_{e\in E(F)}
                \bigl(\HH(W\cap e)-\Prob(W\cap e\ne\varnothing)\bigr)\\
                &\quad+\frac1d\sum_v a_vh(p_v)+\frac{np}{k}-\frac{T}{kd}.
            \end{aligned}
        \end{equation}

        On each edge, use the weights from the proof of Claim~\ref{clm:vacancy-degree-loss} at $\lambda=1$.
        The resulting estimate~\eqref{eq:edge-entropy-hit} holds for every law of $W\cap e$.
        On a singleton, Lemma~\ref{lem:entropy-counting} with weights $w(0)=1$ and $w(1)=\mathrm e d^{-1/k}$ gives
        \[
            h(p_v)\le-\left(1-\frac{\log d}{k}\right)p_v+\log(1+\mathrm e d^{-1/k}).
        \]
        To sum the edge and singleton estimates, use $k|E(F)|=nd-T$ and
        \[
            \sum_{e\in E(F)}\E[|e\cap W|]+\sum_v a_vp_v=dnp.
        \]
        Substituting into~\eqref{eq:maxdegree-shearer-bound} gives
        \[
            \begin{aligned}
                \HH(W)-\frac{\E[S_W]}{kd}
                &\le-np\left(1-\frac{\log d}{k}\right)+\frac{nd-T}{d}d^{-1/k}\\
                &\quad+\frac Td\log(1+\mathrm e d^{-1/k})+\frac{np}{k}-\frac{T}{kd}+O_k(nd^{-2/k})\\
                &=\frac{np}{k}(\log d-k+1)+nd^{-1/k}\\
                &\quad+\frac Td\left(\log(1+\mathrm e d^{-1/k})-d^{-1/k}-\frac1k\right)+O_k(nd^{-2/k}).
            \end{aligned}
        \]
        The coefficient in parentheses tends to $-1/k$ as $d\to\infty$.
        Since $T\ge0$, its contribution is nonpositive for sufficiently large $d$, proving the claim.
    \end{claimproof}

    Applying Claim~\ref{clm:maxdegree-loss} to~\eqref{eq:maxdegree-before-degree-loss} cancels the terms proportional to $np$ and gives
    \[
        R\le nd^{-1/k}+O_k(nd^{-2/k}+nd^{-1/(k-1)})
        =nd^{-1/k}+O_k(nd^{-1/(k-1)}),
    \]
    since $2/k>1/(k-1)$.
    Adding back $\frac nk(\log d-k+1)$ proves the asserted count.
    This completes the proof.
\end{proof}


\end{document}